\documentclass{article}

\usepackage[margin=1in]{geometry}
\usepackage{latexsym}
\usepackage{indentfirst}
\usepackage{graphicx}
\usepackage{amsmath,amsfonts}
\usepackage{amssymb}
\usepackage{amsthm}
\usepackage{tikz}
\usepackage{bm,bbm}
\usepackage{accents}
\usepackage{changepage}
\usepackage{color}
\usepackage{bigints}
\usepackage{enumerate}
\usepackage{amsmath,nccmath}
\usepackage{yhmath}
\usepackage{dirtytalk}
\usepackage{comment}
\usepackage{enumitem}
\usepackage{hyperref}
\usepackage{cleveref}
\usepackage{comment}
\usepackage{caption}
\usepackage{subcaption}
\usepackage{float}
\usepackage{svg}

\crefname{subsection}{subsection}{subsections}

\newtheorem{prop}{Proposition}[section]
\newtheorem*{defi*}{Definition}
\newtheorem{defi}[prop]{Definition}
\newtheorem{exam}[prop]{Example}
\newtheorem{lem}[prop]{Lemma}
\newtheorem{rem}[prop]{Remark}
\newtheorem{thm}[prop]{Theorem}
\newtheorem{coro}[prop]{Corollary}
\newtheorem*{main}{Theorem}

\def\XXint#1#2#3{{\setbox0=\hbox{$#1{#2#3}{\int}$ }
\vcenter{\hbox{$#2#3$ }}\kern-.6\wd0}}

\def\eps{{\varepsilon}}

\newtheorem{theorem} {\sc  Theorem\rm} [section]

\newtheorem{definition}[theorem]{\sc  Definition\rm}

\newcounter{marnote}

\DeclareFontFamily{OT1}{rsfs}{}
\DeclareFontShape{OT1}{rsfs}{m}{n}{ <-7> rsfs5 <7-10> rsfs7 <10-> rsfs10}{}
\DeclareMathAlphabet{\mycal}{OT1}{rsfs}{m}{n}

\def\supp{{\rm Supp}\,}

\def\be{\begin{equation}}
\def\ee{\end{equation}}

\def\supp{{\rm Supp}\,}

\newcommand{\R}{\mathbb{R}}

\newcommand{\po}{{\pi_{\partial D}}}
\newcommand{\Oe}{{D_{\epsilon_0}}}

\def\be{\begin{equation}}
\def\ee{\end{equation}}
\def\bea#1\eea{\begin{align}#1\end{align}}
\def\non{\nonumber}

\numberwithin{equation}{section}
\title {From Monte Carlo to neural networks approximations\\ of boundary value problems

\vspace{0.5 cm}
}

\author{
Lucian Beznea$^{2,1}$ \thanks{E-mails: \texttt{lucian.beznea@imar.ro}}
\and 
Iulian C\^{i}mpean$^{1,2}$\thanks{E-mails: \texttt{iulian.cimpean@unibuc.ro;iulian.cimpean@imar.ro}}
\and
Oana Lupa\c scu-Stamate $^{3}$\thanks{E-mails: \texttt{oana.lupascu$\_$stamate@yahoo.com}}
\and
Ionel Popescu $^{1,2}$\thanks{E-mails: \texttt{ioionel@gmail.com,ionel.popescu@fmi.unibuc.ro}}
\and
Arghir Zarnescu $^{4,2,5}$\thanks{E-mails: \texttt{azarnescu@bcamath.org}}
\vspace{0.5 cm}
}

\date{\small
$^1$University of Bucharest, Faculty of Mathematics and Computer Science,\\ 
14 Academiei str., 70109, Bucharest, Romania\\
$^2$ Simion Stoilow Institute of Mathematics of the Romanian Academy,\\ P.O. Box 1-764, RO-014700 Bucharest, Romania \\
$^3$``Gheorghe Mihoc~--~Caius Iacob" Institute of Mathematical Statistics and Applied Mathematics\\ of the Romanian Academy, 13~Calea 13 Septembrie, 050711~Bucharest, Romania\\
$^4$ BCAM, Basque Center for Applied Mathematics, Mazarredo 14, E48009 Bilbao, Bizkaia, Spain\\
$^5$ IKERBASQUE, Basque Foundation for Science, Maria Diaz de Haro 3, 48013, Bilbao, Bizkaia, Spain
}

\begin{document}

\maketitle

\begin{abstract}
In this paper we study probabilistic and neural network approximations for solutions to Poisson equation subject to H\" older data in general bounded domains of $\mathbb{R}^d$.
We aim at two fundamental 
goals.

The first, and the most important, we show that the solution to Poisson equation can be numerically approximated in the sup-norm by Monte Carlo methods, { and that this can be done highly efficiently if we use 
 a modified version} of the walk on spheres algorithm { as an acceleration method. This provides estimates which are efficient with respect to the prescribed approximation error and with polynomial complexity in the dimension and the reciprocal of the error.}   {A crucial feature is that} the overall number of samples does not not depend on the point at which the approximation is performed. 

As a second goal, we show that the obtained Monte Carlo solver renders { in a constructive way} ReLU deep neural network (DNN) solutions to Poisson problem, whose sizes depend at most polynomialy in the dimension $d$ and in the desired error.  In fact we show that the random DNN provides with high probability a small approximation error and low polynomial complexity in the dimension.  

\medskip


\vspace{2mm}

\end{abstract}

\newpage

\tableofcontents

\section{Introduction}
Partial differential equations provide the most commonly used framework for modelling a large variety of phenomena in science and technology. Using these models in practice requires fast, accurate and stable computations of solutions of PDEs. Broadly speaking there exist two large classes of simulations: deterministic and stochastic. The deterministic methods (e.g. finite differences, finite element methods, etc) are very effective in globally approximating the solutions but their computational effort grows exponentially with respect to the dimension of the space. On the other hand, the probabilistic methods manage to overcome the dimensionality issue, but they are usually employed to obtain approximations at a given point and changing the point requires a different approximation. 

\smallskip
In recent years, another powerful class of methods has been developed, namely the (deep) neural network models (in short DNN).  These have been able to provide a remarkable number of achievements in the technological realm, such as: image classification, language processing and time series analysis, to name only a few. However, despite their remarkable achievements their rigorous understanding is still in its infancy.

On the theoretical side it is known that DNNs are capable of providing good approximation properties for continuous functions \cite{cybenko1989approximation,EDGB_2019,pinkus1999approximation,yarotsky2017error}. For a more recent in depth analysis see \cite{EDGB_2019} and the references in there. However it should be noted that these approximations are not constructive and indeed the issue of constructibility and error estimates is the crucial one from the practical point of view. 

In what concerns the approximation of solutions of PDEs there are numerical treatments in low dimensions, as for example in \cite{lagaris1998artificial,lagaris2000neural,malek2006numerical},  which propose schemes for solving PDEs using some form of neural networks. However, these do not provide any error estimates. Their approaches depend on a grid discretization of the space and for the most challenging case, namely that of simulations in high dimensions, there are no theoretical guarantees that the methods would scale well in high dimensions. Typically, these approximations, though convincing, remain at the level of numerical experiments.

On the other hand there is a rigorous body of literature which treats the approximation of the solution where error estimates are provided for PDEs with neural networks as, for example, in
\cite{raissi2019physics,SIRIGNANO20181339,malek2006numerical,NEURIPS2021_7edccc66,GrHe21,han2018solving} but typically these scale poorly in high dimensions. A discussion of some of these is provided in the next subsection.

Our aim in this paper will be to study approximations of Poisson equation and provide DNNs built on stochastic approaches that address some of these shortcomings mentioned before,  significantly advancing the state of the art and providing tools that can be extended to more general equations, as the cost of a suitable increase in the technical details. A discussion of our contribution is provided in subsection ~\ref{ss:contribution}.

\medskip

\subsection{{ Brief technical description of the main results.}} 
Our starting point is suggested by \cite{GrHe21} which essentially builds on the stochastic representation of the solution to the Poisson equation.  In turn, the stochastic representation is then followed by the standard \emph{walk on spheres} (WoS) method to accelerate the computation of the integrals of the Brownian trajectory.  This is then used to construct neural networks approximations.  In the present paper we fundamentally expand, clarify, simplify and explicitly construct starting from some of the ideas pointed out in \cite{GrHe21}.  { An extension of \cite{GrHe21} to the fractional Laplacian has been developed in \cite{valenzuela2022new}, so the refined methods proposed here should essentially apply to \cite{valenzuela2022new} as well.}

Throughout, $D\subset \mathbb{R}^d$ is at least a bounded open set; we write $h\in C(D)$ to say that $h:D\rightarrow \mathbb{R}$ is continuous on $D$.
$L^p(D)$ is the standard Lebesgue space with norm denoted by $|\cdot|_{L^p(D)}$.
For a bounded function $h:D\rightarrow \mathbb{R}$, that is for $h\in L^\infty(D)$, we shall denote by $|h|_\infty$ 
the essential sup-norm of $h$.
For $\alpha \in [0,1]$ and $h:D\rightarrow \mathbb{R}$ an $\alpha$-H\"older function for $\alpha\in (0,1)$( or Lipschitz, for $\alpha=1$), we set $|h|_\alpha:=\sup\limits_{x,y\in D}\dfrac{|h(x)-h(y)|}{|x-y|^\alpha}$.

Now we descend into the description of our main results. We study the Poisson boundary value problem 
\begin{equation}\label{e:0}
\begin{cases}
\frac{1}{2}\Delta u=-f \,\textrm{ in } D  \\
u|_{\partial D}=g, 
\end{cases}
\end{equation}
where $D$ is a bounded and open set in $\mathbb{R}^d$, whilst $f:D\to\R$ and $ g:\overline{D}\to\R$ are given continuous functions. { It is well known that there exists a unique solution $u\in C(D)\cap H^1_{loc}(D)$ to \eqref{e:0}, see \cite[Theorem 6]{Ge92} or \Cref{thm:representation} from below. 
The fact that $u\in C(D)\cap H^1_{loc}(D)$ is a solution to \eqref{e:0} means that
\begin{align}
    &\frac{1}{2}\int\limits_D \langle \nabla u, \nabla \varphi\rangle\;dx = \int\limits_D f\varphi \;dx \quad \forall \varphi\in C_c^\infty(D)\\
    &\lim\limits_{D\ni x\to x_0\in \partial D}u(x)=g(x_0), \quad \forall x_0\in \partial D \mbox{ which is a regular point for } D. 
\end{align}
}
In \cite{GrHe21} the domain $D$ is taken to be convex.  We treat several layers of generality for the domain $D$ which lead to different final results.

All the random variables employed in the sequel are assumed to be on the same probability space $(\Omega,\mathcal{F},\mathbb{P})$, whilst the expectation is further denoted by $\mathbb{E}$.
Further, let us consider the generalised walk on sphere process defined as 
\begin{equation}\label{e:2} 
{X}^{x}_0:=x\in D, \quad 
{X}^{x}_{n+1}:={X}^{x}_n+\widetilde{r}({X}^{x}_n)U_{n+1},\; n\geq 0, 
\end{equation}
where $x$ is the starting point in the domain $D$, the function $\widetilde{r}$ denotes the replacement of the { distance function $r$ to the boundary $\partial D$,} and $U_i$ are drawn independent and identically on the unit sphere in $\R^d$.  
Essentially we use a Lipschitz $\widetilde{r}$ such that $ \widetilde{r}\le r$ on the whole domain and $\beta r \le \widetilde{r}$ as long as $r\ge\varepsilon_0>0$.   
We call such a candidate a $(\beta,\varepsilon_0)$-distance and these constants play an important role in the estimates below.  In all cases we can take  $\beta\ge1/3$ as we point out in Remark~\ref{rem:adistance} below.
With this process at hand,  { we introduce the Monte Carlo estimator $u_M^N$ of the solution $u$ to problem \eqref{e:0} by} 
\begin{equation}\label{e:5}
    u_{M}^{N}(x):=\frac{1}{N}\sum_{i=1}^N\left[ g(X^{x,i}_{M})+\frac{1}{d}\sum\limits_{k= 1}^M \widetilde{r}^2(X^{x,i}_{k-1})f\left(X^{x,i}_{k-1}+\widetilde{r}(X^{x,i}_{k-1})Y^i\right)\right], \quad x\in D, M,N\geq 1. 
\end{equation}
Here, the random variable $ U_{n,i}$ and $Y^i$, ${n,i\ge 0}$, are all independent, $U_{n,i}$ is drawn uniformly on the unit sphere, $X_n^{x,i}$ is given by \eqref{e:2} with $U_n$ replaced by $U_{n,i}$, whilst $Y^i$ is drawn on the unit ball in $\R^d$ from the distribution { $\mu$} which has an explicit density proportional to $|y|^{2-d}-1, |y|<1$ if $d\geq 3$, and which is in fact the (normalized) Green kernel of the Laplacian on the unit ball with pole at $0$. { It is easy to see that if $R$ and $Z$ are independent random variables such that $R$ has distribution $Beta\left(\frac{4-d}{d-2},2\right)$ on $[0,1]$ and $Z$ is uniformly distributed on the unit sphere in $\mathbb{R}^d$, then $R^{\frac{1}{d-2}}Z$ has distribution $\mu$.}

{ The first part of the main result of this paper is the following:

\begin{main}[Part I; see \Cref{thm:main} for the full quantitative version]
Fix a small $\varepsilon_0>0$, $\beta\in (0,1]$, $\widetilde{r}$ a $(\beta,\varepsilon_0)$-distance, and consider $u_M^N$ given by \eqref{eq:MCestimator}.
Also, assume that $f$ and $g$ are $\alpha$-H\" older on $D$ for some $\alpha\in (0,1]$.
Then, for any compact subset $F\subset D$, for all $N, M, K \geq 1$, $\gamma>0$ and $\varepsilon \geq \varepsilon_0$, there are explicit quantities $A(F,M,K,d, 
\varepsilon)$ and $B(M,K,d)$ given in terms of the boundary regularity, the parameters $N,M,K, d,\varepsilon$, the set $F$ and the data $f,g$
such that 
\begin{equation}\label{eq:main:e:i}
\mathbb{E}\left\{\sup_{x\in F} \left | u(x)-u_M^{N}(x)\right|\right\}\le A(F,M,K,d,\varepsilon)+\frac{B(M,K,d)}{\sqrt{N}}.
\end{equation}
Moreover, for an arbitrary bounded and open set $D$, for any compact subset $F\subset D$, we have that 
\begin{equation}\label{e:main:f:i}
\lim_{\varepsilon_0\to 0}\lim_{M\to \infty}\lim_{N\to\infty}\mathbb{E}\left\{\sup_{x\in F} \left | u(x)-u_M^{N}(x)\right|\right\}=0.
\end{equation}

In addition, if the domain $D$ satisfies the uniform ball condition, we can take $F=D$ in the estimates \eqref{eq:main:e:i} and \eqref{e:main:f:i}.  
\end{main}
For completeness, we recall that a bounded and open set $D\subset \mathbb
{R}^d$ is said to satisfy a uniform exterior ball condition if there exists a radius $r>0$ such that for every $y\in \partial D$ there exists $a\in \mathbb{R}^d$ for which $B(a,r)\subset \mathbb{R}^d\setminus D$ and $y\in \partial D\cap \partial B(a,r)$.
}

{
\begin{rem} A few comments are in place here.  
\begin{enumerate}[label=(\roman*)]
\item The estimates are true for any arbitrary domain (bounded and open), both in expectation and also in the tail.   However, in this very general case we do not get any quantitative estimates, only asymptotic convergence guarantees on compact subsets. 

\item  The full power of the result is a little more technical and states that we have a  tail estimate in the form
\begin{equation} \label{eq:main1:i}
\mathbb{P}\left( \sup_{x\in F} \left | u(x)-u_M^{N}(x)\right| \geq \gamma \right) \leq 2 \exp\left(C_1(M,K,d)-\frac{\left((\gamma-A(F,M,K,d,\varepsilon))^+\right)^2}{C_2(M,d)}N \right),
\end{equation}
where 
\begin{equation}\label{e:A,B,C}
\begin{split}
C_1(M,K,d):&=d\left(\lceil M/\alpha\rceil\log(2+|\widetilde{r}|_1)+\log(K)\right),\quad
C_2(M,d):=|g|_\infty+\frac{M}{d}{\sf diam}(D)^2|f|_\infty \\
A(F,M,K,d,\varepsilon)&:=2\left(|g|_\alpha+\frac{{\rm diam}(D)^2|f|_\alpha+2{\rm diam}(D)|f|_\infty}{d}\right)\left(\frac{{\sf diam}(D)}{K}\right)^\alpha \\
&\quad + d^{\alpha/2}|g|_\alpha v^{\alpha/2}(F,\varepsilon)+|f|_\infty v(F,\varepsilon) + (4|g|_\infty+\frac{2}{d}{\sf diam}(D)^2|f|_\infty) e^{-\frac{\beta^2\varepsilon^2}{4{\sf diam}(D)^2}M}.
\end{split}
\end{equation}

\item The parameter  $\varepsilon$ measures the closeness to the boundary and can be taken arbitrarily such that $\varepsilon\geq \varepsilon_0$.   The function $v(x,\varepsilon)$ and $v(F,\varepsilon)$ defined in \eqref{def:wepsilon} measure the geometry of the boundary from a rather stochastic viewpoint.  We can upper-bound $v(F,\varepsilon)$ by a more more tractable and analytical version of it, namely $|v|_\infty(\varepsilon)$ (see \eqref{def:wepsilon}).   Moreover, if the domain $D$ satisfies the exterior ball condition, we can replace the compact set $F$ with the whole domain $D$ and $v(D,\varepsilon)$ with $\varepsilon {\sf adiam}(D)$.  

Furthermore, for the case of $\delta$-defective convex domains $D$ (see \eqref{eq:dconvex}) we can replace $e^{-\frac{\beta^2\varepsilon^2}{4{\sf diam}(D)^2}M}$ by $\left(1-\frac{\beta^2(1-\delta)}{4d}\right)^M \sqrt{\frac{{\sf diam}(D)}{\varepsilon}}$.  
  
\item There are many parameters in  \eqref{eq:main:e} and \eqref{eq:main1:i}. $M$ stands for the number of steps the walk on spheres is allowed to take.  $N$ is the number of Monte Carlo simulations.

The mysterious constant $K$ comes from a grid discretisation used for the estimate.  Notice that the left hand side of \eqref{eq:main:e:i} or \eqref{eq:main1:i} does not depend on $\varepsilon$ or $K$.  Furthermore,  the larger the $K$ and $M$ the larger $C_1(M,K,d)$, nevertheless this is compensated by $N$ and the dependency of $A(F,M,K,d,\varepsilon)$, which becomes smaller for large $K$ and $M$. Therefore, the strategy is to optimize the right hand sides of  \eqref{eq:main:e:i} or \eqref{eq:main1:i} to obtain the best estimate.  

\item We can take the constant 
\[
B(M,K,d):=C_2(M,d)(\sqrt{C_1(M,K,d)+\log(2)}+1).
\]
which, as well as the constants $C_1(M,K,d)$ and $C_2(M,d)$, does not depend on $N$ nor $\varepsilon$. Thus the limit over $N$ in \eqref{e:main:f:i} and then over $M$ leaves the right hand side of \eqref{eq:main:e:i} dependent only on $K$ and $\varepsilon$.  However, the limit in \eqref{e:main:f:i} is independent of $K$ and $\varepsilon$.  Letting $K\to\infty$ and $ \varepsilon_0\leq \varepsilon\to0$ yields \eqref{e:main:f:i} for general domains. Quantitative versions can be obtained by carefully tuning all the parameters, $N$, $M$, $K$, $\varepsilon$.

\item The estimate \eqref{eq:main1:i} is the key to guaranteeing that the error is actually small with high probability.  This is more useful than the expectation result \eqref{e:main:f:i}.

\item Finally, we point out that these rather intricate terms are the key in choosing the dependency of all the constants $N,M,K,\varepsilon$ in terms of dimension $d$ for large $d$.  Nevertheless, the estimates  are very useful also for small dimensions.  

\item When $g\in C^2(\partial D)$ then the above estimate can be improved; for more details we refer the reader to the extended arxiv version of this paper.
\end{enumerate}
\end{rem}
}

After all the remarks above, the goal is to make the right hand side of \eqref{eq:main1:i} small.  Assuming that the domain $D$ { satisfies the uniform exterior ball condition}, and both $f,g$ have $\alpha$-H\"older regularity, { then we can ensure that
\begin{equation}
\mathbb{P}\left( \sup_{x\in D} \left | u(x)-u_M^{N}(x)\right| \geq \gamma \right)\leq \eta,
\end{equation}
by sampling $N$ times a number of $M$ steps of the WoS chain \eqref{e:2}, with complexity (see more details below in Remark~\ref{r:choice})
\begin{equation*}
M=  \mathcal{O}\left(\frac{d^2\log(1/\gamma)}{\gamma^{4/\alpha}}\right), \quad 
N= \mathcal{O}\left(\frac{d^2\log(1/\gamma)^2\left[d^2\log(1/\gamma)+\gamma^{2/\alpha}\log(1/\eta)\right]}{\gamma^{2+4/\alpha}}\right).
\end{equation*}
Moreover, if $D$ is defective convex (see \eqref{eq:dconvex}) then we get the significant improvement
\begin{equation*}
M= \mathcal{O}\left(d\log(d/\gamma)\right),  \quad 
N= \mathcal{O}\left(\frac{\log^2(d/\gamma)\left(d^2\log(d/\gamma)+\log(1/\eta)\right)}{\gamma^2}\right).
\end{equation*}
}
The expressions inside the $\mathcal{O}$ symbols represent the dependency on $f,g, {\sf diam}(D)$ and the geometry of $D$.  As these expressions make it clear, the dependency on the dimension and $\gamma$ is poly-logarithmic, with better choices in the case the domain $D$ has better geometry. 
We conclude by pointing out that the total number of flops to compute the approximation $u_M^N$ is $\mathcal{O}(d^2MN)$, thus in total a polynomial complexity in $d$.

{
Before we proceed to the DNN theoretical counterpart, let us briefly present a single numerical test meant to support the estimates in {\bf Theorem}(Part I). For more numerical tests in this sense we refer the reader to \Cref{S:numerics}. Here, \Cref{fig:Intro} below depicts the evolution w.r.t. $N$ of the mean errors $\mathbb{E}\left\{|D|^{-1}\int_D \left | u-u_M^{N}\right| \; dx\right\}$  and $\mathbb{E}\left\{\sup_{x\in D} \left | u(x)-u_M^{N}(x)\right|\right\}$ for the Poisson problem \eqref{e:0}, where $d=100$, $D$ is the coronal cube $D_{\sf ac}$ defined in \Cref{S:numerics}, whilst the exact solution $u$ is given by \eqref{eq:exact_u} from \Cref{S:numerics}. 
The above expectations have also been approximated by Monte Carlo simulations, for more details see \Cref{S:numerics}, more precisely \Cref{fig:N_steps_d=100_ac}.

\begin{figure}[H]

\begin{subfigure}{0.5\textwidth}
\centering
\includegraphics[scale=0.4]{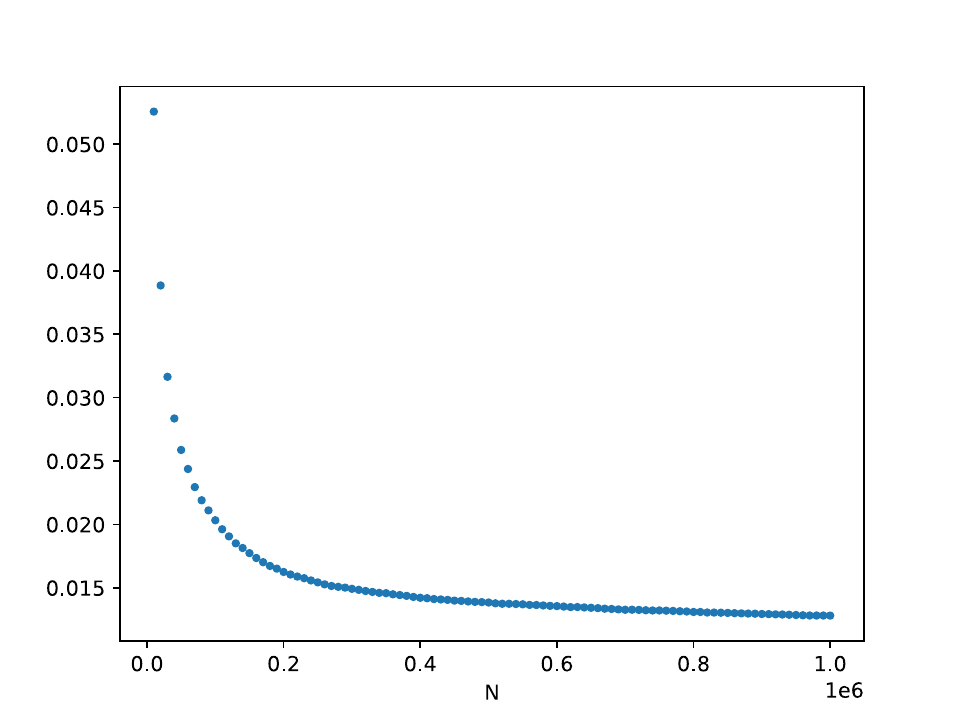}
\caption{$\mathbb{E}\left\{1/|D|\int_D \left | u(x)-u_M^{N}(x)\right| \; dx\right\}$}
\label{fig:Intro_a}
\end{subfigure}
\begin{subfigure}{0.5\textwidth}
\centering
\includegraphics[scale=0.4]{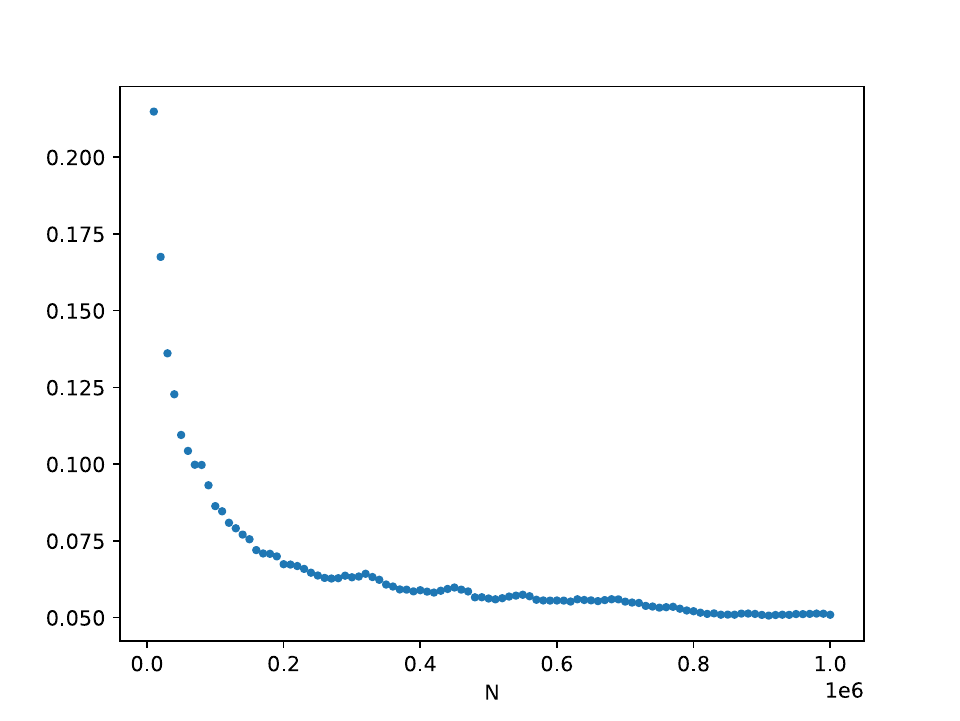}
\caption{$\mathbb{E}\left\{\sup_{x\in D} \left | u(x)-u_M^{N}(x)\right|\right\}$}
\label{fig:Intro_b}
\end{subfigure}
\caption{ The evolution of (Monte Carlo estimates of) $\mathbb{E}\left\{1/|D|\int_D \left | u(x)-u_M^{N}(x)\right| \; dx\right\}$ and $\mathbb{E}\left\{\sup_{x\in D} \left | u(x)-u_M^{N}(x)\right|\right\}$ w.r.t. $N$, for $D=D_{\sf ac}$ (see \Cref{S:numerics}), $d=100$, $M=500$. The computed errors are decreasing to a small value as the number of WoS trajectories $N$ increases. The limit error attained when $N$ goes to infinity is not zero as it depends on $M$, but it decreases to zero as the latter parameter is also increased to infinity; see \Cref{S:numerics} for more details.}
\label{fig:Intro}
\end{figure}
}

{ The second part of the main result of this paper} is dedicated to the construction of a neural network approximation. In fact, in this paper we deal only with the special class of feed-forward neural networks whose activation function is given by ReLU (rectified linear unit), and such a network shall further be referred to as a ReLU DNN.
The essential part of this construction is the formula \eqref{e:5}. 
The key fact is that we choose $\widetilde{r}$ to be already a ReLU DNN. 
This is the main reason for our modification of the walk on spheres algorithm with $\widetilde{r}$ instead of the usual distance to the boundary.
In addition, having some ReLU DNN approximations of the data $f$ respectively $g$, we can use these building blocks together with some basic facts about the ReLU DNNs in conjunction with \eqref{eq:main1:i} to get the following result.  
{
\begin{main}[Part II; see Theorem~\ref{thm:mainNN} for a more explicit version]\label{thm:part II}
Under the same context as {\rm {\bf Theorem} (Part I)}, assume that $D$ satisfies the uniform exterior ball condition, and that we are given ReLU DNNs $\phi_f:D\rightarrow\mathbb{R},\phi_g:\overline{D}\rightarrow\mathbb{R},\phi_r:D\rightarrow \mathbb{R}$ such that
\begin{equation*}
    |f-\phi_f|_\infty\leq \epsilon_f\leq |f|_\infty, \quad |g-\phi_g|_\infty\leq \epsilon_g, \quad  |r-\phi_r|_\infty\leq \epsilon_r.
\end{equation*} 
If $\gamma>0$ is chosen such that $a.1)-a.3)$ from \Cref{thm:mainNN} hold, and $0<\eta<1$, then we can construct an explicit (see \Cref{thm:mainNN} for the precise construction) (random) ReLU DNN $\mathbb{U}(x)$ such that 
\begin{equation*}
\mathbb{P}\left(\sup\limits_{x\in D} \left| u(x)-\mathbb{U}(\cdot,x)\right| \leq \gamma \right) \geq 1-\eta,
\end{equation*}
with
\begin{equation*}
{\rm size}(\mathbb{U}(\omega,\cdot))= \mathcal{O}\left(d^7\gamma^{-16/\alpha-4}\log^4\left(\frac{1}{\gamma}\right)\left[d^3\gamma^{-4/\alpha}\log\left(\frac{1}{\gamma}\right)+\log\left(\frac{1}{\eta}\right)\right]\log(2+|\phi_r|_1){\rm S}\right),
\end{equation*}
where 
\[
{\rm S}:=\left[\max(d,\mathcal{W}(\phi_r),\mathcal{L}(\phi_r))+{\rm size}(\phi_r)+{\rm size}(\phi_g)+{\rm size}(\phi_f)\right].
\]
Furthermore, if $D$ is defective convex (see \eqref{eq:dconvex}) then we get the significant improvement
\begin{equation*}
{\rm size}(\mathbb{U}(\omega,\cdot))=\mathcal{O}\left(\frac{d^3}{\gamma^2}\log^4\left(\frac{d}{\gamma}\right)\left[d^2\log\left(\frac{d}{\gamma}\right)+\log\left(\frac{1}{\eta}\right)\right]\log(2+|\phi_r|_1){\rm S}\right).
\end{equation*}
Here $\mathrm{size}$ denotes the number of non-zero parameters in the neural network, $\mathcal{W}(\phi)$ and $
\mathcal{L}(\phi)$ represent the width respectively the length of the neural network $\phi$.  
The implicit constants depend on $|g|_\alpha,|g|_\infty,|f|_\infty,{\rm diam}(D),\\ 
{\rm adiam}(D), \delta,\alpha$. 
\end{main}  
}    

\subsection{Further extensions}

This paper has outlined a novel method with the potential to extend to a general principle for solving high-dimensional partial differential equations (PDEs) without encountering the curse of dimensionality. The approach is threefold:

\begin{enumerate}[label=(\roman*)]
    \item Initially, we establish a probabilistic representation of the solution via a stochastic process. Notably, such processes are already known to exist for a wide range of PDEs, ranging from linear types \cite{LeStTr} to fully non-linear models \cite{bsde, Gall, beznea, pardoux1999forward, gobet2016monte}.
    \item The subsequent phase involves the rapid simulation of this stochastic process, aligning closely with the PDE solution representation. This step is harmonized with the Monte-Carlo method for approximating the solution.
    \item The final phase focuses on approximating the solution at specific points within an auxiliary grid of the domain. Following this, we employ a concentration inequality to extend this approximation to the continuous domain, thus compensating for the complexity introduced by the grid.
\end{enumerate}

Such a strategy finds immediate application in the realm of parabolic PDEs solvable by diffusion processes, as discussed in \cite{BaTa} and \cite{BuChSy}, especially within time-dependent domains. It also applies to PDEs characterized by (sticky-)reflecting boundary conditions, in line with the frameworks established in \cite{Hsu} and \cite{Max}.

Furthermore, this methodology opens new avenues in molecular dynamics, particularly in scenarios where the operator exhibits discontinuous coefficients, as explored in \cite{talay2, talay3, lejay2013new}. Another promising direction involves the extension of this approach to (sub-)Riemannian or metric measure spaces, utilizing tools as the ones in \cite{Stroock, Hsu2, von2005transport, arnli, Baudoin, thal}.

\subsection{\bf On related previous work.} For the sake of comparison (see \Cref{ss:contribution}), let us present here a short review of existing works that we find strongly connected to our paper, trying to point out some limitations that are fundamentally addressed in our present work. 
For a comprehensive overview of deep learning methods for PDEs we refer to \cite{BeHuJeKu23}.

\paragraph{Monte Carlo methods for PDEs} The Monte Carlo method for solving linear PDEs has been well understood and intensively used for a long time. Let us mention e.g. \cite{TaTu90} for the case of linear parabolic PDEs on $\mathbb{R}^d$, and \cite{SaTa95} for linear elliptic PDEs in bounded domains; one issue that is encountered in all these classical works and which is particularly crucial to us, is that the theoretical errors are point-dependent, in the sense that there is no guarantee that one can use the same Monte Carlo samples uniformly for all the locations in a bounded and open Euclidian set where the solution is aimed to be approximated.
Recently, it was shown in \cite{JuJeKrAnWu20} and \cite{HuJeKu19} (see also \cite[Theorem 1.1]{HuJeWu20}) that a multilevel Picard Monte Carlo method can be derived in order to numerically approximate the solutions to semilinear parabolic PDEs without suffering from the curse of high dimensions; here as well, the obtained probabilistic errors are derived for a given fixed location where the exact solution is approximated.
Monte Carlo methods have also been extended to fully nonlinear parabolic PDEs in $\mathbb{R}^d$, as for example in \cite{FaToWa11}, where a mixture between the Monte Carlo method and the finite difference scheme is proposed; moreover, a locally uniform (in space) convergence of the proposed numerical scheme is obtained, but the issue that matters to us is that the required number of samples grows like $h^{-d}$ (see \cite[Example 4.4]{FaToWa11}) where $d$ is the dimension whilst $h$ is the time discretization parameter. 
Also, because the space discretization is performed by a finite difference scheme on a uniform grid and the Monte Carlo sampling needs to be performed for each point in the grid, the algorithm complexity once more suffers exponentially with respect to the dimension.

\paragraph{DNNs for the Dirichlet problem on bounded domains} We mention two directions in the literature that aim at rigorously proving that DNNs can be used as numerical solvers without suffering from the curse of high dimensions. 
The first one is proposed in \cite{GrHe21} and it is the most related to our present work, so we shall frequently refer to it in what follows. The approach goes through the stochastic representation of the solution, and aims at designing a Monte Carlo sampler that can be used uniformly for all locations in the domain where the solution is approximated. However, some important issues remained open, like the fact that the obtained estimates depend on the volume of the domain which can nevertheless grow exponentially with respect to the dimension; these are going to be discussed in detail later.
The second approach is in \cite{NEURIPS2021_7edccc66} which is rather different and constructs the neural network progressively using a gradient descent method and then calculates the polynomial complexity of neural network approximation. 
A main feature of this second approach is that the construction of the DNN solver uses the theoretical spectral decomposition of the differential operator which is numerically not available, hence the obtained existence of the DNN approximator for the exact solution is of theoretical nature. 

\paragraph{DNNs for (linear) Kolmogorov PDEs} In the case of linear parabolic PDEs, DNNs solvers based on probabilistic representations have been proposed and numerically tested in \cite{BeBeGrJe21}.
A theoretical proof that DNNs are indeed able to approximate solutions to a class of linear Kolmogorov PDEs without suffering from the curse of dimensions has been provided in \cite{jentzen2018proof}. 
The strategy fundamentally aims at minimizing the error in $L^2(D;\lambda/\lambda(D))$ ($\lambda$ is the Lebesgue's measure), hence the existence of a DNN that approximates the solution without the curse of dimensions is in fact obtained on domains whose volumes increase at most polynomially with respect to the space dimension. 
In  \cite{BeBeGrJe21}, the authors include some numerical evidence that the $L^\infty(D)$-error also scales well with respect to the dimension, but some caution needs to be taken as mentioned in \cite[4.7 Conclusion]{BeBeGrJe21}. 
In the case of the heat equation on $\mathbb{R}^d$ it was proved in \cite{GoGrJeKoSi22} that any solution with at most polynomial growth can be approximated in $L^\infty([a,b]^d)$ by a DNN whose size grows at most polynomially with respect to $d$, the reciprocal of the prescribed approximation error, and $\max(|a|,|b|)$; the authors use heavily the representation of the solution by a standard Brownian motion shifted at the location point where the solution is approximated. 
In our paper we deal with the Poisson problem in a bounded domain in $\mathbb{R}^d$, and the main difficulty and difference at the same time, comes from the fact that the representing process (namely the Brownian motion stopped when it exits the domain) depends in a nonlinear way on the starting point where the solution is represented, and this dependence is strongly influenced by the geometry of the domain.

\paragraph{DNNs for semilinear parabolic PDEs on $\mathbb{R}^d$}  In the case of semilinear heat equation on $\mathbb{R}^d$ with gradient-independent nonlinearity, a rigorous proof that DNNs can be used as numerical solvers that do not suffer from the curse of dimensions can be found in \cite[Theorem 1.1]{HuJeKrNg20}. 
The approximating errors are considered in the $L^2([0,1]^d)$ sense, so in general, by a scaling argument, they depend on the volume of the domain where the solution is approximated. 
Deep learning methods for general semilinear parabolic PDEs on $\mathbb{R}^d$ have been proposed and efficiently tested in high dimensions in \cite{han2018solving}. 
Rigorous proofs that these type of deep solvers are indeed capable of approximating solutions to general PDEs without suffering from the curse of dimensions are still waiting to be derived.

\subsection{\bf Our contribution}\label{ss:contribution}
In the present work we are primarily interested in studying Monte Carlo and DNN numerical approximations for solutions to the Poisson boundary value problem \eqref{e:0} in bounded domains in $\mathbb{R}^d$, explicitly tracking the dependence of the Monte Carlo estimates as well as the size of the corresponding neural networks in terms of the spatial dimension, the reciprocal of the accuracy, the regularity of the domain and the prescribed source and boundary data.
There are several key issues that we tackle throughout the paper, so let us briefly yet systematically point them out here:

\paragraph{Overcoming the curse of high dimensionality} We outline here that a very important consequence of our results concerns the breaking of the curse of high dimensions in the sense that the size of the neural network approximating the solution $u$ to problem \eqref{e:0} adds at most a (low degree) polynomial complexity to the overall complexity (see {\rm {\bf Theorem} (Part II)} from above) of the approximating networks for the distance function and the data. Moreover, as typical in machine learning, we also show that despite the fact that the neural network construction is random, it breaks the dimensionality curse with high probability. 
In terms of the dimension $d$, our main results state, in particular, that if the domain is sufficiently regular (e.g. convex) then the complexity of the Monte Carlo estimator of the exact solution to \eqref{e:0} scales at most like $ d^3\log^4(d)$, whilst the DNN estimator of the same solution scales at most like $d^5\log^5 (d) {\rm S}$, where ${\rm S}$ is a cumulative size of the DNNs used to approximate the given data and the distance function to the boundary of the domain.
In contrast with the results from \cite{GrHe21}, our construction of the DNN approximators is explicit and such that their sizes do not depend on the volume of the domain.  
The herein obtained estimates should also be compared with the conclusion from \cite{NEURIPS2021_7edccc66} where the size of the network is $\mathcal{O}(d^{\log(1/\gamma)})$, where $\gamma$ represents the accuracy of the approximation, thus the degree grows with the allowed error. Also, the construction adopted in \cite{NEURIPS2021_7edccc66} is of theoretical nature in the sense that it guarantees the existence of DNNs with the desired properties, but is unclear how it could be implemented in practice.
In contrast, our schemes can be easily implemented using GPU computing, as discussed in \Cref{S:numerics}.

\paragraph{Low dimensions improvements for general bounded domains with Dirichlet data}  We should point out that our approach also has interesting consequences in low dimensions.  The key is that we can reuse the samples for the Monte Carlo solver to simultaneously approximate the solution for all points in the domain, and furthermore the number of the steps required by the designed Walk-on-Spheres algorithm does not depend on the starting point; these two features make the proposed scheme highly parallelizable, and GPU computing can be very efficiently employed. Moreover, this works for arbitrary bounded domains with quantitative estimates while for more regular domains one obviously obtain improved estimates.

\paragraph{General bounded domains and H\"older continuous data} 
Recall that in \cite{GrHe21} the domain is assumed to be convex, whilst in \cite{NEURIPS2021_7edccc66} it is of class $C^\infty$. 
In the present paper it is shown that the curse of high dimensions can be overcome for a general class of domains, namely those that satisfy a uniform exterior ball condition. 
As a matter of fact, our results are even more general, covering the case of a arbitrary bounded domain in $\mathbb{R}^d$ (see \Cref{thm:main}), but then the estimates are given in terms of the behavior of the function $v_D$ defined in \eqref{def:v} in the proximity of the boundary of the domain. 
Concerning the regularity of the source and boundary data, our assumption is that they are merely H\"older continuous. Recall that in \cite{GrHe21} the source and boundary data are assumed to be twice continuous differentiable.

\paragraph{$L^\infty(D)$ estimates} Recall that in \cite{GrHe21}, the accuracy is measured with respect to the $L^2$  norm. 
However, as pointed out in \cite{NEURIPS2021_7edccc66}, the estimates depend actually on the volume of the domain $D$, hence they implicitly exhibit an exponential dependence on the dimension for sufficiently large domains.  
In the present work we estimate the errors in the uniform norm which gives on one hand much better results. 
On the other hand, we prove that the uniform norm of the error is small with large probability,  whilst the approximation complexity depends on $D$ merely through its (annular) diameter. 
We emphasize that obtaining reliable estimates for the expectation or the tail probability of the Monte Carlo error computed in the uniform norm is in general highly nontrivial. 
For example, such estimates have been only very recently been obtained for the (linear) heat equation in $\mathbb{R}^d$ in \cite{GoGrJeKo22}, after some considerable effort. In our case, the difficulty arises mainly from the fact that we work in bounded domains. Nevertheless, our method is in some sense much simpler and could be easily transferred in other settings as well.

\paragraph{Walk-on-Spheres acceleration revisited} 
The walk on spheres (WoS) was introduced in \cite{Mu56} as a way of accelerating the calculation of integrals along the paths of Brownian motion.  
The standard walk-on-spheres uses the following scheme:  take some $x\in D$; 
then, instead of simulating the entire Brownian motion trajectory, we start by uniformly choosing a point on the sphere of maximal radius inside $D$. 
Then, this step is repeated until the current position enters some thin neighborhood of the boundary ($\varepsilon$-shell), where the chain is stopped. Thus, the chain so constructed  is used as an approximation for the Brownian motion started from $x$ and stopped at the exit time from the domain $D$; recall that the distribution of the latter variable is precisely the harmonic measure with pole $x$, hence it is used to represent the solution $u(x)$ to problem \eqref{e:0} with $f\equiv 0$.
Here, we modify this algorithm in two respects. 
Firstly, our stopping rule for the walk on sphere is deterministic, namely, we run the walk-on-spheres chain for a given number of steps, uniformly for all trajectories and all points in the domain. 
This differs from the stopping rule used in \cite{GrHe21}, and, in fact, to the one usually adopted in the literature. 
It has a number of advantages, but perhaps the most important one is that the neural network construction outlined here is explicit. 
In order to understand how large one should take such a deterministic stopping time in order to achieve the desired estimates, we have to investigate different estimates in less or more regular domains for the number of steps needed for the walk on spheres chain to reach the $\varepsilon$ neighborhood of the boundary. 
Secondly, our walk-on-spheres scheme is performed with the maximal radius replaced by a more general radius, which is not necessarily maximal and is compatible with ReLU DNNs.
Overall, we develop a generalized walk-on-spheres algorithm which is of self interest and which is much more compatible with parallel computing, when compared to the classical scheme. 

\paragraph{The core ideas are surprisingly simple and of general nature} 
Putting aside the WoS acceleration algorithm, the crucial ingredients are the following: the first one is that the Monte Carlo approximation $u_M^N(x)$ given by \eqref{e:5} 
of the exact solution $u$ to problem \eqref{e:0} is a.s. H\"older continuous with respect to $x$, yet with a H\"older constant which is exponentially large with respect to the number $M$ of WoS steps. The second one is to consider a uniform grid discretization of the domain $D$ and to approximate the solution $u$ in the sup-norm merely on this grid. 
The third ingredient is to employ the (otherwise very poor) H\"older regularity of $u_M^N$ to extrapolate the approximation from the grid to the entire domain. 
The grid needs to be taken extremely refined, and thus leads immediately to an exponential complexity in terms of $M$ and the diameter of $D$. 
Now the fourth ingredient comes into play, namely we use H\"offding's inequality combined with the union bound inequality in order to efficiently compensate for the exponential complexity induced by the uniform grid. 
We emphasize that the grid discretization is just an instrument to prove the main result which remains in fact grid-independent. This approach that uses an auxiliary uniform grid whose induced complexity is compensated by a concentration inequality is, as a matter of fact, simple and of general nature. Therefore, we expect that out approach can be easily employed for other classes of PDEs.

\paragraph{Universality with respect to given data} 
One useful feature of the estimator explicitly constructed in this paper is that it essentially approximates the  operator that maps the data (source and boundary) of problem  \eqref{e:0} into the corresponding solution $u$. 
In particular, it means that the DNN solvers constructed herein consist of the composition of two separate neural networks: one which approximates the  the source and boundary data and one for the above-mentioned operator.
In this light, the present method could be interpreted as an {\it operator learning} method, and once the operator is learned, the source and boundary data can be varied very easily, without any further training.

\paragraph{Explicit construction of the approximation}  One key element of our approach is the explicit formulation of the approximation.  This is reflected in the formula \eqref{e:5} where all elements are fully determined.  We should also insist that \eqref{e:5} is much simpler than a neural network and does not need any training.  On the other hand, this structure can be exploited to initialize a DNN with significantly less complexity than the guaranteed Monte Carlo construction we provide.  Once this initialization is done, we can train this network in order to further decrease the approximation error.

\subsection{ Organization of the paper} 
In Section~\ref{s:main} we present the main notations, the important quantities and the main results.  For the sake of readability of the paper we moved the proofs to Section~\ref{s:proofs}.  Section~\ref{s:main} in turn contains several subsections which we discuss now because they show the main approach.   Subsection~\ref{ss:representation} details the first probabilistic representation of the solution and various types of estimates based on the annular diameter of a domain, which is a certain one dimensional characterisation of the domain.  In Subsection~\ref{ss:nsteps}  the walk on spheres and its modified version enters the scene and we give the main estimates on the number of steps needed to get in the proximity of the boundary for general bounded and open sets in $\mathbb{R}^d$. Also here we introduce the class of  \emph{$\delta$-defective convex} domains, a class of bounded and open sets in $\mathbb{R}^d$ for which we provide better estimates for the number of steps to the proximity of the boundary.  Subsection~\ref{ss:errors} provides the main analysis of the modified walk-on-spheres chain stopped at a deterministic time (thus uniformly for all points in the domain and all samples) as opposed to the one in \cite{GrHe21} which is random and very difficult to control.  Here we estimate first the sup norm of  $u-u_{M}$ { where $u$ is the solution to \eqref{e:0} represented probabilistically by \eqref{eq:representation}, whilst $u_M$ is given by \eqref{eq:um}; the estimates are given in terms of the regularity of the data, the geometry of the boundary and the parameter $M$.}  If the annular diameter is finite then in fact the estimates depend only on the diameter and the annular diameter or the parameter $\delta$ (the convex defectivness).  Furthermore, Subsection~\ref{ss:meantailestimates} introduces the Monte Carlo estimator $u_{M}^{N}$ (see \eqref{e:5}) for $u_{M}$  and contains the main results, namely Theorem~\ref{thm:main} and  Corollary~\ref{coro:noKepsilon}.  
Subsection~\ref{ss:extension} contains the main extensions we need to deal with the walk on sphere algoritm.  Fundamentally, in order to construct either \eqref{eq:um} or \eqref{e:5} we need to extend the values of boundary data $g$ inside the domain and we do this under some regularity conditions.  

Section~\ref{s:NN} contains the neural network consequences of the main probabilistic results and benefits from the very careful preliminary construction of the modified walk on spheres.   We should only point out the key fact that from \eqref{e:5} namely that once we replaced $g$, $f$ and $\tilde{r}$ by  neural networks then the function $u_{M}^{N}$ also becomes  a neural network.  A word is in place here about the distance function $r$, the distance to the boundary.  In the original walk on spheres, one uses $r$ for the construction of the radius of the spheres.  We modified this into $\tilde{r}$.  The benefit is that if we have an approximation of $r$ by a neural network, we can easily construct $\tilde{r}$ which is already a  neural network.  This avoids complications which arise in \cite{GrHe21} from the approximation of $r$ by a neural network after the Monte-Carlo estimator is constructed.  The rest of the section here is judicious counting of the size of the neural network obtained for $u_{M}^{N}$ replacing $f$ and $g$ by their corresponding approximating networks.  

{ Finally, \Cref{S:numerics} is devoted to several numerical tests based on the Monte Carlo approach proposed and analyzed in the previous sections. Some key theoretical bounds are numerically validated and the PDE \eqref{eq:pd} below is numerically solved for some relevant domains in $\mathbb{R}^d$ for $d=10$ and $d=100$.}

\section{The presentation of the main results}\label{s:main}

{ As we already discussed}, this work concerns probabilistic representations and their DNN counterpart for the solution $u$ to problem
\begin{equation} \label{eq:pd}
\left\{\begin{array}{ll} \frac{1}{2}\Delta u=-f &\,\textrm{ in } D\subset \mathbb{R}^d\\[1mm]
\phantom{\frac{1}{2}\Delta}u=g &\,\textrm{ on }\partial D.
\end{array}\right.
\end{equation}
Throughout this paper, $D$ is bounded and open set in $\mathbb{R}^d$, $f$ is bounded on $D$ and $g$ is continuous on the boundary $\partial D$.
Further regularity shall also be imposed on $D\subset \mathbb{R}^d$, $f$ and $g$, so let us fix some notations:
{ We say that a set $D\subset \mathbb{R}^d$ is of class $C^k$ if its boundary $\partial D$ can be locally represented as the graph of a $C^k$ function.
}

\medskip
Before we proceed, let us remark that one could also consider the anisotropic operator $\nabla\cdot K\nabla$ instead of $\Delta$, where $K$ is a (homogeneous) positive definite symmetric matrix, without altering the forthcoming results.
This can be done mainly due to the following straightforward change of variables lemma.  For completeness, we also include its short proof in the smooth case.

\begin{lem} \label{lem:K}
For any given $K$ a $d\times d$ symmetric, positive-definite matrix, assume that $u$ is a classical solution to \eqref{eq:pd} with $\Delta$ replaced by $\nabla \cdot K\nabla$. 
If we take a $d\times d$ matrix $A$ such that $AA^T=K$, then denoting $D_A:=A^{-1}(D)$ and $v(x):=u(Ax), f_A(x)=f(Ax), g_A(x)=g(Ax), x\in D_A$, we have
{\begin{equation*}
\frac{1}{2}\Delta v=-f_A \,\textrm{ in } D_A, \qquad
v=g_A \textrm{ on }\partial D_A.
\end{equation*}}
\medskip
\noindent{Proof in \Cref{pf:1}.}
\end{lem}

\subsection{Probabilistic representation for Laplace equation and exit time estimates} \label{ss:representation}
We fix $B^0(t),t\geq0$,  to be a standard Brownian motion on $(\Omega, \mathcal{F}, \mathcal{F}_t, \mathbb{P})$ which starts from zero, that is $\mathbb{P}(B(0)=0)=1$.
Then, we set
\begin{equation*}
    B^x(t):=x+B^0(t), \quad t\geq0, \ x\in \mathbb{R}^d,
\end{equation*}
and recall that the law $\mathbb{P}^x:=\mathbb{P}\circ (B^x(\cdot))^{-1}$ on the path-space { $C([0,\infty);\mathbb{R}^d)$} is precisely the law of the Brownian motion starting from $x\in \mathbb{R}^d$.

By $\tau_{D^c}:=\tau^x_{D^c}$ we denote the first hitting time of $D^c:=\mathbb{R}^d\setminus D$ by $(B^x(t))_{t\geq 0}$, namely
\begin{equation*}
    \tau_{D^c}(\omega):=\inf\{t>0 : B^x(t,\omega)\in D^c\}, \quad \omega\in \Omega.
\end{equation*}

The following result is the fundamental starting point of this work. 
It is standard for sufficiently regular data, but under the next assumptions we refer to \cite[Theorem 6]{Ge92}.
\begin{thm}[\cite{Ge92}] \label{thm:representation}
Let $g\in C(\partial D)$ and $f\in L^\infty(D)$.
Then there exists a unique function $u\in C(D)\cap H^1_{loc}(D)$ such that $u$ is a (weak) solution to problem \eqref{eq:pd}.
Moreover, it is given by
\begin{equation}\label{eq:representation}
   u(x)=\mathbb{E}\{g(B^x(\tau_{D^c}))\}+\mathbb{E}\left\{\int_0^{\tau_{D^c}}f(B^x(t))\;dt\right\}, \quad x\in D.  
\end{equation}
\end{thm}

Let $v_D:\mathbb{R}^d\longrightarrow \mathbb{R}_+$ be given by
\begin{equation} \label{def:v}
    v_D(x):=\mathbb{E}\{\tau^x_{D^c}\}, \quad x\in \mathbb{R}^d.
\end{equation}
Most of the main estimates obtained in this paper are expressed in terms of the sup-norm of $v_D$ in the proximity of the boundary of $D$ (see e.g. \eqref{def:vepsilon}). 
In a following paragraph we shall explore such estimates for $v_D$ for domains that { satisfy the uniform exterior ball condition}.
Before that, let us start with the following consequence of \Cref{thm:representation}:
\begin{coro} \label{coro:v_inf}
The following assertions hold for $v_D$ given by \eqref{def:v}. 

\begin{enumerate}
    \item[i)] $|v_D|_\infty \leq {\sf diam}(D)^2/d$.
    \item[ii)] $v_D\in C^\infty(D)$ is the solution to the Poisson problem
{\be\label{eq:v}
-\frac{1}{2} \Delta v_D=1\,\textrm{ in } D,\qquad v_D=0\, \textrm{ on }\partial D.
\ee}
\end{enumerate}

\noindent{Proof in \Cref{pf:1}.}
\end{coro}

\paragraph{The annular diameter of a set in $\R^d$ and exit time estimates}
\label{ss:exit}

Now we explore more refined bounds for $v_D(x):=\mathbb{E}\{\tau^x_{D^c}\}$, $x\in D$, aiming at providing a general class of domains $D\subset \mathbb{R}^d$ for which $v_D(x)\leq {\rm adiam}(D)d(x,\partial D)$, where ${\rm adiam}(D)$ is a sort of \textit{annular diameter} of $D$ which in particular is a one dimensional parameter that depends on $D$.  For more details,  see \eqref{e:d1} below for the precise definition, and \Cref{prop:extime} for the precise result.

The first step here in understanding the exit problem from a domain which is not necessarily convex is driven by our first model, namely the annulus defined by $0<R_0<R_1$ as 
\[
 A(a,R_0,R_1)=\{ x\in\R^n : R_0<|x-a|<R_1 \}.  
\]
We set $A(R_0,R_1):=A(0,R_0,R_1)$.  

The reason we look at annular domains first is that they will allow us to extend the next result to any domain that satisfy the uniform exterior ball condition.

\begin{prop}
\label{prop:exitannulus}   
Take $d\ge 3$ and $D:= A(R_0,R_1)$. 
Then, for every point $x\in D$
\begin{equation}\label{e:1}
v_D(x)=\mathbb{E}\{ \tau^x_{D^c}\}\le d(x,\partial D)\frac{(R_1-R_0)R_1}{R_0}.
\end{equation}

\noindent
Proof in \Cref{pf:1}
\end{prop}

{ Now we extend the above result to a larger class of domains, namely those that satisfy the uniform exterior ball condition. 
To do it in a quantitative way, we need to} define first a notion of {\it annular diameter} of a set as follows:  
\begin{definition}
For a bounded and open set $D\subset \mathbb{R}^d$, and $x\in\partial D$, we set
\begin{equation}\label{e:i:6}
{\sf adiam}(D)_x:=\inf\left\{\frac{(R_{1}-R_{0})R_{1}}{R_{0}}: \exists a\in \mathbb{R}^d \mbox{ such that } x\in\partial B(a,R_0) \mbox{ and } D\subset A(a,R_0,R_1) \right\},
\end{equation}
with the convention $\inf\emptyset:=\infty$.
{Furthermore, we set} 
\begin{equation}\label{e:d1}
{\sf adiam}(D):=\sup_{x\in \partial D}{\rm adiam}(D)_x. 
\end{equation}
\end{definition}
Note that for a point $x\in \partial D$ and every pair $(R_0,R_1)$ such that $\exists a\in \mathbb{R}^d \mbox{ such that } x\in\partial B(a,R_0) \mbox{ and } D\subset A(a,R_0,R_1)$, the latter holds for the pair $(R_0,R_0+{\sf diam}(D))$.
Consequently,
\begin{equation}\label{e:i:loc:1}
{\sf adiam}(D)_x\le {\sf diam}(D)\left(1+\frac{{\sf diam}(D)}{R_{0}}\right), \quad \mbox{hence} \quad {\sf adiam}(D)_x\le {\sf diam}(D)\left(1+\frac{{\sf diam}(D)}{R_{0,x}}\right),
\end{equation}
where $R_{0,x}\in [0,\infty]$ denotes the radius of the largest closed ball in $\mathbb{R}^d\setminus D$ which contains $x$. 

Note that $D$ satisfies the exterior ball condition at $x\in \partial D$ if and only if ${\rm adiam}(D)_x<\infty$.

{ The above discussion leads to the following formal statement.
\begin{prop}
A bounded and open set $D\subset \mathbb{R}^d$ satisfies the uniform exterior ball condition if and only if ${\rm adiam}(D)<\infty$. 
If the radius of the exterior ball is at least $r_0>0$, we can estimate 
\[
{\sf adiam}(D)\le {\sf diam}(D)\left(1+\frac{{\sf diam}(D)}{r_0}\right).
\]
\end{prop} 
}
However, for instance a ball centered at $0$ from which we remove a cone with the vertex at the center does not have a finite {\sf adiam}. 
On the other hand, obviously, a convex bounded domain has finite {\sf adiam} and in fact, this is actually equal to the diameter of the set.  Indeed, one can see this by taking a tangent ball of radius $R_0$ and taking $R_1=R_0+{\sf diam}(D)$.  Letting $R_0$ tend to infinity we deduce that for a convex set $D$ we actually have ${\sf adiam}(D)={\sf diam}(D)$.  

Now we can present the main estimate of this paragraph. 
\begin{prop} \label{prop:extime}
If $D\subset \mathbb{R}^d$ has ${\sf adiam}(D)<\infty$, then for any $x\in D$,  
\[
 v(x):=\mathbbm{E}\{ \tau^x_{D^c}\} \le d(x,\partial D){\sf adiam}(D).  
\]
In particular, using \eqref{e:i:loc:1} and \eqref{e:d1}, we have
\[
 v(x)\le 2d(x,\partial D){\sf diam}(D)\left(1+\frac{{\sf diam}(D)}{R_0}\right), x\in D,  
\]
where recall that $\tau_{D^c}$ is the first exit time from $D$ and $R_0:=\inf\{R_{0,x}:x\in\partial D\}$. 

\medskip
\noindent{Proof in \Cref{pf:1}.}
\end{prop}

\subsection{Walk-on-Spheres (WoS) and \texorpdfstring{$\varepsilon$}{epsilon}-shell estimates revisited}
\label{ss:nsteps}

An important benefit of representation \eqref{eq:representation} is that the solution $u$ may be numerically approximated by the empirical mean of iid realizations of the random variables under expectation, thanks to the law of large numbers.
One way to construct such realizations is by simulating a large number of paths of a Brownian motion that starts at $x\in D$  and is stopped at the boundary $\partial D$.   However, as introduced by Muller in \cite{Mu56}, there is a much more (numerically) efficient way of constructing such realisations, based on the idea that \eqref{eq:representation} does not require the entire knowledge of how Brownian motion reaches $\partial D$.   This is clearer if one considers the case $f\equiv 0$, when the only information required in \eqref{eq:representation} is the location of Brownian motion at the hitting point of $\partial D$.  In this subsection we shall revisit and enhance Muller's method.

For any $x\in D$ let $r(x)\in[0,{\sf diam}(D)]$ denote the distance from $x$ to $\partial D$, or equivalently, the radius of the largest sphere centered at $x$ and contained in $D$, that is
\begin{equation} \label{eq:rdistance}
    r(x):=\inf\{|x-y| : y\in \partial D\}=\sup\{r>0 : B(x,r) \subseteq D\}.
\end{equation}
Clearly, $r$ is a Lipschitz function.

Recall that the standard WoS algorithm introduced by Muller \cite{Mu56} (see also \cite{sabtal}, \cite{DeLe} or \cite{deaconu} for more recent developments) is based on constructing a Markov chain that steps on spheres of radius $r(x)$ where $x\in D$ denotes its current position.
However, it is often the case, especially in practice, that merely an approximation $\overline{r}$ of the distance function is available, and not the exact $r$. 
This is the case, for example, if for computational reasons, $r$ is simply approximated with the (computationally cheaper) distance function to a polygonal surrogate for the domain $D$. 
Another situation, which is in fact central to this study, is given in \Cref{s:NN} below, when $r$ is approximated by a neural network.
In both cases, $r$ is approximated by a function $\overline{r}$ with a certain error. 
It turns out that considering the chain which walks on spheres of radius $\overline{r}(x),x\in D$ exhibits certain difficulties regarding the error analysis and also the construction of the chain itself.   We do not go into more details at this point, but we refer the reader to \Cref{rem:estimator}, iii)-iv) below for a more technical explanation of these issues.

Our strategy is to solve both of the above difficulties at once, by developing from scratch the entire analysis in terms of (a modification of) $\overline{r}$, and not in terms of $r$ as it is typically done. 
This motivates the following concept of distance.

\begin{defi}\label{defi:adistance}
Let $D\subset \mathbb{R}^d$ be a bounded open set and $r$ the distance function to the boundary.
Given $\varepsilon_0\geq 0$ and $\beta\in (0,1]$,  a Lipschitz function $\widetilde{r}:D\rightarrow [0,{\sf diam}(D)]$
is called a $(\beta,\varepsilon_0)$-distance on $D$ if
\begin{equation*}
i)\; 0\leq \widetilde{r}\leq r \;  \mbox{ on } D \quad \mbox{ and }\quad
ii)\; \widetilde{r}(x)\geq \beta r(x) \; \mbox{ for } x\in D, r(x)\geq \varepsilon_0.
\end{equation*}
When $\varepsilon_0=0$ we say that $\widetilde{r}$ is a $\beta$-distance, and if in addition $\beta =1$ then it is obvious that $\widetilde{r}=r$. 
\end{defi}

Notice that if $\tilde{r}$ is a $(\beta, \varepsilon_0)$-distance, then it is also a $(\beta,\varepsilon')$-distance for any $\varepsilon'\geq\varepsilon_0$. 
Thus if we fix a $\varepsilon_0$ and $\tilde{r}$ is a $(\beta,\varepsilon_0)$-distance, then, $\tilde{r}$ is also a $(\beta, \varepsilon)$-distance for any $\varepsilon\geq\varepsilon_0$. 

\begin{rem}\label{rem:adistance}  
Suppose that $\phi_r:D\rightarrow [0,\infty)$ is a Lipschitz function such  that $|\phi_r-r|_{\infty}\leq \epsilon$.
If $\varepsilon_0>2\epsilon$ and $0<\beta\leq 1-\frac{2\epsilon}{\varepsilon}$, then a simple computation yields that
$$
\widetilde{r}(x):=(\phi_r(x)-\epsilon)^{+}, \quad x\in D
$$
is a $(\beta,\varepsilon_0)$-distance on $D$.
For example if we take $\varepsilon_0=3\epsilon$, we can choose $\beta=1/3$, therefore, given any $\phi_r$, a Lipschitz $\epsilon$ approximation of $r$, there exists $\widetilde{r}$ which is $(1/3,3\epsilon)$-distance.  The moral is that we can always work with a $(\beta,\varepsilon_0)$-distance  with $\beta\ge 1/3$.
 
\end{rem}

This example already anticipates the amenability of $\tilde{r}$ from \Cref{rem:adistance} to the ReLU neural networks. 
Indeed, the positive function $x^{+}$ is precisely the non-linear activation function and from this standpoint, if $\phi_{r}$ is a neural network, we can argue that $\tilde{r}$ becomes also a ReLU neural network.   More on this in Section~\ref{s:NN}.

\paragraph{$\widetilde{r}$-WoS chain.} Let $U_n:\Omega\rightarrow S(0,1),\; n\geq1$ be defined (for simplicity) on the same probability space $(\Omega,\mathcal{F},\mathbb{P})$, independent and uniformly distributed, where $S(0,1)\subset \mathbb{R}^d$ is the sphere centered at the origin with radius $1$.
Let $(\widetilde{\mathcal{F}}_n)_{n\geq 0}$ be the filtration generated by $(U_n)_{n\geq 0}$, where $U_0=0$, namely
\begin{equation*}
    \widetilde{\mathcal{F}}_n:=\sigma(U_i \;:\; i\leq n), \;n\geq 0.
\end{equation*}

Also, let $\varepsilon_0\geq 0$, $\beta\in (0,1]$, and $\widetilde{r}$ be a $(\beta,\varepsilon_0)$-distance on $D$. 
For each $x\in D$, we construct the chain $({X}^{x}_n)_{n\geq 0}$ recursively by
\begin{align} 
{X}^{x}_0&:=x \label{wos0}\\
{X}^{x}_{n+1}&:={X}^{x}_n+\widetilde{r}({X}^{x}_n)U_{n+1},\; n\geq 0. \label{wosn}
\end{align}
Clearly, $({X}^{x}_n)_{n\geq 0}$ is a homogeneous Markov chain in $D$ with respect to the filtration $(\widetilde{\mathcal{F}}_n)_{n\geq 0}$, which starts from $x$ and has transition kernel given by
\begin{equation*}
    P f(x)=\int_{S(0,1)}f(x+\widetilde{r}(x)z) \;\sigma(dz), \quad x\in D, 
\end{equation*}
where $\sigma$ is the normalized surface measure on $S(0,1)$; that is, $\mathbb{E}\{f(X_n^{x})\}=P^nf(x), x\in D, f$ bounded and measurable.
We name it an $\widetilde{r}$-WoS chain.

We return now to problem \eqref{eq:pd}
\begin{equation*}
\begin{array}{ll} \frac{1}{2}\Delta u=-f \textrm{ in } D 
\text{ with }u=g \textrm{ on }\partial D,
\end{array}
\end{equation*}
which by Theorem \ref{thm:representation} admits the probabilistic representation \eqref{eq:representation}.
Further, we consider the following sequence of stopping times
\begin{equation}
    \tau_0^{x}=0,\quad \tau^{x}_{n+1}=\inf\{t>\tau_n^x : |B^x(t)-B^x(\tau^x_n)|\geq \widetilde{r}(B^x(\tau^x_n))\}, \;n\geq 1.
\end{equation}

\begin{rem} \label{rem:beta0distance}
It is straightforward that if $\varepsilon_0=0$, i.e. $\widetilde{r}$ is a $\beta$-distance, then
$$
\lim\limits_n\tau_n^{x} = \tau_{D^c}^x \quad \mathbb{P}^0\mbox{-a.s.} \quad\mbox{as well as} \lim\limits_{M\to\infty}X^x_{M}=B^x_{\tau^x_{\partial D}} \quad \mbox{in distribution}. 
$$
More generally, if $\varepsilon_0\geq 0$ and $\widetilde{r}$ is a $(\beta,\varepsilon_0)$-distance, then $\tau^{x}_{n}$ converges a.s. to some $\tau_\infty^x$ such that $\tau^x_{D_{\varepsilon_0}}\leq \tau_\infty^x\leq \tau^x_{D^c}$.
\end{rem}

The following result is a generalization of Lemma 3.4 from \cite{GrHe21}:
\begin{coro} \label{coro:representation}
Let $D\subset \mathbb{R}^d$ be a bounded open set, $g\in C(\overline{D}), f\in L^\infty(D)$, and $u$ be the solution to \eqref{eq:pd}.
Also, for $\varphi:B(0,1)\rightarrow \mathbb{R}$ bounded and measurable set
\begin{equation*}
K_0\varphi:=\mathbb{E}
\left\{ \int_0^{\tau^0_{B(0,1)^c}}\varphi(B^0(t)) \; dt \right\}.
\end{equation*}
If $\widetilde{r}$ is a $\beta$-distance ( i.e. it is a $(\beta,0)$-distance) then the following assertions hold:
\begin{enumerate}
    \item[i)] For all $x\in D$ we have \begin{equation*}
    u(x)=\lim\limits_{M\to\infty}\mathbb{E}\{g(X^x_{M}))\}+\mathbb{E}\left\{\sum\limits_{k\geq 1}\widetilde{r}^2(X^x_{k-1})K_0F_{x,k}\right\},\end{equation*}
    where $F_{x,k}(y)=f(X^x_{k-1}+\widetilde{r}(X^x_{k-1})y), y\in B(0,1)$, whilst $K_0$ acts on $F_{x,k}$ with respect to the $y$ variable.
    \item[ii)] The mapping $\mathcal{B}(B(0,1))\ni A\mapsto \mu(A):=d K_0 1_A\in [0,1]$ renders a probability measure $\mu$ on $B(0,1)$ with density $d G(0,y),y\in D$, where $G(x,y)$ is the Green function associated to $-\frac{1}{2}\Delta$ on $B(0,1)$. More explicitly, for $d\geq 3$ we have that $G(0,y)$ is proportional to $|y|^{2-d}-1, |y|<1.$
    \item[iii)] Let $Y$ be a real valued random variable defined on $\left(\Omega, \mathcal{F}, \mathbb{P}\right)$, with distribution $\mu$, such that $Y$ is independent of $(U_n)_n$.
    Then for all $x\in D$ we have \begin{equation}\label{eq:representationY}
    u(x)=\lim\limits_{M\to\infty}\mathbb{E}\{g(X^x_{M}))\}+\frac{1}{d}\mathbb{E}\left\{\sum\limits_{k\geq 1}\widetilde{r}^2(X^x_{k-1})f(X^x_{k-1}+\widetilde{r}(X^x_{k-1})Y)\right\}.\end{equation}
\end{enumerate}

\noindent{Proof in \Cref{pf:1}.}
\end{coro}

Consider the $\widetilde{r}$-WoS chain described above, and for each $x\in D$ let us define the required number of steps to reach the {\it $\varepsilon$-shell} of $\partial D$ by
\begin{equation}\label{eq:N_eps_x}
    {N}_\varepsilon^{x}:=\inf\{n\geq 0 \;:\; r
    ({X}_n^{x})<\varepsilon\},
\end{equation}
which is clearly an $(\widetilde{\mathcal{F}}_n)$-stopping time.

The estimates to be obtained in the next subsection, and consequently the size of the DNN that we are going to construct in order to approximate the solution to \eqref{eq:pd}, depend on how big $N_\varepsilon^{x}$ is.   The goal of this section is to provide upper bounds for $N^{x}_\varepsilon$, the number of steps the walk on spheres needs to get to the $\varepsilon$-shell.  
These estimates are first obtained for general bounded and open sets in $\mathbb{R}^d$, and then improved considerably for "defective convex" domains which are introduced in \Cref{defi:dconvex}, below.  We provide what are, to our knowledge, the strongest estimates when compared with the currently available literature, as well as rigorous proofs that rely on the general technique of Lyapunov functions.   Some of these results have some something in common with the results in \cite{binder2012rate}, though the estimates in there are not clearly determined in terms of the dimension.

It is essentially known that for a general bounded and open set $D\subset\mathbb{R}^d$, the average of $N_\varepsilon^{x}$ grows with respect to $\varepsilon$ at most as $({\sf diam}(D)/\varepsilon)^2$ (see \cite[Theorem 5.4]{Ky17} and its subsequent discussion), hence, by Markov inequality, $\mathbb{P}(N_\varepsilon^{x}\geq M)$ can be bounded by $({\sf diam}(D)/\varepsilon)^2/M$.
The next result shows that, in fact, $\mathbb{P}(N_\varepsilon^{x}\geq M)$ decays exponentially with respect to $M$, and independent of the dimension $d$; we do this for $(\beta,\varepsilon_0)$-distances.

\begin{prop}\label{prop:nsteps_general} Let $D\subset \mathbb{R}^d$ be bounded and open, $\varepsilon_0>0$, $\beta\in (0,1]$ and $\widetilde{r}$ be a $(\beta,\varepsilon_0)$-distance. 
Let $\varepsilon\geq \varepsilon_0 \mbox{ and } 1\leq p<1+\frac{\sqrt{d}}{2}$.
Then
\begin{equation*}
\mathbb{E}\left\{e^{\frac{(p-1)\beta^2\varepsilon^2}{p \left[{\sf  diam}(D)^2\right]}N_\varepsilon^x }\right\}
    \leq 
    \left(\frac{1}{1-\frac{4(p-1)^2}{d}}\right)^{1/(2p)} 2^{(p-1)/p}, \quad x\in D.
\end{equation*}
In particular, for $d\geq 2$ and $p=1+\frac{1}{2}$ we obtain
\begin{equation*}
\mathbb{E}\left\{e^{\frac{\beta^2\varepsilon^2}{3 \left[{\sf  diam}(D)^2\right]}N_\varepsilon^x }\right\}
    \leq 
    \left(\frac{1}{1-\frac{1}{d}}\right)^{1/3} 2^{1/3}\leq 2, \quad x\in D.
\end{equation*}
hence
\begin{equation*}
    \mathbb{P}(N^{x}_\varepsilon\geq M)\leq 2e^{-\frac{\beta^2\varepsilon^2}{3\left[{\sf diam}(D)^2\right]}M} \quad \mbox{for all }M\in \mathbb{N},
\end{equation*}
where ${\sf diam}(D)$ denotes the diameter of $D$.
\end{prop}
\noindent{\it Proof in \Cref{pf:2}.}

\medskip
If the domain $D$ is convex, then the average number of steps required by the WoS to reach the $\varepsilon$-shell is of order $\log (1/\varepsilon)$. 
This was shown by Muller in \cite{Mu56}, and also  reconsidered in \cite{binder2012rate}. 
{ However, it is of high importance to the present work to track an explicit constant in front of $\log (1/\varepsilon)$ in terms of the space dimension $d$.}
In this subsection we aim to clarify (the proof of) this result and extend it from convex domains to a larger class of domains, resembling the technique of Lyapunov functions from ergodic theory.
More importantly, as in the case of \Cref{prop:nsteps_general}, we are strongly interested in tail estimates for $N_\varepsilon^x$, not in its expected value.
In a nutshell, the idea is to show that the square root of the distance function to $\partial D$ is pushing the WoS chain towards the boundary at geometric speed.
This is a technique meant to be easily extended for more general operators, in a further work.

The following definition settles the class of domains for which the aforementioned estimate is going to hold, in the more general setting where the chain is constructed using a $(\beta,\varepsilon_0)$-distance instead of the true distance to the boundary $r$.
\begin{defi} \label{defi:dconvex}
We say that a Lipschitz bounded and open set $D\subset \mathbb{R}^d$ is "$\delta$-defective convex" if $\delta<1$ and
\begin{equation} \label{eq:dconvex}
    \Delta r\leq\frac{\delta}{2 \;\rm{rad}(D)} \quad \mbox{weakly on } D,
\end{equation}
where  we recall that $r$ is the distance function to the boundary $\partial D$, whilst $\rm{rad}(D):=\sup\limits_{x\in D} r(x)$.
\end{defi}

\begin{rem}
Recall that by \cite{Kuran}, $D$ is convex if and only if the signed distance function $r_s$ is superharmonic on $\mathbb{R}^d$, where
\begin{equation*}
    r_s:= r \mbox{ on } \overline{D} \mbox{ and } r_s:=-r \mbox{ on } \mathbb{R}^d\setminus D.
\end{equation*}
Hence a defective convex domain as defined by \eqref{eq:dconvex} is more general than a convex domain; in fact, even if \eqref{eq:dconvex} holds (on $D$) with $\delta=0$, it is not necessarily true that $D$ is convex, as explained in \cite{Kuran}.
To give an intuition on how a defective convex domain could differ from a convex domain, imagine a ball in $3D$ which is deformed into a defective convex domain by squeezing it slightly, or a straight cylinder in $3D$ which is bent mildly.
However, a defective convex domain can differ seriously from a convex domain, as revealed by the next two examples.
\end{rem}

\begin{exam} \label{ex:defective-anulus}
Let $A(R_1,R_2)\subset \mathbb{R}^d$ be an annulus of radii $R_1<R_2$,  namely $A(R_1,R_2)=\{x\in\R^d;R_1<|x|<R_2\}$.
If $\frac{R_2}{R_{1}}<1+\frac{\delta}{d-1}$ for some $\delta<1$, then $A(R_1,R_2)$ is a $\delta$-defective convex domain.

\medskip
\noindent{Proof in \Cref{pf:2}.}
\end{exam}  

\begin{exam}\label{ex:defective-tube}
We take $\Gamma$ to be a connected, compact orientable $C^2$ hypersurface in $\mathbb{R}^d$, with $d \geq 2$, endowed with the Riemannian metric $g$ induced by the embedding.  We denote by 
$k_1(x)\leq k_2(x)\le \ldots\leq k_{d-1}(x)$ 
the principal curvatures at $x\in\Gamma$.   The orientation is specified by a globally defined unit normal vector field $n:\Gamma\to \mathbb{S}^{d-1}$.
Then there exists a positive thickness $\eps$,
such that the tubular neighbourhood given by
\begin{equation}\label{layer.intro}
  D_\eps := \left\{x+\eps\,t\,n(x) \in \mathbb{R}^d \ \big| \ 
  (x,t) \in \Gamma \times (0,1) \right\}
  \,
\end{equation}
is $\delta$-defective convex.
In fact, $\eps$ can be chosen explicitly in terms of the principal curvatures of $\Gamma$.  

\medskip
\noindent{Proof in \Cref{pf:2}.}
\end{exam}

\begin{prop} \label{prop:logepsilon}
Let $D\subset \mathbb{R}^d$ be $\delta$-defective convex as in \eqref{eq:dconvex}.
Let $\varepsilon_0>0$, $\beta\in (0,1]$ and $\widetilde{r}$ be a $(\beta,\varepsilon_0)$-distance, and consider $P$ the transition kernel of the $\widetilde{r}$-WoS Markov chain $(X_n^{\cdot})_{n\geq 0}$.
If we set $V(x):=r(x)^{1/2}, x\in D$, then
\begin{equation} \label{eq:lyapunov}
    \forall \varepsilon\geq \varepsilon_0: \quad P V(x)\leq \left(1-\frac{\beta^2(1-\delta)}{4d}\right)V(x), \quad x\in D.
\end{equation}
As a consequence, 
\begin{equation} \label{eq:exptail}
\forall \varepsilon\geq \varepsilon_0: \quad\mathbb{P}(N_\varepsilon^{x}> M)\leq\mathbb{P}(r(X_M^x)\geq \varepsilon) \leq \left(1-\frac{\beta^2(1-\delta)}{4d}\right)^M \frac{V(x)}{\sqrt{\varepsilon}}, \quad x\in D, 
\end{equation}
and if $\delta_d:=1-\frac{\beta^2(1-\delta)}{4d}$,  then for any $1<a<1/\delta_d$
\begin{equation}\label{eq:logepsilon}
    a^{\mathbb{E}\left\{ N_\varepsilon^{x}\right\} }\leq 
    \mathbb{E}\left\{ a^{N_\varepsilon^{x}}\right\}
    \leq 1+\frac{a}{1-a\delta_d}\frac{V(x)}{\sqrt{\varepsilon}}, \quad x\in D.
\end{equation}

\medskip
\noindent{Proof in \Cref{pf:2}.}
\end{prop}

\begin{rem}\label{rem:any D is defective}
{ It can be shown that \Cref{prop:logepsilon} still holds if condition \eqref{eq:dconvex} is satisfied merely in some strict neighbourhood of the boundary. In particular, in view of \Cref{ex:defective-tube}, \Cref{prop:logepsilon} holds for all bounded and open sets with smooth boundary, but the estimates would also depend on the thickness of the neighbourhood where condition \eqref{eq:dconvex} is fulfilled. Going even further, any bounded and open set that can be uniformly approximated from inside by smooth domains enjoys, for fixed $\varepsilon$, a similar estimate with respect to $M$ and $D$ as in \eqref{eq:exptail}, with $\delta$ possibly depending on $\varepsilon$. This behavior is numerically confirmed by {\bf Test 1} in \Cref{S:numerics} for annular hypercubes, and it is going to be analyzed theoretically in a forthcoming work.}
\end{rem}

\subsection{WoS stopped at deterministic time and error analysis}
\label{ss:errors}

Throughout this subsection we assume that $u$ is the solution to problem \eqref{eq:pd}, hence \Cref{thm:representation} and \Cref{coro:representation} are applicable.
Moreover, we keep all the notations from the previous subsections.
Before we proceed with the main results of this subsection, let us emphasize several aspects that are essential to this work.  To make the explanation simple, assume that $\widetilde{r}=r$, that is the WoS chain is constructed using the $(1,0)$-distance $r$.
Given $x\in D$, a usual way to employ WoS Markov chain in order to approximate $u(x)$ through the representation furnished by \Cref{coro:representation}, is to start the chain from $x$ and run it until it reaches the $\varepsilon$-shell, for some given $0<\varepsilon<<1$.
In other words, $(X_k^x)_{k\geq 0}$ is usually stopped at the (random) stopping time $N_\varepsilon^x$ and $u$ represented by \eqref{eq:representation} is then approximated with
\begin{equation*}
    u_\varepsilon(x):=\mathbb{E}\{g(X^x_{N_\varepsilon^x})\}+\frac{1}{d}\mathbb{E}\left\{\sum\limits_{ k=1}^{ N_\varepsilon^x}r^2(X^x_{k-1})f(X^x_{k-1}+r(X^x_{k-1})Y)\right\}.
\end{equation*}
The intuition behind is that stopping WoS chain at the $\varepsilon$-proximity of the boundary should provide a good approximation of how the Brownian motion first hits the boundary $\partial D$, and estimates that certify this fact are in principle well known.
As discussed in the previous subsection, the number of steps required to reach the $\varepsilon$-shell is small, especially if the domain is (defective) convex, which eventually leads to a fast numerical algorithm.

From the point of view of this work (also of \cite{GrHe21}), the fundamental inconvenience of the above stopping rule is that it depends strongly on the starting point $x$, mainly through $N_\varepsilon^x$.
In other words, although computationally efficient for estimating a single value $u(x)$, the above \textit{point estimate} $u_\varepsilon$ of $u$ is expected to fail at overcoming the curse of high dimensions for solving \eqref{eq:pd} globally in $D$.
Moreover, if one aims at constructing a (deep) neural network architecture based on the above representation, as considered in \cite{GrHe21} and also in \Cref{s:NN} below, $x$ would be the input while $N_\varepsilon^x$ would give the number of layers; however, the latter should be independent of $x$, which is obviously not the case. 
To deal with this architectural impediment, in \cite{GrHe21} the authors proposed $\sup\limits_{x\in D}N_\varepsilon^x$ as a random time to stop the WoS chain. 
However, beside the measurability and the stopping time property issues for $\sup\limits_{x\in D}N_\varepsilon^x$, it is still unclear, at least to us, that $\mathbb{E}\left\{\sup\limits_{x\in D}N_\varepsilon^x\right\}<\infty$ and that this expectation does not depend on the dimension $d$.

Anyway, our approach is consistently different.  Instead of stopping the WoS chain at a random stopping time, be it $N_\varepsilon^x$ or $\sup\limits_{x\in D}N_\varepsilon^x$, the idea is to stop the chain after a deterministic number of steps, say $M$, independently of the starting point $x$. 
Such a choice turns out to be feasible, and it not only avoids the above mentioned issues concerning $\sup\limits_{x\in D}N_\varepsilon^x$, but it eventually provides a way to break the curse of high dimensions for solving \eqref{eq:pd}, merely using WoS algorithm but in a global fashion.
Furthermore, in terms of neural networks, this strategy would also render a way of explicitly constructing a corresponding DNN architecture, that could be easily sampled, and why not, further trained.
Also, as already mentioned in \Cref{rem:adistance}, when we shall deal with neural networks in \Cref{s:NN} the distance to the boundary $r$ needs to be replaced by an approximation given by a DNN, with a certain error.
Therefore, in light of \Cref{defi:adistance} and \Cref{rem:adistance}, we shall work instead with a $(\beta,\varepsilon_0)$-distance $\widetilde{r}$ on $D$, for some $\varepsilon_0>0$ and $\beta\in(0,1]$ properly chosen. 
Having all these in mind, the aim of this subsection is to estimate the error of approximating the solution $u$ with
\begin{equation} \label{eq:um}
    u_M(x):=\mathbb{E}\left\{ g(X_M^{x})+\sum\limits_{k= 1}^M \widetilde{r}^2(X^{x}_{k-1})K_0f(X^{x}_{k-1}+\widetilde{r}(X^{x}_{k-1})\cdot)\right\} \quad \mbox{ for all } x\in D,
\end{equation}
for a given (deterministic) number of steps $M\geq 1$ that does not depend on $x\in D$.

To keep the assumption on the regularity of $D$ as general as possible, the forthcoming estimates shall be obtained in terms of the function $v$ defined by \eqref{def:v} or \eqref{eq:v}, more precisely in terms of the behavior of $v$ near the boundary measured for each $\varepsilon >0$ (in fact we shall consider only $\varepsilon\geq \varepsilon_0>0$) by
\begin{equation} \label{def:vepsilon}
    |v| _{\infty}(\varepsilon ):=\sup\{v(x): r(x):=d(x,\partial D)\leq \varepsilon \}.
\end{equation}
{
Though this is our primary measure of the geometry of the boundary, we can in fact refine things by defining  
\begin{equation} \label{def:wepsilon} 
\begin{split}
v(x,\varepsilon):&=\mathbb{E}\left\{v(B^x_{\tau_{N_\varepsilon^{x}}^{x}})\right\} \\
v(F,\varepsilon):&=\sup_{x\in F}v(x,\varepsilon) \text{ for }F\subset D. 
\end{split}
\end{equation}

We include here a small result which reveals the main properties we need further on. 

\begin{prop}\label{p:vs}
We have for every $\varepsilon\geq \varepsilon_0$ that
\begin{equation}\label{e:w<v}
v(x,\varepsilon)\le |v|_\infty(\varepsilon). 
\end{equation}
If $\varepsilon_0=0$, then for any bounded and open set $D$ and any compact $F\subset D$,   
\begin{equation}
\lim_{\varepsilon\to 0}v(F,\varepsilon)=0.
\end{equation}
\end{prop}
Since the first inequality is clear, let us only point how one can prove the second part, by using the observation that for any stopping time $\tau$,  $v(B_{\tau\wedge t})$ is a bounded and non-negative right-continuous supermartingale, thus converges for $t\to\infty$.  As a consequence, we obtain that $v(x,
\varepsilon)$ converges to $0$.   On the other hand, we just observe that $v(x,\varepsilon)$ is non-increasing in $\varepsilon$, thus by Dini's theorem, we get the uniform convergence to $0$.   
}

Let us first consider the case of homogeneous boundary conditions, namely $g\equiv 0$.
\begin{prop} \label{prop:errorspoisson}
Let $D\subset \mathbb{R}^d$ be bounded and open.
Let $\varepsilon_0>0$, $\beta\in (0,1]$, $\widetilde{r}$ be a $(\beta,\varepsilon_0)$-distance, $f\in L^\infty(D)$, $M\in \mathbb{N}^\ast$ and $u$ be the solution to \eqref{eq:pd} with $g\equiv 0$.
If $u_M$ is given by \eqref{eq:um}
then 
\begin{equation}\label{e:f:0}
    |u(x)-u_M(x)|\leq |f|_\infty \left[ v(x,\varepsilon) + \frac{2}{d} {\sf diam}(D)^2 e^{-\frac{\beta^2\varepsilon^2}{3{\sf diam}(D)^2}M}\right ] \quad \mbox{ for all } \varepsilon\geq \varepsilon_0.
\end{equation}
In particular, 
\begin{equation}\label{e:f:1}
    \sup_{x\in D}|u(x)-u_M(x)|\leq |f|_\infty \left[ |v|_{\infty}(\varepsilon) + \frac{2}{d} {\sf diam}(D)^2 e^{-\frac{\beta^2\varepsilon^2}{3{\sf diam}(D)^2}M}\right ] \quad \mbox{ for all } \varepsilon\geq \varepsilon_0.
\end{equation}

\medskip
\noindent{Proof in \Cref{pf:3}.}
\end{prop}

Let us treat now the in-homogeneous Dirichlet problem, this time taking $f\equiv 0$.

\begin{prop} \label{prop:errorsdirichlet}
Let $D\subset \mathbb{R}^d$ be bounded and open.
Let $\varepsilon_0>0$, $\beta\in (0,1]$, $\widetilde{r}$ be a $(\beta,\varepsilon_0)$-distance, $u$ be the solution to \eqref{eq:pd} with $g\in C(\overline{D})$ and $f\equiv 0$.
Further, for each $M\in \mathbb{N}_\ast$ consider that $u_M$ is given by \eqref{eq:um}
If $g$ is $\alpha$-H\"{o}lder on $\overline{D}$ for some $\alpha\in [0,1]$ then
\begin{equation}\label{e:g:0}
|u(x)-u_M(x)|\leq 2d^{\alpha/2}|g|_\alpha\cdot v(x,\varepsilon)^{\alpha/2} + 4|g|_\infty e^{-\frac{\beta^2\varepsilon^2}{3{\sf diam}(D)^2}M}\quad \mbox{ for all } \varepsilon\geq \varepsilon_0,
\end{equation}
and in particular, 
\begin{equation}\label{e:g:1}
\sup\limits_{x\in D}|u(x)-u_M(x)|\leq 2d^{\alpha/2}|g|_\alpha\cdot |v|^{\alpha/2}_{\infty}(\varepsilon) + 4|g|_\infty e^{-\frac{\beta^2\varepsilon^2}{3{\sf diam}(D)^2}M}\quad \mbox{ for all } \varepsilon \geq \varepsilon_0.
\end{equation}

\medskip
\noindent{Proof in \Cref{pf:3}.}
\end{prop}

We can now superpose \Cref{prop:errorspoisson} and \Cref{prop:errorsdirichlet} to obtain the following key result:

\begin{thm} \label{thm:stepdet}
Let $D\subset \mathbb{R}^d$ be bounded and open. 
Let $\varepsilon+0>0$, $\beta\in (0,1]$, $\widetilde{r}$ be a $(\beta,\varepsilon_0)$-distance, and $u$ denote the solution to \eqref{eq:pd} with $g\in C(\overline{D})$ and $f\in L^\infty(D)$.
Further, for each $M\in \mathbb{N}_\ast$ let $u_M$ be given by \eqref{eq:um}.
If $g$ is $\alpha$-H\"{o}lder on $\overline{D}$ for some $\alpha\in [0,1]$, then for all $\varepsilon\geq \varepsilon_0$ we have
$$
|u(x)-u_M(x)|\leq 2d^{\alpha/2}|g|_\alpha\cdot v^{\alpha/2}(x,\varepsilon)+ |f|_\infty v(x,\varepsilon) + (4|g|_\infty+\frac{2}{d}{\sf diam}(D)^2|f|_\infty) e^{-\frac{\beta^2\varepsilon^2}{3{\sf diam}(D)^2}M},
$$
in particular, we get that 
\[
\sup_{x\in D}|u(x)-u_M(x)|\leq 2d^{\alpha/2}|g|_\alpha\cdot |v|_\infty^{\alpha/2}(\varepsilon)+ |f|_\infty |v|_\infty(\varepsilon) + (4|g|_\infty+\frac{2}{d}{\sf diam}(D)^2|f|_\infty) e^{-\frac{\beta^2\varepsilon^2}{3{\sf diam}(D)^2}M}.
\]
\end{thm}

Let us point out that if ${\sf adiam}(D)<\infty$ (see \eqref{e:i:6} below), then $v(x,\varepsilon)$ and $|v|_\infty(\varepsilon)$
involved above may be replaced by $\varepsilon \; {\sf adiam}(D)$.
Furthermore,
when the domain is defective convex (see \cref{ss:nsteps}), we can improve considerably the above error estimates with respect to the required number of WoS steps, $M$. 
This can be done analogously to the proofs of \Cref{prop:errorspoisson} and \Cref{prop:errorsdirichlet}, just by replacing the tail estimate given by \Cref{prop:nsteps_general} with the one provided by \Cref{prop:logepsilon}.
Therefore, we give below the precise statement, but we skip its proof.

\begin{coro} \label{coro:adiam1}
In the context of \Cref{thm:stepdet}, the following additional assertions hold:

\begin{enumerate}
    \item[i)] If { $D$ satisfies the uniform exterior ball condition}, then by \Cref{prop:extime}, $|v|_\infty(\varepsilon)$
involved in the above estimate may be replaced with $\varepsilon \; {\sf adiam}(D)$, { where recall that ${\sf adiam}(D)$ is given by \eqref{e:d1}.}
    \item[ii)] If the domain $D$ is $\delta$-defective convex ($\delta<1$, see \eqref{eq:dconvex} for the definition) so that the conclusion from \Cref{prop:logepsilon} is in force, then the factor $e^{-\frac{\beta^2\varepsilon^2}{3{\sf diam}(D)^2}M}$ from the above estimate can be replaced by \\ $\left(1-\frac{\beta^2(1-d)}{4d}\right)^M \sqrt{\frac{{\sf diam}(D)}{\varepsilon}}$.
\end{enumerate}
\end{coro}

\subsection{Monte-Carlo approximations: mean versus tail estimates}
\label{ss:meantailestimates}
We place ourselves in the same framework as before, namely: $D\subset \mathbb{R}^d$ is {bounded and open}, $g\in C(\overline{D})$, $f \in L^\infty(D)$, and $u$ is the solution to \eqref{eq:pd}.

Further, let $(U^i_n)_{n\geq 1, i\geq 1}$ be a family of independent and uniformly distributed random variables on $S(0,1)$, $\widetilde{r}$ be a $(\beta,\varepsilon_0)$-distance on $D$ for some $\beta \in (0,1]$ and $\varepsilon_0>0$ (see \Cref{defi:adistance}), and set:
\begin{align*}
    X^{x,i}_{n+1}&:=X_n^{x,i}+\widetilde{r}(X_n^{x,i})\cdot U^i_{n+1},\;\;n\geq 0, i\geq 1.
\end{align*}
On the same probability space as $(U^i_n)_{n\geq 1, i\geq 1}$, let $(Y^i)_{i\geq 1}$ be iid random variables with distribution $\mu$ given by \Cref{coro:representation}, ii), such that the family $(Y^i)_{i\geq 1}$ is independent of $(U^i_n)_{n\geq 1, i\geq 1}$.

For $N,M \in \mathbb{N}^\ast$ let $u_M$ be given by \eqref{eq:um}, and consider the Monte Carlo estimator
\begin{equation}\label{eq:MCestimator}
    u_{M}^{N}(x):=\frac{1}{N}\sum_{i=1}^N\left[ g(X^{x,i}_{M})+\frac{1}{d}\sum\limits_{k= 1}^M \widetilde{r}^2(X^{x,i}_{k-1})f\left(X^{x,i}_{k-1}+\widetilde{r}(X^{x,i}_{k-1})Y^i\right)\right], x\in D.
\end{equation}
As in \Cref{coro:representation}, iii), we have
\begin{equation*}
   \mathbb{E}\left\{ u_{M}^{N}(x)\right\}=u_M(x), \quad x\in D, N\geq 1.
\end{equation*}

\begin{rem} \label{rem:estimator}
At this point we would like to point out that estimator $u_M^{N}$ is different than the one employed in \cite{GrHe21}, Proposition 4.3 in several main aspects:
\begin{enumerate}
    \item[i)] The first aspect was already anticipated in the beginning of \Cref{ss:errors}, namely instead of stopping the WoS chain at $\sup\limits_{x\in D}N_\varepsilon^x$ which is a random time that is difficult to handle both theoretically and practically, we simply stop it a deterministic time $M$ which is going to be chosen according to the estimates obtained in \Cref{thm:main} below and its two subsequent corollaries.    
    \item[ii)] The second aspect is that the estimator used in \cite{GrHe21} considers $N$ iid samples drawn from $\mu$, for each of the $N$ iid samples drawn from $(U_n)_{n\geq 1}$, leading to a total of $N^2$ samples. 
In contrast, $u_M^{N}$ requires merely $N$ samples, because $Y$ and $(U_n)_n$ are sampled simultaneously (and independently), on the same probability space.
    \item[iii)] The third aspect is more subtle: In \cite{GrHe21}, the Monte Carlo estimator of type \eqref{eq:MCestimator} is constructed based on a given DNN approximation $\overline{r}$ of the distance to the boundary $r$, for any prescribed error, let us say $\eta$. 
    Then, the approximation error of the solution is obtained based on the error of the Monte Carlo estimator constructed with the exact distance $r$, and on how such an estimator varies when $r$  is replaced with $\overline{r}$. However, the latter source of error scales like $2^{N_{\varepsilon}}\eta$, where $N_{\varepsilon}:=\sup\limits_{x\in D}N_\varepsilon^x$.
    To compensate this explosion of error, $\eta$ has to be taken extremely small, and to do so, in \cite{GrHe21} it is assumed that $\overline{r}$ can be realized with complexity $O(\log(1/\eta))$; The authors show that such a complexity can indeed be attained for the case of a ball or a hypercube in $\mathbb{R}^d$, and probably can be extended to other domains with a nice geometry. 
    Our approach is different and the key ingredient is to rely on the notion of $(\beta,\varepsilon)$-distance introduced in $\Cref{defi:adistance}$. 
    More precisely, using \Cref{rem:adistance} we can replace $\overline{r}$ by some $(\beta,\varepsilon)$-distance $\widetilde{r}$ at essentially no additional cost, and rely on the herein developed analysis for $\widetilde{r}$-WoS. 
    This approach turns out to avoid the additional error of order $2^{N_{\varepsilon}}\eta$ mentioned above, in particular we shall be able to consider domains whose distance function to the boundary may be approximated by a DNN merely at a polynomial complexity with respect to the approximation error.
    \item[iv)] Another issue regards the construction of the WoS chain itself. 
    Because $\overline{r}$ from iii) may be strictly bigger than $r$, for a given position $x\in D$, the sphere of radius $\overline{r}(x)$ might exceed $D$, so there is a risk that the WoS chain leaves the domain $D$. 
    In particular, if one constructs the WoS chain based on $\overline{r}$, then in order to make the analysis rigorous the boundary data $g$ and the source $f$ should be extended also to the complement of the domain $\overline{D}$.
    Fortunately, this issue is completely avoided by considering $\widetilde{r}$-WoS chains (as it is done in this work), since by definition $\widetilde{r}\leq r$ on $D$.
\end{enumerate}
\end{rem}

Let us begin with the following mean estimate in $L^2(D)$:
\begin{prop}\label{prop:meanestimates}
Let $\varepsilon_0>0$, $\beta\in (0,1]$, and $\widetilde{r}$ be a $(\beta,\varepsilon_0)$-distance. 
Then for all $N, M\in \mathbb{N}$, $\gamma \geq 0$,
\begin{equation} \label{eq:meanestimate}
\mathbb{E}\left\{ \left | u (\cdot)-u_M^{N}(\cdot) \right |^2_{L^2(D)} \right\}
\leq 2\lambda (D)\left[ \sup\limits_{x\in D} |u(x)-u_M(x)|^2+\frac{2(|g|^2_\infty+\frac{1}{d^3}M|f|^2_\infty {\sf diam}(D)^4}{N}\right],
\end{equation}
where $\lambda$ is the Lebegue measure on $\mathbb{R}^d$, whilst $u_M$ and $u_M^N$ are given by \eqref{eq:um} and \eqref{eq:MCestimator}. 
In particular, the above inequality can be made more explicit by employing the estimates for $\sup\limits_{x\in D} |u(x)-u_M(x)|$ obtained in \Cref{thm:stepdet} and \Cref{coro:adiam1}, depending on the regularity of $D$ and $g$.

\medskip
\noindent{Proof in \Cref{pf:4}.}
\end{prop}

\begin{rem}
Note that as in \cite{GrHe21}, Section 4, the above error estimate depends on the volume $\lambda(D)$. 
When $\lambda(D)$ scales well with the dimension (e.g. at most polynomially), then \eqref{eq:meanestimate} can be employed to overcome the curse of high dimensions; in fact, if $D$ is a subset of a hypercube whose side has length less then some $\delta<1$, then $\lambda(D)\leq \delta^d$, hence, in this case, the factor $\lambda(D)$ improves the mean squared error exponentially with respect to $d$.
However, $\lambda(D)$ may also grow exponentially with respect to $d$, and the above estimate can not be used to construct a neural network whose size scales at most polynomially with respect to the dimension. Therefore, our next (and in fact) main goal is to solve this inconvenience, by looking at tail estimates for the Monte-Carlo error; one key idea is to quantify the error using the sup-norm instead of the $L^2(D)$-norm.
\end{rem}

Before we move forward, we recall the notion of a regular domain.  We say that a bounded and open set $D$ is regular if for any continuous function $g$ on the boundary, the harmonic function $u$ with the boundary condition $f$ is continuous in $\bar{D}$.

As announced in the above remark, we conclude now with the central result of this paper.

\begin{thm}\label{thm:main}
Keep the same framework and notations as in the beginning of this subsection.  Fix a small $\varepsilon_0>0$, $\beta\in (0,1]$, $\widetilde{r}$ a $(\beta,\varepsilon_0)$-distance, and consider $u_M$ and $u_M^N$ given by \eqref{eq:um} and \eqref{eq:MCestimator}.
{  Also, assume that $f$ and $g$ are $\alpha$-H\" older on $D$ for some $\alpha\in (0,1]$.
Then, for any compact subset $F\subset D$, for all $N, M, K \geq 1$, $\gamma>0$ and $\varepsilon\geq\varepsilon_0$,  
\begin{equation} \label{eq:main1}
\mathbb{P}\left( \sup_{x\in F} \left | u(x)-u_M^{N}(x)\right| \geq \gamma \right) \leq 2 \exp\left(C_1(M,K,d)-\frac{\left((\gamma-A(F,M,K,d,\varepsilon))^+\right)^2}{C_2(M,d)}N \right),
\end{equation}
where 
\begin{equation}\label{e:C_12}
\begin{split}
C_1(M,K,d):&=d\left(\lceil M/\alpha\rceil\log(2+|\widetilde{r}|_1)+\log(K)\right) \\ 
C_2(M,d):&=|g|_\infty+\frac{M}{d}{\sf diam}(D)^2|f|_\infty \\
\end{split}
\end{equation}
and
\begin{equation}\label{e:A}
\begin{split}
A(F,M,K,d,\varepsilon)&:=2\left(|g|_\alpha+\frac{{\rm diam}(D)^2|f|_\alpha+2{\rm diam}(D)|f|_\infty}{d}\right)\left(\frac{{\sf diam}(D)}{K}\right)^\alpha \\
&\quad + 2d^{\alpha/2}|g|_\alpha\cdot v(F,\varepsilon)^{\alpha/2}+|f|_\infty v(F,\varepsilon) + (4|g|_\infty+\frac{2}{d}{\sf diam}(D)^2|f|_\infty) e^{-\frac{\beta^2\varepsilon^2}{3{\sf diam}(D)^2}M} 
\end{split}
\end{equation}

Moreover, if we set
\begin{equation}\label{e:B}
\begin{split}
B(M,K,d):&=C_2(M,d)(\sqrt{C_1(M,K,d)+\log(2)}+1)\\
&=\left(|g|_\infty+\frac{M}{d}{\sf diam}(D)^2|f|_\infty\right)\left(\sqrt{d\left(\lceil M/\alpha\rceil\log(2+|\widetilde{r}|_1)+\log(K)\right)+\log(2)}+1\right).
\end{split}
\end{equation}
then we also have the estimate on the expectation of the total error in the form 
\begin{equation}\label{eq:main:e}
\mathbb{E}\left\{\sup_{x\in F} \left | u(x)-u_M^{N}(x)\right|\right\}\le A(F,M,K,d,\varepsilon)+\frac{B(M,K,d)}{\sqrt{N}}.
\end{equation}
As a consequence, from \Cref{p:vs}, for any compact set $F\subset D$, 
\begin{equation}\label{e:main:f}
\lim_{M\to \infty}\lim_{N\to\infty}\mathbb{E}\left\{\sup_{x\in F} \left | u(x)-u_M^{N}(x)\right|\right\}=0.
\end{equation}
For regular domains, we can take $F=D$.  

Furthermore, for any domain, we can replace $v(F,\epsilon)$ with $|v|_\infty(\varepsilon)$. Moreover, if the domain satisfies the exterior ball condition, then in \eqref{e:A} we can take $F=D$ and replace $|v|_\infty(\varepsilon)$ by $\varepsilon \; {\rm adiam}(D)$.

If the domain $D$ is $\delta$-defective convex,  we can replace $e^{-\frac{\beta^2\varepsilon^2}{3{\sf diam}(D)^2}M}$ from the definition of A in  \eqref{e:A} with $\left(1-\frac{\beta^2(1-\delta)}{4d}\right)^M \sqrt{\frac{{\sf diam}(D)}{\varepsilon}}$.

}

\medskip
\noindent{Proof in \Cref{pf:4}.}
\end{thm}






%

\begin{rem}
Note that the left hand side of \eqref{eq:main1} does not depend on $K$, and if $\widetilde{r}$ is a $\beta$-distance, then it does not depend on $\varepsilon$ as well.
Therefore, in the right hand side of \eqref{eq:main1} one may take the infimum with respect to $K\geq 1$, and if $\widetilde{r}$ is a $(\beta,\varepsilon_0)$-distance then one may also take the infimum with respect to $\varepsilon\geq \varepsilon_0$; but
optimizing the previously obtained bounds in this way may be cumbersome.
Anyway, convenient bounds can be easily obtained from particular choices of $K$ and $\varepsilon$, so let us do so in the sequel.
\end{rem}

Let us conclude this subsection with the following consequence obtained for some convenient choices for $K$ and $\varepsilon$.
\begin{coro} \label{coro:noKepsilon}
Let $D\subset \mathbb{R}^d$ such that { it satisfies the uniform exterior ball condition}. 
Further, let $f$ and $g$ be $\alpha$-H\"{o}lder on $\overline{D}$ for some $\alpha\in [0,1]$, $\gamma>0$ be a prescribed error, $\eta>0$ be a prescribed confidence, and $\widetilde{r}$ be a $(\beta,\varepsilon_0)$-distance with $\beta\in (0,1]$ and $\varepsilon\geq \varepsilon_0>0$ such that
\begin{align} 
    \varepsilon
    &\leq{\left[1+4(2|g|_\alpha+|f|_{\infty}){\sf adiam}(D)\vee 1\right]}^{-\frac{2}{\alpha}}\gamma^{\frac{2}{\alpha}}d^{-1}. \label{eq:epsilon}\\ 
    \intertext{Also, choose}
   K &:= \left\lceil {\sf diam}(D)\left(\frac{8\left(2|g|_\alpha+\frac{{\rm diam}(D)^2|f|_\alpha+2{\rm diam}(D)|f|_\infty}{d}\right)+1}{\gamma}\right)^{1/\alpha}\right\rceil. \label{eq:K}
\end{align}
Then
 \begin{equation}\label{eq:confidence}
    \mathbb{P}\left( \sup_{x\in D} \left | u(x)-u_M^{ N}(x)\right| \geq \gamma \right) \leq \eta
\end{equation}
whenever we choose
\begin{equation}\label{eq:N_N}
    N\geq \frac{16\left\{d\left[\lceil M/\alpha\rceil\log(2+|\widetilde{r}|_1)+\log(K)\right]+\log(\frac{2}{\eta})\right\}\left[|g|_\infty+\frac{M}{d}{\sf diam}(D)^2|f|_\infty\right]^2}{9\gamma^2}
\end{equation}
and 
{\begin{equation}\label{eq:Mii}
        M\geq \frac{\left[\log(4/\gamma)+\log(4|g|_\infty+\frac{2}{d}{\sf diam}(D)^2|f|_\infty)\right]4{\sf diam}(D)^2}{\beta^2\varepsilon^2}.
    \end{equation} 
Furthermore, if $D$ is $\delta$-defective convex, then $M$ can be chosen as 
\begin{equation}\label{eq:Mi}
    M\geq \frac{4d\left[\log\left(\frac{4}{\gamma}\sqrt{\frac{{\sf diam}(D)}{\varepsilon}}\right)+\log(4|g|_\infty+\frac{2}{d}{\sf diam}(D)^2|f|_\infty)\right]}{\beta^2(1-\delta)}.
\end{equation}
}

\medskip
\noindent{Proof in \Cref{pf:4}.}
\end{coro}

\begin{rem}\label{r:choice}
By some simple computations, $\varepsilon,\varepsilon_0$, $M$,  and $N$ from \Cref{coro:noKepsilon}, exhibit the following asymptotic behaviors:
\begin{enumerate}
    \item[] $\varepsilon,\varepsilon_0 \in \mathcal{O}(\gamma^{\frac{2}{\alpha}}d^{-1})$, 
\end{enumerate}
{and if $D$ is $\delta$-defective convex then
\begin{equation*}
M\in \mathcal{O}\left(\frac{d\log(d/\gamma)}{\beta^2(1-\delta)}\right), \quad N\in \mathcal{O}\left(\frac{\log^2(d/\gamma)\left[d^2\log(d/\gamma)+\beta^2(1-\delta)\log(1/\eta)\right]}{\beta^4\gamma^2(1-\delta)^2}\right),
\end{equation*}
whilst if $D$ { satisfies the uniform exterior ball condition} then
\begin{equation*}
M\in \mathcal{O}\left(\frac{d^2\log(1/\gamma)}{\beta^2\gamma^{4/\alpha}}\right), \quad N\in \mathcal{O}\left(\frac{d^2\log(1/\gamma)^2\left[d^2\log^2(1/\gamma)+\beta^2\gamma^{2/\alpha}\log(1/\eta)\right]}{\beta^4\gamma^{2+4/\alpha}}\right).
\end{equation*}
Here, the Landau symbols tacitly depend on (the regularity of) $f,g, {\sf diam}(D), {\sf adiam(D)},\mbox{ and } \widetilde{r}$.
In particular, only in terms of the dimension $d$, if the domain is $\delta$-convex then $M\in \mathcal{O}(d\log(d))$ and $N\in \mathcal{O}(d^2\log^3(d))$, whilst merely under the uniform exterior ball condition, $M\in \mathcal{O}(d^2)$ and $N\in \mathcal{O}(d^4)$.}
\end{rem}

\subsection{On regular extensions of the boundary data inside the domain}
\label{ss:extension}

Recall that one assumption of the main results in the previous subsections (see e.g. \Cref{thm:main}) is that the boundary data $g$ can be extended as a regular function (H\" older or $C^2$) defined on the entire domain $\overline{D}$. 
This is required by the fact that the data needs to be evaluated at the location where the WoS chain is stopped, see \eqref{eq:MCestimator}, and such stopped position lies in principle in the interior of the domain $D$. 
However, usually in practice, $g$ is measured (hence known) merely at the boundary $\partial D$. 
With this issue in mind, in this subsection we address the problem of extending $g$ regularly from $\partial D$ to $\overline{D}$, in a constructive way which is also DNN-compatible.

We take $D\subset \mathbb{R}^d$ to be a set of class $C^k$, $k=3$ or $k=2$, hence (see \cite[sec. 14.6]{GTbook}) there exists a neighbourhood $D_{\epsilon_0}:=\{x\in D;\textrm{dist}(x,\partial D)<\epsilon_0\}$ of $\partial D$ such that the restriction of the distance function $r:\Oe\to\mathbb{R}_+$ is of class $C^k$, and the nearest point projection $\po:\Oe\to \partial D$ is of class $C^{k-1}$.
We have:
\begin{lem} \label{lem:extension}Let $D\subset \R^d$ to be a set of class $C^2$ hence for any point $x\in\partial D$ there exist a function $\phi_x:\mathbb{R}^{d-1}\to\R$ of class $C^2$ and a radius $r_x>0$ such that 
$$
D\cap B(x,r_x)=\{y=(y_1,\dots,y_d)\in B(x,r_x); y_d<\phi_x(y_1,\dots,y_{d-1}).
$$
We denote by $M:=\sup_{\stackrel{i,j,k\in\{1,\dots,d\}}{x\in\partial D}}\big|\frac{\partial^3 \phi_x}{\partial y_i\partial y_j\partial y_k}(x)\big|$. Furthermore, denoting by $k_1(x)\le k_2(x)\le \dots\le k_{d-1}(x)$ the ordered principal curvatures of $\partial D$ let us take $\eps_0:=\min_{x\in\partial D}k_{d-1}^{-1}(x)$. 
Take $\psi\in C^\infty_c([0,\infty),\R)$ to be such that $\psi\equiv 1$ on $[0,1]$ and $\psi\equiv 0$ on $[3,\infty)$, $|\psi|_{\infty}|,|\psi'|_{\infty}|,|\psi''|_{\infty}\le 1$.
We define the extension $G$ in $\overline{D}$ of the $\alpha$-H\"older function $g$ given on the boundary $\partial D$, for $\alpha\in(0,1]$ as:
\begin{equation}\label{def:Gext}
 G: \overline{D}\to\mathbb{R}, \quad G(x):=\psi\left(\frac{1}{\eps_0}r(x))\right)g(\po(x)),\; x\in \overline{D}.
\end{equation}
Then $G$ is $\alpha$-H\"older on $\overline{D}$ and 
$
|G|_\alpha\le |\nabla\po|_{\infty}^\alpha |g|_\alpha+|g|_{\infty}\eps_0^{-1} |\sf{diam}(D)|^{1-\alpha}.
$

\medskip
\noindent{Proof in \Cref{pf:6}.}
\end{lem}

\begin{rem} One can easily provide a non-constructive $\alpha$-H\"older extension $\widetilde{G}$ on $\overline{D}$ of the $\alpha$-H\"older boundary data $g$ given on $\partial D$ by setting $\widetilde G(x):=\inf\{ g(y)+|g|_\alpha|x-y|^\alpha, y\in\partial D\}, x\in \overline{D}$.
\end{rem}

\section{DNN counterpart of the main results}\label{s:NN}

Let $\sigma:\mathbb{R} \to \mathbb{R}$ be the rectified linear unit (ReLU) activation function, 
that is $\sigma(x) := \max\{0,x\}$, $x\in \mathbb{R}$.
Let $(d_i)_{i=0,\ldots,d}$ be a sequence of positive integers. 
Let $A^i \in \mathbb{R}^{d_{i} \times d_{i-1}}$ and $b^i\in \mathbb{R}^{d_i}$, 
$i=1,\ldots,L$, and set $W^i(x):=A^ix+b^i, x\in\mathbb{R}^d$.
We define the realization of the DNN $\mathbb{R}^{d_0} \ni x\mapsto \phi(x)$ by 
\begin{equation}\label{eq:def_relu_nn}
\mathbb{R}^{d_0}\ni x \mapsto \phi(x):= W^L\circ\sigma\circ W^{L-1}\cdots\circ\sigma\circ W^1(x)\in \mathbb{R}^{d_L}, \quad x\in \mathbb{R}^{d_0}, 
\end{equation}
where $\mathbb{R}^{d} \ni x \mapsto \sigma (x):= (\sigma(x_1),\ldots,\sigma(x_d)) $, $d\in\mathbb{N}$,
is defined coordinatewise.
The \emph{weights} of the ReLU DNN $\phi$ are the entries of $(A^i, b^i)_{i=1,\ldots,L}$. 
The {\it size} of $\phi$ denoted by ${\rm size}(\phi)$ is the number of 
non-zero weights. 
The {\it width} of $\phi$ is defined by ${\rm width}(\phi^L) = \max\{d_0,\ldots,d_L\}$
and $L$ is the depth of $\phi$ denoted by $\mathcal{L}(\phi)$.
In the sequel, we only consider DNNs with ReLU activation function.

For the reader's convenience, before we proceed to the main result of this section (see \Cref{thm:mainNN} below), we present first several technical lemmas following \cite{EDGB_2019} and \cite{yarotsky2017error}, as well as some of their consequences; all these preparatory results are meant to provide a clear and systematic way of quantifying the size of the DNN which is constructed in the forthcoming main result, namely \Cref{thm:mainNN}.

The following lemma is~\cite[Proposition~3]{yarotsky2017error}.

\begin{lem}\label{lem:NN_prod_scalars}
For every $c>0$ and $\delta\in (0,1)$, there exists a DNN $\Pi_\delta^c$ 
such that 
\begin{equation*}
 \sup_{a,b\in [-c,c]}|ab - \Pi_\delta^c(a,b)|
 \leq \delta \quad \mbox {and}\quad {\rm size}(\Pi_\delta^c) = \mathcal{O}(\lceil\log (\delta^{-1})+\log(c)\rceil). 
\end{equation*}
\end{lem}

Now, we recall Lemma II.6 from \cite{EDGB_2019}:
\begin{lem}\label{lem:additions_NNs}
Let $\phi_{i},\;i=1,\ldots, n$,  be ReLU DNNs with the same input dimension $d_0\in\mathbb{N}$ and the same depth $\mathcal{L}:=\mathcal{L}(\phi_i),1\leq i\leq n$.
Let $a_i$, $i=1,\ldots,n$, be scalars.
Then there exists a ReLU DNN $\phi$ such that
\begin{itemize}
    \item[i)] $\phi(x) =\sum_{i=1}^n a_i \phi_i(x)$ for every $x\in \mathbb{R}^{d_0}$,
    \item[ii)] $\mathcal{L}(\phi)=\mathcal{L}$,
    \item[iii)] $\mathcal{W}(\phi)\leq \sum\limits_{1\leq i\leq n}\mathcal{W}(\phi_i)$,
    \item[iv)] ${\rm size}(\phi)\leq \sum\limits_{1\leq i\leq n}{\rm size}(\phi_i)$.
\end{itemize}
\end{lem}

The following lemma is taken from ~\cite[Lemma~II.3]{EDGB_2019}, with the mention that the last assertion iv) brings some improvement which is relevant to our purpose; it is immediately entailed by the proof of the same ~\cite[Lemma~II.3]{EDGB_2019}, so we skip its justification.
\begin{lem}\label{lem:composition_NNs}
Let $\phi_1:\mathbb{R}^{d_1}\rightarrow\mathbb{R}^{d_2}$ and $\phi_2:\mathbb{R}^{d_3}\rightarrow\mathbb{R}^{d_1}$ be two ReLU DNNs.
Then there exists a ReLU DNN $\phi:\mathbb{R}^{d_3}\rightarrow\mathbb{R}^{d_2}$ such that 
\begin{enumerate}
    \item[i)] $\phi(x) = \phi_1(\phi_2(x))$ for every $x\in\mathbb{R}^{d_0^2}$, 
    \item[iii)]$\mathcal{L}(\phi)=\mathcal{L}(\phi_1)+\mathcal{L}(\phi_2)$, 
    \item[iii)] $\mathcal{W}(\phi)\leq \max(\mathcal{W}(\phi_1),\mathcal{W}(\phi_2),2d_1)$, 
    \item[iv)] ${\rm size}(\phi ) \leq  \min\left({\rm size}(\phi_1 ) + {\rm size}(\phi_2 ) + d_1[\mathcal{W}(\phi_1)+\mathcal{W}(\phi_2)], \; 2{\rm size}(\phi_1 ) + 2{\rm size}(\phi_2 )\right)$. 
\end{enumerate}
\end{lem}

The next lemma is essentially \cite[Lemma II.4]{EDGB_2019}. As in the case of the previous lemma, assertion iv) comes with a slight modification of the original result, which can be immediately deduced from the proof of \cite[Lemma II.4]{EDGB_2019}.
\begin{lem}\label{lem:augment}
Let $\phi:\mathbb{R}^{d_0}\rightarrow\mathbb{R}^{d_1}$ be a ReLU DNN such that $\mathcal{L}(\phi)<L$. 
Then there exists a second ReLU DNN $\widetilde{\phi}:\mathbb{R}^{d_0}\rightarrow\mathbb{R}^{d_1}$ such that
\begin{enumerate}
    \item[i)] $\phi(x)=\widetilde{\phi}(x)$ for all $x\in \mathbb{R}^{d_0}$,
    \item[ii)] $\mathcal{L}(\widetilde{\phi})=L$,
    \item[iii)] $\mathcal{W}(\widetilde{\phi})=\max(2d_1,\mathcal{W}(\phi))$,
    \item[iv)] ${\rm size}(\widetilde{\phi})\leq \min\left({\rm size}(\phi)+d_1\mathcal{W}(\phi), \; 2{\rm size}(\phi)\right)+2d_1(L-\mathcal{L}(\phi))$.
\end{enumerate}

\end{lem}

As a direct consequence of \Cref{lem:NN_prod_scalars}, \Cref{lem:composition_NNs}, \Cref{lem:augment}, and \cite[Lemma II.5]{EDGB_2019}, one gets the following approximation result for products of scalar ReLU DNN:
\begin{coro} \label{coro:product_NNs}
Let $\phi_1,\phi_2:\mathbb{R}^{d}\rightarrow\mathbb{R}$ be two ReLU DNNs, $D\subset \mathbb{R}^d$ be a bounded subset, and let $\Pi:=\Pi_{\epsilon_p}^c$ be given by \Cref{lem:NN_prod_scalars} for $c:=\max\left(\sup\limits_{x\in D}\phi_1(x),\sup\limits_{x\in D}\phi_2(x)\right)$ and $\epsilon_p>0$.
Then there exists a ReLU DNN $\phi:\mathbb{R}^{d}\rightarrow\mathbb{R}$ such that 
\begin{enumerate}
    \item[i)] $\phi(x) =\Pi( \phi_1(x),\phi_2(x))$ for every $x\in\mathbb{R}^d$, 
    \item[ii)] 
        $
         \sup\limits_{x\in D}|\phi_1(x)\phi_2(x) - \phi(x)|
         \leq \epsilon_p,
        $
    \item[iii)] ${\rm size}(\phi ) \leq 4{\rm size}(\phi_1)+4{\rm size}(\phi_2)+\mathcal{O}(\lceil\log (\epsilon_p^{-1})+\log(c)\rceil)$.
\end{enumerate}
\end{coro}

The following two lemmas are going to be employed later in order to quantify the size of one generic step of the WoS chain given by \eqref{wos0}-\eqref{wosn}, regarded as an action of a ReLU DNN. 
\begin{lem} \label{lem:sizeonestep}
Let $\phi:\mathbb{R}^{d}\rightarrow\mathbb{R}$ be a ReLU DNN and $v\in \mathbb{R}^d$ be a vector. 
Then there exists a ReLU DNN $\phi_v:\mathbb{R}^{d}\rightarrow\mathbb{R}$ such that
\begin{enumerate}
    \item[i)] $\phi_v(x)=x+\phi(x)v$ for all $x\in \mathbb{R}^d$, 
    \item[ii)] $\mathcal{L}(\phi_v)=\mathcal{L}(\phi)+1$, 
    \item[iii)] $\mathcal{W}(\phi_v)\leq 2d+\max(d,\mathcal{W}(\phi))$, 
    \item[iv)] ${\rm size}(\phi_v)\leq 2{\rm size}(\phi) +2d[\mathcal{L}(\phi)+2]$.
\end{enumerate}

\medskip
\noindent{Proof in \Cref{pf:6}.}
\end{lem}

The following result is easily deduced by employing recursively \Cref{lem:sizeonestep} and \Cref{lem:composition_NNs}, so we omit its proof.
\begin{coro} \label{coro:sizeksteps}
Let $\phi:\mathbb{R}^{d}\rightarrow\mathbb{R}$ be a ReLU DNN and $v_k\in \mathbb{R}^d,k\geq 1$ be a sequence of vectors. 
Then there exist ReLU DNNs $\theta_k:\mathbb{R}^{d}\rightarrow\mathbb{R}^d,k\geq 0$,  such that for every $k\geq 0$
\begin{enumerate}
    \item[i)] $\theta_{k+1}(x)=\phi_{v_{k+1}}(\theta_{k}(x))$ and $\theta_0(x)=x$ for all $x\in \mathbb{R}^d$, where $\phi_{v_{k}}$ is the one constructed in \Cref{lem:sizeonestep},
    \item[ii)] $\mathcal{L}(\theta_{k+1})=(k+1)(\mathcal{L}(\phi)+1)+1$,
    \item[iii)] $\mathcal{W}(\theta_{k+1})\leq 2d+\max(d,\mathcal{W}(\phi))$,
    \item[iv)] ${\rm size}(\theta_{k+1})\leq 2d(k+1)[4d+\mathcal{W}(\phi)+\mathcal{L}(\phi)+2]+d+2(k+1){\rm size}(\phi)$.
\end{enumerate}
\end{coro}

We end this first paragraph by a ReLU DNN extension of a H\"older continuous boundary data $g$ to the entire domain $\overline{D}$.
\begin{coro} \label{cor:NNforGext}
Let $ D\subset \R^d$ to be a set of class $C^2$  and $g$ be $\alpha$-H\"older on $\partial D$. 
For $\eps_0>0$ and $\psi\in C^\infty_c([0,\infty),\R)$ as defined in Lemma~\ref{lem:extension} we assume that
for every $\delta_r,\delta_p, \delta_\psi, \delta_g\in (0,1)$, there exist ReLU DNNs $\phi_{r}$, $\phi_{\pi}$,$\phi_{\psi}$ and $\phi_{g}$ such that
\[
|r -  \phi_{r}|_{\infty}\leq \delta_r,\quad
|\po -  \phi_{\pi}|_{\infty}\leq \delta_\pi,\quad
| \psi -  \phi_{\psi}|_{\infty}\leq \delta_\psi,\quad
|g -  \phi_{g}|_{\infty}\leq \delta_g.
\]
%
%
With $\varepsilon_0$ the one given in \Cref{lem:extension}, set 
\[
\overline{\delta}:=2\left(3\delta_\psi+\frac{\delta_d}{\eps_0}\right)|g|_\infty+ 2\left(3\delta_g+|\nabla g|_{\infty}\delta_\pi \right)(\delta_\psi+1)\in \mathcal{O}(\delta_\psi+\delta_r+\delta_g+\delta_\pi).
\]
If $G$ is the $\alpha$-H\"older extension in $D$ of the boundary data $g$ given by \eqref{def:Gext}, 
then there exists a ReLU DNN $\phi_{G}$
such that
\begin{enumerate}
    \item[i)] $|G-\phi_{G}|_{\infty}\leq\overline{\delta},$\label{G:approx}
    \item[ii)] ${\rm size}(\phi_{G}) \leq 2{\rm size}(\phi_{\psi})+2{\rm size}(\phi_{r})+2{\rm size}(\phi_{g})+2{\rm size}(\phi_{\pi})+\mathcal{O}(\lceil\log (\overline{\delta}^{-1})+\log(|g|_\infty)\rceil)\nonumber$.
\end{enumerate}

\medskip
\noindent{Proof in \Cref{pf:6}.}
\end{coro}

\subsection{DNN approximations for solutions to problem \eqref{eq:main1}}

We are now ready to present the DNN byproduct of \Cref{thm:main}, in fact of \Cref{coro:noKepsilon}.
First, let us state that the $\widetilde{r}$-WoS chain given by \eqref{wos0}-\eqref{wosn} renders a ReLU DNN as soon as $\widetilde{r}$ is a ReLU DNN; this follows from a simple corroboration of \Cref{lem:sizeonestep} and \Cref{coro:sizeksteps}, so we skip its formal proof.
\begin{coro} \label{lem:sizeWoS}
Suppose that $\widetilde{r}$ is a ReLU DNN on the bounded set $D\subset \mathbb{R}^d$ such that $0\leq \widetilde{r}\leq r$, where recall that $r$ specified by \eqref{eq:rdistance} is the distance function to the boundary of $D$.
Further, let $M\geq 0$ and $\left(X_M^x,x\in D\right)$ be the $\widetilde{r}$-WoS chain at step $M$ given by $\eqref{wos0}-\eqref{wosn}$.
Then for each $\omega \in \Omega$ there exists a ReLU DNN defined on $D$ and denoted by $\mathbb{X}_M^\omega(\cdot)$ such that
\begin{enumerate}
    \item[i)] $\mathbb{X}_M^\omega(x)=X_M^x(\omega)$ for all $x\in D$,
    \item[ii)] $\mathcal{L}(\mathbb{X}_M^\omega)=M(\mathcal{L}(\widetilde{r})+1)+1$,
    \item[iii)] $\mathcal{W}(\mathbb{X}_M^\omega)\leq 2d+\max(d,\mathcal{W}(\widetilde{r}))$,
    \item[iv)] ${\rm size}(\mathbb{X}_M^\omega)\leq 2dM[4d+\mathcal{W}(\widetilde{r})+\mathcal{L}(\widetilde{r})+2]+d+2M{\rm size}(\widetilde{r})$.
\end{enumerate}
\end{coro}

The main result of this section is the following, proving that ReLU DNNs can approximate the solution $u$ to problem \eqref{eq:main1} without the curse of dimensions.
\begin{thm} \label{thm:mainNN} 
The statement requires a detailed context, so let us label the assumption and the conclusion separately.

\noindent
\textit{\textbf{Assumption:}}
Let $D\subset \mathbb{R}^d$ be a bounded open set { satisfying the uniform exterior ball condition}, $f$ and $g$ be $\alpha$-H\"older functions on $\overline{D}$ for some $\alpha \in [0,1]$, and $u$ be the solution to \eqref{eq:main1}, as in \Cref{thm:representation}.
Let $\phi_f:D\rightarrow\mathbb{R},\phi_g:\overline{D}\rightarrow\mathbb{R},\phi_r:D\rightarrow \mathbb{R}$ be ReLU DNNs such that
\begin{equation*}
    |f-\phi_f|_\infty\leq \epsilon_f\leq |f|_\infty, \quad |g-\phi_g|_\infty\leq \epsilon_g, \quad  |r-\phi_r|_\infty\leq \epsilon_r,
\end{equation*}
and set 
$
    \widetilde{r}(x):=(\phi_r(x)-\epsilon_r)^{+}, \, x\in D.$ 
Also, let $\Pi:=\Pi_{\epsilon_p}^c$ be the ReLU DNN given by \Cref{lem:NN_prod_scalars}. 

Further, let $\gamma>0$ be a prescribed error, $0<\eta<1$ be a prescribed confidence, and consider the following assumptions on the parameters $\epsilon_f,\epsilon_g, \epsilon_r, \epsilon_p \text{ and } c $:
\begin{enumerate}
    \item[a.1)] $\epsilon_g\leq \gamma/6$,
    \item[a.2)] $\epsilon_f\leq\frac{d\gamma}{6M{\rm diam}(D)^2}$,
    \item[a.3)] $\epsilon_r<\varepsilon_0:=\frac{1}{3}\left[(8|g|_\alpha+|f|_{\infty}){\sf adiam}(D)\vee 1\right]^{-\frac{2}{\alpha}}(\gamma/2)^{\frac{2}{\alpha}}d^{-1}$, so that, by \Cref{rem:adistance}, $\widetilde{r}$ is a $(\beta,\varepsilon_0)$-distance if we choose $\beta=1/3$,
    \item[a.4)] $\epsilon_p=\frac{\gamma d}{6M(1+2|f|_\infty)}$ and $c=\max({\rm diam}(D),2|f|_\infty)$,
\end{enumerate}
where $M\geq 1$ is specified below. 

Further, consider the iid pairs $((X^{\cdot,i}_{k})_{k\geq 1},Y^i), i\geq 1$  on $(\Omega,\mathcal{F},\mathbb{P})$, as in the beginning of \Cref{ss:meantailestimates}, and
\begin{equation}\label{eq:NNMC}
    \widetilde{u}_{M}^{N}(x):=\frac{1}{N}\sum_{i=1}^N\left[ \phi_g(X^{x,i}_{M})+\frac{1}{d}\sum\limits_{k=1}^M \Pi\left(\Pi\left(\widetilde{r}(X^{x,i}_{k-1}),\widetilde{r}(X^{x,i}_{k-1})\right),\phi_f\left(X^{x,i}_{k-1}+\widetilde{r}(X^{x,i}_{k-1})Y^i\right)\right)\right], x\in D.
\end{equation}

{
Let us choose
\begin{align}\label{eq:N}
    &N\geq \frac{64\left\{d\left[\lceil M/\alpha\rceil\log(2+|\phi_r|_1)+\log(K)\right]+\log(\frac{2}{\eta})\right\}\left[|g|_\infty+\frac{M}{d}{\sf diam}(D)|f|_\infty\right]^2}{9\gamma^2},\\
    &M\geq \frac{36\left[\log(8/\gamma)+\log(4|g|_\infty+\frac{2}{d}{\sf diam}(D)^2|f|_\infty)\right]{\sf diam}(D)^2}{\varepsilon_0^2},
\end{align}
where $K := \left\lceil {\sf diam}(D)\left(\frac{16\left(2|g|_\alpha+\frac{{\rm diam}(D)^2|f|_\alpha+2{\rm diam}(D)|f|_\infty}{d}\right)+2}{\gamma}\right)^{1/\alpha}\right\rceil$.

Furthermore, if $D$ is $\delta$-defective convex then $M$ can be chose such that
\begin{equation*}
    M\geq 36 \frac{d\log\left(\frac{8}{\gamma}\sqrt{\frac{{\sf diam}(D)}{\varepsilon_0}}\right)
    +\log(4|g|_\infty+\frac{2}{d}{\sf diam}(D)^2|f|_\infty)}{(1-\delta)}.
\end{equation*}
}

\noindent
\textit{\textbf{Conclusion:}}
Under the above assumption and keeping the same notations, there exits a measurable function $\mathbb{U}_M^N:\Omega\times D\rightarrow \mathbb{R}$ such that
\begin{itemize}
\item[c.1)] $\mathbb{U}_M^N(\omega,\cdot)$ is a ReLU DNN for each $\omega\in \Omega$, $\mathbb{U}_M^N(\cdot,x)=\widetilde{u}_M^N(x), \quad x\in D,$ and $$
\mathbb{P}\left(\sup\limits_{x\in D} \left| u(x)-\mathbb{U}_M^N(\cdot,x)\right| \geq \gamma \right) \leq \eta.$$
\item[c.2)] For each $\omega\in \Omega$ we have that
\[
{\rm size}(\mathbb{U}_M^N(\omega,\cdot))\in \mathcal{O}\left(MN\left[dM\max(d,\mathcal{W}(\phi_r),\mathcal{L}(\phi_r))+M{\rm size}(\phi_r)+{\rm size}(\phi_f)+\left\lceil\log \left(\frac{1}{\gamma d}\right)\right\rceil\right]\right),
\]
where the tacit constant depends only on $\max({\rm diam}(D),|f|_\infty)$.
\end{itemize}
{
In particular, using \eqref{eq:N}, $\mathbb{U}_M^N(\omega,\cdot)$ can be realized such that
\[
{\rm size}(\mathbb{U}_M^{ N}(\omega,\cdot))\in \mathcal{O}\left(d^7\gamma^{-16/\alpha-4}\log^4\left(\frac{1}{\gamma}\right)\left[d^3\gamma^{-4/\alpha}\log\left(\frac{1}{\gamma}\right)+\log\left(\frac{1}{\eta}\right)\right]\log(2+|\phi_r|_1){\rm S}\right), 
\]
where 
\[
{\rm S}:=\left[\max(d,\mathcal{W}(\phi_r),\mathcal{L}(\phi_r))+{\rm size}(\phi_r)+{\rm size}(\phi_g)+{\rm size}(\phi_f)\right]
\]
and the tacit constant depends on $|g|_\alpha,|g|_\infty,|f|_\infty,{\rm diam}(D),{\rm adiam}(D),\delta,\alpha$.

Furthermore, if $D$ is $\delta$-defective convex then
if $a.5)$ holds then
\[
{\rm size}(\mathbb{U}_M^{ N}(\omega,\cdot))\in \mathcal{O}\left(\frac{d^3}{\gamma^2}\log^4\left(\frac{d}{\gamma}\right)\left[d^2\log\left(\frac{d}{\gamma}\right)+\log\left(\frac{1}{\eta}\right)\right]\log(2+|\phi_r|_1){\rm S}\right).
\]
}

\medskip
\noindent{Proof in \Cref{pf:6}.}
\end{thm}

We end the exposition of the main results with the remark that it is sufficient to prescribe the Dirichlet data $g$ merely on $\partial D$ (not necessarily extended to $\overline{D}$), as expressed by the following direct consequence of \Cref{cor:NNforGext}.
\begin{coro} If the open and bounded set $D$ is of class $C^2$, and if $g$ is given merely on $\partial D$ and it is $\alpha$-H\"older there, then $g$ can be constructively extended to an $\alpha$-H\"older function on $\overline{D}$.
Furthermore, a ReLU DNN approximation $\phi_g$ can be constructed as in \Cref{cor:NNforGext}, so \Cref{thm:main} and \Cref{thm:mainNN} fully apply.
\end{coro}

\section{Proofs of the main results}
\label{s:proofs}

\subsection{Proofs for \Cref{ss:representation}}
\label{pf:1}

\begin{proof}[\textit{\textbf{Proof of \Cref{lem:K}}}]

Let $v(x):=u(Ax)$. Then $\partial_k v(x)=\sum_{j,l=1}^n \partial_j u(Ax) \partial_k(A_{jl}x_l)=\sum_{j=1}^n \partial_j u(Ax) A_{jk}$ and

\begin{equation*}
\sum_{k=1}^n \partial_k\partial_k v(x)=\sum_{i,j,k,l=1}^n \partial_i\partial_j u(Ax)A_{jk}\partial_k (A_{il}x_l)=\sum_{i,j,k=1}^n \partial_i\partial_j u(Ax) A_{jk}A_{ik}. 
\end{equation*}
Thus we need to determine $A$ such that $AA^T=K$.
We know that there exists a rotation matrix $R\in O(3)$ (hence $RR^T=Id$) such that $RKR^t=\textrm{diag}(\lambda_1,\dots,\lambda_n)$. 
Since $K$ is positive definite we have $\lambda_i>0$ for $i\in\{1,\dots,n\}$. 
Then we have:

\begin{equation}\label{rel:RARS}
RAR^TRA^TR^T=RKR^T=\textrm{diag}(\lambda_1,\dots,\lambda_n).
\end{equation}

We denote $B:=RAR^T$ and observe that \eqref{rel:RARS} can be rewritten as $BB^T=\textrm{diag}(\lambda_1,\dots,\lambda_n)$. We can now take $B:=\textrm{diag}(\sqrt{\lambda_1},\dots,\sqrt{\lambda_n})$. Thus $A:=R^T\textrm{diag}(\sqrt{\lambda_1},\dots,\sqrt{\lambda_n})R$.
\end{proof}

\begin{proof}[\textit{\textbf{Proof of \Cref{coro:v_inf}}}]
The second assertion follows directly from \Cref{thm:representation} and from classical regularity theory for Poisson equation, so let us prove the first one.
To this end, note that without loss of generality we may assume that $x=0\in D$, so by Ito's formula we get that 
$(|B^0(t)|^2-dt)_{t\geq 0}$
is a martingale and
\begin{equation*}
    v(0)=\mathbb{E}\left\{\tau^0_{D^c}\right\}=\frac{1}{d}\mathbb{E}\left\{|B^0(\tau_{D^c})|^2\right\}\leq {\sf diam}(D)^2/d.
\end{equation*}
\end{proof}

\begin{proof}[\textit{\textbf{Proof of \Cref{coro:representation}}}]

\noindent
i). Using \Cref{rem:beta0distance} we have that,
\begin{align*}
  \mathbb{E}\left\{\int_0^{\tau^x_{D^c}}f(B^x(t))\;dt\right\}&=\sum\limits_{n\geq 0}\mathbb{E}\left\{\int_{\tau^{x}_{n}}^{\tau^{x}_{n+1}}f(B^x(t))\;dt\right\}\\
  &=\sum\limits_{n\geq 1}\mathbb{E}\left\{\widetilde{r}^2(B^x_{\tau^{x}_{n-1}})K_0f(B^x_{\tau^{x}_{n-1}}+\widetilde{r}(B^x_{\tau^{x}_{n-1}})\cdot)\right\}\\
  &=\mathbb{E}\left\{\sum\limits_{k\geq 1}\widetilde{r}^2(X^{x}_{k-1})K_0f(X^{x}_{k-1}+\widetilde{r}(X^{x}_{k-1})\cdot)\right\}, \quad x\in D,
\end{align*}
where the second equality follows by the strong Markov and the scaling properties of Brownian motion, whilst the last equality follows from the fact that the law of $(X^x_n)_{n\geq 0}$ under $\widetilde{\mathbb{P}}$ and the law of $(B^x_{\tau^x_{n}})_{n\geq 0}$ under $\mathbb{P}^0$ are the same.
Therefore, the statement follows by \Cref{thm:representation} and also \Cref{rem:beta0distance}.

\medskip
\noindent
ii). The claim follows from the fact that  $\mu(A)=w(0)$, where $w$ solves
\begin{equation*}
\left\{\begin{array}{ll} -\frac{1}{2}\Delta w=d1_A &\,\textrm{ in } B(0,1)\\[1mm]
\phantom{\Delta}w=0 &\,\textrm{ on }\partial B(0,1)
\end{array}\right.
\end{equation*}
for all $A\in \mathcal{B}(B(0,1))$.

\medskip
\noindent
iii). We use conditional expectation, namely for every $x\in D$
\begin{align*}
   \mathbb{E}&\left\{\sum\limits_{k\geq 1}\widetilde{r}^2(X^{x}_{k-1})f(X^{x}_{k-1}+\widetilde{r}(X^{x}_{k-1})Y)\right\}=\sum\limits_{k\geq 1}\mathbb{E}\left\{\widetilde{r}^2(X^{x}_{k-1})f(X^{x}_{k-1}+\widetilde{r}(X^{x}_{k-1})Y)\right\}\\ &=\sum\limits_{k\geq 1}\mathbb{E}\left\{\mathbb{E}\left[\widetilde{r}^2(X^{x}_{k-1})f(X^{x}_{k-1}+\widetilde{r}(X^{x}_{k-1})Y)|\,  X_{k-1}^{x}\right]\right\}=
   d\sum\limits_{k\geq 1}\mathbb{E}\left\{\widetilde{r}^2(X^{x}_{k-1})K_0f(X^{x}_{k-1}+\widetilde{r}(X^{x}_{k-1})\cdot)\right\},
\end{align*}
where for the last equality we used that $Y$ has distribution $\mu$ and is independent of $X_{k-1}^x$.
\end{proof}

\begin{proof}[\textit{\textbf{Proof of \Cref{prop:exitannulus}}}]
Recall that $\mathbb{E}\{\tau^x_{D^c}\}=v(x),x\in D$, where $v$ is given by \eqref{def:v}.
The idea is to explicitly solve for $v(x)=-2w(|x|)$ in radial form as 
\[
 w''(r)+\frac{d-1}{r}w'(r)=1 \text{ with }w(R_0)=0, w(R_1)=0 
\]
which is explicitly solved as 
\[
 w(r)=\frac{r^2}{2d}+C_1r^{2-d}+C_2\text{ with }C_1=\frac{R_1^2-R_0^2}{2d(R_0^{2-d}-R_1^{2-d})}, C_2=\frac{R_0^d-R_1^d}{2d(R_1^{d-2}-R_0^{d-2})}.
\]
Now, if we start from a point $x$ we can use the intermediate value theorem to obtain first that 
\[
 \mathbbm{E}_{x}\{ \tau_{D^c}\} =-2w(|x|)\le 2d(x,\partial A(R_0,R_1))|w'(r_0)|  
\]
for some point $r_0\in[R_0,R_1]$. 
Therefore our task now is to estimate the above derivative, which we can compute explicitly as 
\[
w'(r)=\frac{r}{d}+\frac{2-d}{2d}\frac{R_1^2-R_0^2}{R_0^{2-d}-R_1^{2-d}}r^{1-d}.  
\]
To estimate this in a transparent way we set $R_1=\rho R_0$ and $r=tR_0$ for $1\le t\le \rho$. In these new notations we have 
\[
v'(r)=\frac{R_0}{d}\left(t-\frac{(2-d)(\rho^2-1)}{2(\rho^{2-d}-1)}t^{1-d} \right).
\]
Now, it is an elementary matter that for two continuous functions $f,g:[a,b]\to\R$ which are differentiable on $(a,b)$, with $g(b)\neq g(a)$ and $g'$ non-vanishing, we can find a point $\xi\in(a,b)$ such that 
\[
\frac{f(b)-f(a)}{g(b)-g(a)}=\frac{f'(\xi)}{g'(\xi)}.
\]
Using this fact for $a=1,b=\rho$, $f(x)=x^2$ and $g(x)=x^{2-d}$ we argue that for some point $\xi\in(1,\rho)$ we have
\[
w'(r)=\frac{R_0 t}{d}\left(1-(\xi/t)^{d} \right).
\]
As a function of $t\in [1,\rho]$ the above function is increasing, thus we have 
\begin{equation}\label{e:i:0}
\frac{R_0}{d}\left(1-\xi^{n} \right)\le w'(r)\le \frac{R_0 \rho}{d}\left(1-(\xi/\rho)^{d} \right).
\end{equation}
Therefore in order to control $|w'(r)|$ it suffices to control the absolute values of the two bounds above. 
The right hand side bound is easy because for any $x\in(0,1)$, thus we obtain that 
\begin{equation}\label{e:i:1}
R_0\rho\frac{1-(\xi/\rho)^n}{n}\le R_0\rho (1-\xi/\rho)=R_0(\rho-\xi)\le R_1-R_0.  
\end{equation}
The left hand side of \eqref{e:i:0} in absolute value is bounded by 
\begin{equation}\label{e:i:2}
\begin{split}
\left|\frac{R_0(1-\xi^d)}{d}\right|&=\frac{R_0(\xi^d-1)}{d}=\frac{R_0}{d}\left( \frac{(d-2)(\rho^2-1)}{2(1-\rho^{2-d})}-1 \right)=\frac{R_0}{d}\left( \frac{(d-2)(\rho^2-1)\rho^{d-2}}{2(\rho^{d-2}-1)}-1 \right) \\ 
& =\frac{R_0}{d}\left( \frac{(d-2)(\rho^2-1)}{2}+\frac{(d-2)(\rho^2-1)}{2(\rho^{d-2}-1)}-1 \right) \\ 
&\le \begin{cases} 
      \frac{R_0}{3}\left( \frac{(\rho^2-1)}{2}+\frac{\rho-1}{2} \right), & d=3 \\
      \frac{R_0}{2d}(\rho^2-1), & d\ge 4
     \end{cases} 
\\ & \le \frac{R_0(\rho-1)\rho}{2} \\
&\le \frac{(R_1-R_0)R_1}{2R_0}, 
\end{split}
\end{equation}
where we go back to the fact that $\xi^d=\frac{(d-2)(\rho^2-1)}{2(1-\rho)}$ and in the second line we used that $2(\rho^{d-2}-1)\ge (d-2)(\rho^2-1)$ because $\rho\ge1$ for $d\ge4$. Thus combining \eqref{e:i:1} and \eqref{e:i:2} we get that 
\[
 |w'(r)|\le \frac{(R_1-R_0)R_1}{2R_0}
\]
which is our claim. \qedhere 
\end{proof}

\begin{proof}[\textit{\textbf{Proof of \Cref{prop:extime}}}]  
For $x\in D$, we pick a point $y\in\partial D$ such that $d(x,y)=d(x,\partial D)$ and notice that  $\tau_{D^c}\le \tau_{A(a_y,R_{0,y},R_{1,y})^c}$.  
At the same time, since $d(x,\partial A(y,R_{0,y},R_{1,y}))=d(x,\partial D)$, we employ \Cref{prop:exitannulus} to deduce that 
\[
 \mathbbm{E}\left\{ \tau^x_{D^c}\right\} 
 \le \mathbbm{E}\left\{ \tau^x_{A(y,R_{0,y},R_{1,y})^c} \right\} 
 \le 2d(x,\partial D){\sf adiam}(D).  
\]
\end{proof}

\subsection{Proofs for \Cref{ss:nsteps}}
\label{pf:2}

\begin{proof}[\textit{\textbf{Proof of \Cref{prop:nsteps_general}}}]
Let $\varepsilon\geq \varepsilon_0$ and $x\in D$ such that $r(x)\geq \varepsilon$. 
Also, in order to lighten the notation and without loss, let us assume that the walk-on-spheres process $(X_n^x)_{n\geq 0}$ is realized through the trace of the Brownian motion $(B_t^x)_{t\geq 0}$, on the spheres of radius $\widetilde{r}$, namely
\begin{equation*}
    X_n^x=B^x(\tau_n^x), \quad n\geq 0.
\end{equation*}
Recalling that $N_\varepsilon^x$ is given by \eqref{eq:N_eps_x}, we can easily see that
\begin{equation*}
    \sum\limits_{k=0}^{N_\varepsilon^x-1}\left\|B^x(\tau_{k+1}^x)-B^x(\tau_k^x)\right\|^2\geq N_\varepsilon^x \beta^2\varepsilon^2.
\end{equation*}
Consequently, we have
\begin{equation*}
    \mathbb{E}\left\{e^{\lambda\beta^2\varepsilon^2N_\varepsilon^x }\right\}\leq \mathbb{E}\left\{e^{\lambda\sum\limits_{k=0}^{\infty}\left\|B^x(\tau_{k+1}^x)-B^x(\tau_k^x)\right\|^2}\right\}, \quad \lambda>0.
\end{equation*}
The rest of the proof focuses on estimating the right hand side of the above inequality, whilst $\lambda$ will be properly chosen afterwards.
To this end, notice first that by Ito's formula we have
\begin{equation*}
    \left\|B^x(\tau_{k+1}^x)-B^x(\tau_k^x)\right\|^2=2\int_{\tau_k^x}^{\tau_{k+1}^x}\left\langle B^x(t)-B^x(\tau_k^x), \; dB^x(t)\right\rangle+d\left(\tau_{k+1}^x-\tau_k^x\right), \quad k\geq 0.
\end{equation*}
Therefore, if we set
\begin{equation*}
    H(t):=\sum\limits_{k\geq 0}\left[B^x(t)-B^x(\tau_k^x)\right]1_{\left[\tau_k^x,\tau_{k+1}^x\right)}(t), \quad t\geq 0,
\end{equation*}
we have that $H$ is $\mathcal{F}(t)$-adapted, $\left(\int_0^t \left\langle H(s),\;dB^x(s)\right\rangle\right)_{t\geq 0}$ is a continuous local martingale, and
\begin{align*}
    \sum\limits_{k=0}^{N_\varepsilon^x-1}\left\|B^x(\tau_{k+1}^x)-B^x(\tau_k^x)\right\|^2
    &= 2\int_0^{\tau^x_{N_\varepsilon^x}} \left\langle H(s), \;dB^x(s)\right\rangle+d\sum\limits_{k=0}^{N_\varepsilon^x-1} \left(\tau_{k+1}^x-\tau_k^x\right)\\
    &\leq 2\int_0^{\tau^x_{N_\varepsilon^x}} \left\langle H(s), \;dB^x(s)\right\rangle+d\tau^x_{D_\varepsilon}.
\end{align*}
Thus, for $1/p+1/q=1, \;p,q\geq 1,$ we have
\begin{equation*}
    \mathbb{E}\left\{e^{\lambda\beta^2\varepsilon^2N_\varepsilon^x }\right\}
    \leq \mathbb{E}\left\{e^{2\lambda\int_0^{\tau^x_{N_\varepsilon^x}} \left\langle H(s), \;dB^x(s)\right\rangle+\lambda d\tau^x_{D_\varepsilon}}\right\}
    \leq \mathbb{E}\left\{e^{2p\lambda\int_0^{\tau^x_{N_\varepsilon^x}} \left\langle H(s), \;dB^x(s)\right\rangle}\right\}^{1/p}\mathbb{E}\left\{e^{q\lambda d\tau^x_{D_\varepsilon}}\right\}^{1/q}, \quad \lambda>0.
\end{equation*}
Now, using \cite[Theorem 43]{Protter2005} and employing Fatou's lemma, we get
\begin{equation*}
    \mathbb{E}\left\{e^{2p\lambda\int_0^{\tau^x_{N_\varepsilon^x}} \left\langle H(s), \;dB^x(s)\right\rangle}\right\}
    \leq\mathbb{E}\left\{e^{8p^2\lambda^2\int_0^{\tau^x_{N_\varepsilon^x}} \left\| H(s)\right\|^2\;ds}\right\}^{1/2}
    \leq \mathbb{E}\left\{e^{8p^2\lambda^2{{\sf diam}(D)^2}\tau^x_{\partial D}}\right\} ^{1/2}.
\end{equation*}
Consequently,
\begin{equation}\label{eq:p-q}
     \mathbb{E}\left\{e^{\lambda\beta^2\varepsilon^2N_\varepsilon^x }\right\}
     \leq \mathbb{E}\left\{e^{8p^2\lambda^2{{\sf diam}(D)^2}\tau^x_{\partial D}}\right\} ^{1/(2p)}
     \mathbb{E}\left\{e^{q\lambda d\tau^x_{D_\varepsilon}}\right\}^{1/q}.
\end{equation}
Now, let $x\in D$ note that
\begin{equation*}
\text{If } v(z):=1-\frac{\gamma}{2d}\|z-x\|^2, \text{ for } 0<\gamma<\frac{2d}{{\sf diam}(D)^2}, \text{ then }  0\leq v\leq 1 \mbox{ and }  \Delta v(z)\le -\gamma v(z),\quad z\in B(x,{\sf diam}(D)).
\end{equation*}
As a consequence of this, Ito's formula, and the fact that $\nabla v$ is  uniformly bounded in $B(x,{\sf diam}(D))$, we deduce that $\left(e^{\gamma  t}v\left(B^x_{t\wedge\tau^x_{\partial B(x,{\sf diam}(D))}}\right)\right)_{t\geq 0}$ is a non-negative continuous supermartingale, so by Doob's stopping theorem in conjunction with Fatou's lemma we get
\begin{equation}\label{eq:poincare}
    \mathbb{E}\left\{e^{\gamma\tau^x_{\partial D}}\right\}
    \leq \mathbb{E}\left\{e^{\gamma\tau^x_{\partial B(x,{\sf diam}(D))}}\right\}
    \leq \frac{1}{1-\gamma\frac{{\sf diam}(D)^2}{2d} } \quad \mbox{for } 0<\gamma<\frac{2d}{{\sf diam}(D)^2}.
\end{equation}
Therefore, taking $\lambda=\frac{1}{q{\sf diam}(D)^2}$ in \eqref{eq:p-q} we get
\begin{equation*}
    \mathbb{E}\left\{e^{\frac{(p-1)\beta^2\varepsilon^2}{p\left[{\sf  diam}(D)^2\right]}N_\varepsilon^x }\right\}
    =
    \mathbb{E}\left\{e^{\frac{\beta^2\varepsilon^2}{q{\sf diam}(D)^2}N_\varepsilon^x }\right\}
     \leq 
     \mathbb{E}\left\{e^{\frac{8p^2}{q^2 {\sf diam}(D)^2}\tau^x_{\partial D}}\right\} ^{1/(2p)}
     2^{1/q}
     =
     \mathbb{E}\left\{e^{\frac{8(p-1)^2}{ {\sf diam}(D)^2}\tau^x_{\partial D}}\right\} ^{1/(2p)}
     2^{(p-1)/p}.
\end{equation*}
Note that the condition $p<1+\frac{\sqrt{d}}{2}$ is equivalent to $8(p-1)^2<2d$, hence by \eqref{eq:poincare} we get
\begin{equation*}
    \mathbb{E}\left\{e^{\frac{(p-1)\beta^2\varepsilon^2}{p \left[{\sf  diam}(D)^2\right]}N_\varepsilon^x }\right\}
    \leq 
    \left(\frac{1}{1-\frac{4(p-1)^2}{d}}\right)^{1/(2p)} 2^{(p-1)/p}.
\end{equation*}
The last part of the statement is obtained based on the first estimate and Markov inequality:
\begin{align*} 
    \mathbb{P}\left(N^{x}_\varepsilon\geq M\right)
    &=\mathbb{P}\left(
    e^{\frac{\beta^2\varepsilon^2}{3{\sf diam}(D)^2}N^{x}_\varepsilon}\geq e^{\frac{\beta^2\varepsilon^2}{3{\sf diam}(D)^2}M}
    \right)
    \leq \mathbb{E}\left\{e^{\frac{\beta^2\varepsilon^2}{3{\sf diam}(D)^2}N^{x}_\varepsilon}\right\}e^{-\frac{\beta^2\varepsilon^2}{3{\sf diam}(D)^2}M} \\
    &\leq 2 e^{-\frac{\beta^2\varepsilon^2}{3{\sf diam}(D)^2}M}, \quad M>0, x\in D, r(x)\geq \varepsilon\geq \varepsilon_0.
\end{align*}
\end{proof}

\begin{proof}[\textit{\textbf{Proof of \Cref{ex:defective-anulus}}}]
The $\delta$-defectiv convexity condition for this region amounts to the inequality:

\begin{equation*}
-\int_{\mathcal{A}_{R_1,R_2}} \nabla\varphi(x)\nabla r(x)\,dx \le \frac{\delta}{R_2-R_1}\int_{A(R_1,R_2)}\varphi(x)\,dx \quad \mbox{ for all } 0\leq\varphi\in C_c^\infty(A(R_1,R_2)),
\end{equation*} 
where we used the fact that the distance function $r$ is Lipschitz,  hence one can integrate by parts and discard the boundary terms due to  the fact that $\varphi\in C_c^\infty(A(R_1,R_2),\mathbb{R}_+)$. 

We use the fact that on $A\left(R_1,\frac{R_1+R_2}{2}\right)$ and $A\left(\frac{R_1+R_2}{2},R_2\right)$ the function $r$ is in fact smooth and we integrate by parts on each region to obtain that the left hand side above becomes:

\begin{align*}
\int_{A\left(R_1,\frac{R_1+R_2}{2}\right)}\varphi(x)\Delta r(x)\,dx+\int_{S^2}\int_{r=\frac{R_1+R_2}{2}}\varphi(\sigma,s)r'(s),s^{d-1}\,dsd\sigma+\non\\
+\int_{\mathcal{A}_{R_1,\frac{R_1+R_2}{2}}}\varphi(x)\Delta r(x)\,dx-\int_{S^2}\int_{r=\frac{R_1+R_2}{2}}\varphi(\sigma,s)r'(s)\,s^{d-1}\,dsd\sigma\non\\
=\int_{A\left(R_1,\frac{R_1+R_2}{2}\right)} \varphi(x)\Delta r(x)\,dx,
\end{align*} 
where we used spherical coordinates on the boundary and denoted by prime the derivative in the radial direction. 
Also,  the $\Delta r(x)$ is well defined, in a classical sense, on $A(R_1,R_2)$ except for the set of measure zero that in spherical coordinates is given by $\{(r,\sigma); r=\frac{R_1+R_2}{2},\sigma\in \mathbb{S}^2\}$.
Noting that:
$$
r(\sigma,s)=\left\{\begin{array}{ll}s-R_1 &\mbox{ for all } \sigma\in\mathbb{S}^2,s\in (R_1,\frac{R_1+R_2}{2}), \\[2mm]
R_2-s &\mbox{ for all } \sigma\in\mathbb{S}^2,s\in (\frac{R_1+R_2}{2},R_2)
\end{array}\right.
$$ we have $r'\equiv 1$ for all $x\in A\left(R_1,\frac{R_1+R_2}{2}\right)$ and $r'\equiv -1$ for all $x\in A\left(\frac{R_1+R_2}{2},R_2\right)$, hence  the $\delta$-defective convexity condition amounts to the inequality:

$$
(d-1)\bigg(\int_{\mathbb{S}^2}\int_{R_1}^{\frac{R_1+R_2}{2}} \varphi(\sigma,s)s^{d-2}ds-\int_{\mathbb{S}^2}\int_{\frac{R_1+R_2}{2}}^{R_2} \varphi(\sigma,s)s^{d-2}ds\bigg)\le \frac{\delta}{R_2-R_1}\int_{\mathbb{S}^2}\int_{R_1}^{R_2}\varphi(\sigma,s)s^{d-1}\,ds\,d\sigma, 
$$ which (taking into account that $\varphi\ge 0$) is satisfied if for instance $\frac{R_2}{R_{1}}<1+\frac{\delta}{d-1}$.
\end{proof}

\begin{proof}[\textit{\textbf{Proof of \Cref{ex:defective-tube}}}]
Arguing similarly as in the previous example, the $\delta$-defective convexity condition becomes:
\begin{equation}\label{cond:deltadefconv}
\int_{D_\eps}\varphi(x)\Delta r(x)\,dx \le \frac{\delta}{\eps}\int_{D_\eps}\varphi(x)\,dx\quad \mbox{ for all } 0\leq\varphi\in C_c^\infty(D_\varepsilon), 
\end{equation} 
where $\Delta r(x)$ is well-defined, in a classical sense, except on the set of measure zero $\Gamma_\eps := \big\{x+\frac{1}{2}\eps\,\,n(x) \in \mathbb{R}^d \ \big| \ 
x \in \Gamma  \big\}$. 
We also denote $ D_\eps^+ := \big\{x+\eps\,t\,n(x) \in \mathbb{R}^d \ \big| \
(x,t) \in \Gamma \times (\frac 12,1) \big\}$ respectively $ D_\eps^- := \big\{x+\eps\,t\,n(x) \in \mathbb{R}^d \ \big| \ 
 (x,t) \in \Gamma \times (0,\frac 12) \big\}$.
  
We recall (see for instance, \cite{GTbook}, Lemma $14.17$, p. 355) that 
  
\begin{equation*}
  \Delta r(x,t)=\left\{\begin{array}{ll} \sum_{i=1}^{d-1} \frac{k_i(x)}{1-k_i(x) r(x,t)} \,&\textrm{ if }t\in(0,\frac{1}{2})\\
-\sum_{i=1}^{d-1} \frac{k_i(x)}{1-k_i(x)(\eps- r(x,t))} \,&\textrm{ if }t\in(\frac{1}{2},1)\end{array}, \right.
\end{equation*} hence we have:
  
\begin{equation*}
  \int_{D_\eps}\varphi(x)\Delta r(x)\,dx\le \int_{D_\eps^-}\varphi(x)\left(\sum_{i=1}^{d-1} \frac{2k_i(x)}{2-k_i(x) \eps}\right)\,dx-
  \int_{D_\eps^+}\varphi(x)\left(\sum_{i=1}^{d-1} k_i(x)\right)\,dx
\end{equation*} thus the condition \eqref{eq:dconvex} 
holds for suitably small $\eps>0$.
\end{proof}

\begin{proof}[\textit{\textbf{Proof of \Cref{prop:logepsilon}}}]
We split the proof in two steps.

{{\bf Step I} (Regularization).}
Let $D_n:=\{x\in D : r(x)>1/n \}$ for each $n\geq n_0$, where $n_0$ is such that $D_{n_0}\neq \emptyset$.
Further, let $\rho\geq 0$ be a (smooth) mollifier on $\mathbb{R}^d$ such that $\supp{\rho}\subset B(0,1)$, and set $\rho_t(\cdot):=t^{-d}\rho(\cdot/t)$, $t>0$.
In particular,
\begin{equation} \label{eq:supp}
    \supp{\rho_t} +D_n\subset D \quad \mbox{ for all } t<1/n_0.
\end{equation}
Now let us consider that we extend $r$ from $\overline{D}$ to $\mathbb{R}^d$ by setting $r\equiv 0$ on $\mathbb{R}^d \setminus \overline{D}$, and set
$$
V_t=\sqrt{r_{t}}, \; \mbox{ where } \; r_{t}:=\rho_t \ast r \in C_c^2(\mathbb{R}^d).
$$
We claim that 
\begin{enumerate}
    \item[i)] $\Delta r_{t} \leq \frac{\delta}{2 \;\rm{rad}(D)}$ on $D_n$ for all $n\geq n_0$ and $0<t < 1/{n_0}$.
    \item[ii)] $\Delta V_t \leq -\frac{1}{4}{r_{t}}^{-3/2}[|\nabla r_{t}|^2 -\delta]$ on $D_n$ for all $n\geq n_0$ and $0<t < 1/{n_0}$.
\end{enumerate}
To prove the claim, note first that by a simple calculation we get
\begin{equation*}
\Delta V_t=-\frac{1}{4}{r_{t}}^{-3/2}[|\nabla r_{t}|^2 - 2r_t\Delta r_{t}],
\end{equation*}
hence ii) follows from i) and the fact that $r_t\leq \rm{rad}(D)$.
So, it remains to prove the first assertion of the claim: Let $0\leq \varphi \in C_c^\infty(D_n), n\geq n_0, t<1/n_0$, and proceed with integration by parts and Fubini's theorem as follows
\begin{align*}
    \int_{D_n}\varphi \Delta r_{t} \;dx&=-\int_{D_n}\nabla r_{t} \cdot \nabla \varphi \; dx = -\int_{\mathbb{R}^d}\rho_t(y)\int_{D_n}\nabla r(x-y) \cdot \nabla \varphi(x) \; dx\;dy\\
    &=-\int_{\mathbb{R}^d}\rho_t(y)\int_{\mathbb{R}^d}\nabla r(x-y) \cdot \nabla \varphi(x) \; dx\;dy\\
    &=-\int_{\supp{\rho_t}}\rho_t(y)\int_{\mathbb{R}^d}\nabla r(x) \cdot \nabla \varphi(\cdot+y)(x) \; dx\;dy,\\
    \intertext{so, by \eqref{eq:supp} 
    and then using \eqref{eq:dconvex} we can continue with}
    &=-\int_{\supp{\rho_t}}\rho_t(y)\int_{D}\nabla r(x) \cdot \nabla \varphi(\cdot+y)(x) \; dx\;dy \leq \frac{\delta}{2 \;{\rm rad}(D)}\int_{\supp{\rho_t}}\rho_t(y)\int_{D} \varphi(x+y) \; dx\;dy\\
    &= \frac{\delta}{2 \;{\rm rad}(D)}\int_{\supp{\rho_t}}\int_{\mathbb{R}^d} \varphi(x+y) \; dx\;dy= \frac{\delta}{2 \;{\rm rad}(D)}\int_{\mathbb{R}^d} \varphi(x) \; dx\\
    &=\frac{\delta}{2 \;{\rm rad}(D)}\rho_t(y)\int_{D_n} \varphi(x) \; dx,
\end{align*}
which proves i) and hence the entire claim.

Now, by Ito's formula in corroboration with the claim proved above yield for $x\in D_n$, $r(x)\geq \varepsilon_0$:
\begin{align}
     \mathbb{E}\left\{ V_t\left(B^x_{\tau_{(D_n\cap B(x,\widetilde{r}(x)))^c}}\right)\right\}&
    =V_t(x)+\mathbb{E}
    \left\{\int_0^{\tau_{(D_n\cap B(x,\widetilde{r}(x)))^c}}\Delta V_t(B_s^x) \; ds\right\} \nonumber\\
    & \leq V_t(x)-\frac{1}{4}\mathbb{E}\left\{ \int_0^{\tau_{(D_n\cap B(x,\widetilde{r}(x)))^c}}(r^t)^{-3/2}(B_s^x)\left[|\nabla r^t|^2(B_s^x) -\delta\right]\; ds\right\} \label{eq:V_epsilon}.
\end{align}
{{\bf Step II} (Passing to the limit in \eqref{eq:V_epsilon}).}
The next step is to let $t \to 0$ and then $n\to \infty$ in \eqref{eq:V_epsilon}.
To this end, note that because $r\in C(\overline{D})\cap W^{1,\infty}(D)$, we have
\begin{enumerate}
    \item[iii)] $\lim\limits_{t \to 0}r_{t}=r$ uniformly on $D$,
    \item [iv)] $\lim\limits_{t \to 0}\nabla r_{t}=\nabla r$ a.e. and boundedly on $D$.
\end{enumerate}
In particular, because $\inf\limits_{x\in D_n}r(x)\geq 1/n$ for large enough $n_{0}$ and $n\ge n_{0}$, we have that 
\begin{enumerate}
    \item[v)] $\lim\limits_{t \to 0}r_{t}^{-3/2}=r^{-3/2}$ boundedly on $D_n$, for each $n\geq n_0$.
\end{enumerate}
We are now in the position to let $t\to 0$ in \eqref{eq:V_epsilon} to get that for $x\in D_n$, $r(x)\geq \varepsilon_0$,
\begin{align*}
    \mathbb{E}\left\{ V\left(B^x_{\tau_{(D_n\cap B(x,\widetilde{r}(x)))^c}}\right)\right\} 
    &\leq V(x)-\frac{1 -\delta}{4}\mathbb{E}
    \left\{ \int_0^{\tau_{(D_n\cap B(x,\widetilde{r}(x)))^c}}r^{-3/2}(B_t^x)\; dt\right\},
\end{align*}
where we have used that $|\nabla r|=1$ a.e.  The fact that $|\nabla r|=1$, follows from two basic facts.  On one hand, $r$ is $1$-Lipschitz so $|\nabla r|\le 1$.  On the other hand, from Lipschitz conditions and Rademacher theorem, $r$ is differentiable almost everywhere.  If $x$ is a point where $r$ is differentiable, and $y\in \partial D$ such that $d(x,y)=r(x)$, and $v=y-x$, then it is easy to see that the derivative of $r$ in the direction $v$ is constant $1$, thus the claim.   

Now, on the one hand, since $r\leq 2r(x)$ on $B(x,\widetilde{r}(x))\subset B(x,r(x))$ we deduce that
\begin{align}\label{eq:D_n}
    \mathbb{E}\left\{ V\left(B^x_{\tau_{(D_n\cap B(x,\widetilde{r}(x)))^c}}\right)
    \right\}
    &\leq V(x)-\frac{1-\delta}{4r(x)^{3/2}}\mathbb{E}\left\{ \tau_{(D_n\cap B(x,\widetilde{r}(x)))^c}
    \right\}, \quad x\in D_n.
\end{align}
On the other hand, $D_n\mathop{\nearrow}\limits_n D$, so $\tau_{(D_n\cap B(x,\widetilde{r}(x)))^c} \mathop{\longrightarrow}\limits_n \tau_{B(x,\widetilde{r}(x))^c}$ a.s., hence for all $x\in D$, $r(x)\geq \varepsilon_0$,
\begin{align*}
    P V(x)&=\mathbb{E}\left\{ V\left(B^x_{\tau_{B(x,\widetilde{r}(x))}}\right)\right\} 
    =\lim\limits_n\mathbb{E}\left\{ V\left(B^x_{\tau_{(D_n\cap B(x,\widetilde{r}(x)))^c}}\right)\right\},\\
    \intertext{hence letting $n\to \infty$ in \eqref{eq:D_n} we can continue with}
    P V(x)
    &\leq V(x)-\frac{1-\delta}{4r(x)^{3/2}}\mathbb{E}\left\{ \tau_{B(x,\widetilde{r}(x))^c}\right\} =V(x)-\frac{1-\delta}{4r(x)^{3/2}}\frac{\widetilde{r}(x)^2}{d}\\
    &\leq\left(1-\frac{\beta^2(1-\delta)}{4d}\right)V(x), \quad x\in D, r(x)\geq \varepsilon_0,
\end{align*}
which proves \eqref{eq:lyapunov}.  Notice that for the last inequality we used that $\widetilde{r}$ is a $(\beta,\varepsilon_0)$-distance, see \Cref{defi:adistance}, hence $\beta r(x)\leq \widetilde{r}(x), x\in D, r(x)\geq \varepsilon_0$.

Now, we modify the kernel $P$ and $V$ as follows:
\begin{equation*}
    P_{\varepsilon_0}:=1_{D\setminus D_{\varepsilon_0}}P+1_{D_{\varepsilon_0}}\delta, \quad V_{\varepsilon_0}=1_{D\setminus D_{\varepsilon_0}} V,\quad \mbox{where } D_{\varepsilon_0}:=\{x\in D : r(x)\leq \varepsilon_0\}.
\end{equation*}
It is easy to see that
\begin{equation*}
    P_{\varepsilon_0}V_{\varepsilon_0}\leq \left(1-\frac{\beta^2(1-\delta)}{4d}\right) V_{\varepsilon_0} \quad \mbox{ on } D.
\end{equation*}
Note that on the event $X^{x}_0,...,X^{x}_M\in D\setminus D_{\varepsilon_0}$, the chain $(X^x_i)_{1\leq i\leq M}$ has transition kernel $P_{\varepsilon_0}$.

The tail estimate \eqref{eq:exptail} can be deduced from \eqref{eq:lyapunov} and Markov inequality, as follows:
\begin{align*}
   \forall \varepsilon\geq \varepsilon_0:\quad
    \mathbb{P}(N_\varepsilon^{x}> M)
    &
    \leq \mathbb{P}(r(X^{x}_M)\geq \varepsilon)=\mathbb{P}(V(X^{x}_M)\geq \sqrt{\varepsilon}; X^{x}_0,...,X^{x}_M\in D\setminus D_{\varepsilon_0})\\
    &\leq \frac{1}{\sqrt{\varepsilon}}\mathbb{E}
    \left\{ V(X^{x}_M); X^{x}_0,...,X^{x}_M\in D\setminus D_{\varepsilon_0}\right\} 
    =\frac{1}{\sqrt{\varepsilon}}P_{\varepsilon_0}^MV_{\varepsilon_0}(x)\\
    &\leq \left(1-\frac{\beta^2(1-\delta)}{4d}\right)^M \frac{V(x)}{\sqrt{\varepsilon}}, \quad x\in D, r(x)\geq \varepsilon_0, \quad \mbox{hence for all } x\in D.
\end{align*}

Let us finally conclude the proof by showing \eqref{eq:logepsilon}:
\begin{align*}
    \mathbb{E}\left\{ a^{N_\varepsilon^{x}}\right\} &
    =\sum\limits_{k\geq 0} a^k \mathbb{P}(N_\varepsilon^{x}=k)\leq 1+\sum\limits_{k\geq 1} a^k \mathbb{P}(N_\varepsilon^{x}>k-1)\\
    &\leq 1+ \sum\limits_{k\geq 1} a^k \left(1-\frac{\beta^2(1-\delta)}{4d}\right)^{k-1} \frac{V(x)}{\sqrt{\varepsilon}}
    =1+a\frac{V(x)}{\sqrt{\varepsilon}} \sum\limits_{k\geq 0}  a^k\left(1-\frac{\beta^2(1-\delta)}{4d}\right)^{k}\\
    &=1+\frac{a}{1-a\delta_d}\frac{V(x)}{\sqrt{\varepsilon}}, \quad x\in D.
\end{align*}
\end{proof}

\subsection{Proofs for \Cref{ss:errors}}
\label{pf:3}

\begin{proof}[\textit{\textbf{Proof of \Cref{prop:errorspoisson}}}]
First of all, note that by similar arguments to those used in the proof of \Cref{coro:representation},
\begin{equation*}
    u_M(x)=\mathbb{E}\left\{ \sum\limits_{n=0}^{M-1}\int_{\tau_n^{x}}^{\tau_{n+1}^{x}} f(B^x_t)\; dt\right\}=\mathbb{E}\left\{ \int_0^{\tau_{M}^{x}} f(B^x_t) \; dt \right\}, \quad x\in D.
\end{equation*}
Therefore,
\begin{align*}
    |u(x)-u_M(x)|&=\left| \mathbb{E}\left\{ \int_{\tau_{M}^{x}}^{\tau_{D^c}} f(B^x_t) \; dt \right\} \right|\leq |f|_\infty \mathbb{E}\left\{ \tau_{D^c}-\tau_{M}^x \right\}= |f|_\infty \mathbb{E}^x\left\{ \mathbb{E}^{B^x_{\tau_M^{x}}}\left\{ \tau_{D^c}\right\} \right\}\\
    &=|f|_\infty \mathbb{E}\left\{ v(B^x_{\tau_M^{x}})\right\}=|f|_\infty \mathbb{E}\left\{ v(B^x_{\tau_M^{x}}) ; N_\varepsilon^{x} \leq M\right\}+|f|_\infty \mathbb{E}\left\{ v(B^x_{\tau_M^{x}}) ; N_\varepsilon^{x} > M\right\}\\
    &\leq |f|_\infty \left [\mathbb{E}\left\{ v(B^x_{\tau_{M\vee N_\varepsilon^{x}}^{x}}) ; N_\varepsilon^{x} \leq M\right\}+ |v|_\infty \mathbb{P}(N_\varepsilon^{x} > M)\right]\\
    &\leq |f|_\infty \left [\mathbb{E}\left\{ v(B^x_{\tau_{M\vee N_\varepsilon^{x}}^{x}}) ; N_\varepsilon^{x} \leq M\right\}+ \frac{1}{d}{\sf diam}(D)^2 2e^{-\frac{\beta^2\varepsilon^2}{3 {\sf diam}(D)^2}M}\right], \quad x\in D,
\end{align*}
where the last inequality follows by \Cref{coro:v_inf} and \Cref{prop:nsteps_general}.
Now, one can see that $(v(B_{\tau_n^x}))_{n\geq 0}$ is a bounded and non-negative supermartingale, hence by Doob's stopping theorem we get that
\begin{align*}
    \mathbb{E}\left\{ v(B^x_{\tau_{M\vee N_\varepsilon^{x}}^{x}}) ; N_\varepsilon^{x} \leq M\right\}=\mathbb{E}\left\{ 1_{[N_\varepsilon^{x} \leq M]} \mathbb{E}\left[ v(B^x_{\tau_{M\vee N_\varepsilon^{x}}^{x}}) \mid \mathcal{F}_{N_\varepsilon^{x}}\right]\right\}\leq \mathbb{E}\left\{v(B^x_{\tau_{N_\varepsilon^{x}}^{x}})\right\}=v(x,\varepsilon)\le  v_\infty(\varepsilon),
\end{align*}
hence the desired estimate is now fully justified.
\end{proof}

Further, we need the following lemma:

\begin{lem}\label{prop6}
Let $\varepsilon_0>0$, $\beta\in (0,1]$, $\widetilde{r}$ be a $(\beta,\varepsilon_0)$-distance, $x\in D$, and
$(X_n^x)_{n\geq 0}$ be the corresponding $\widetilde{r}$-WoS chain.
If $\varepsilon\geq \varepsilon_0$ and $\tau'\geq \tau$ are finite stopping times such that $\tau \geq N_{\varepsilon}^{x}$,
whilst $g$ is $\alpha$-H\"{o}lder on $\overline{D}$ for some $\alpha\in [0,1]$, then
$$
\mathbb{E}\left\{ \left|g(X^x_{\tau'})-g(X^x_{\tau})\right|\right\} 
\leq d^{\alpha/2}|g|_\alpha\cdot v(x,\varepsilon)^{\alpha/2} \leq  d^{\alpha/2}|g|_\alpha\cdot |v|^{\alpha/2}_{\infty}(\varepsilon).
$$
\end{lem}

\begin{proof}
It is easy to see that for $z\in\mathbb{R}^d$, $(\langle X^{x}_n,z \rangle - \langle x,z \rangle)_{n\geq 1}$ is a bounded martingale, hence 
\begin{equation*}
    \mathbb{E}\{\langle X^{x}_{\tau'}, X^{x}_{\tau}\rangle\} = \mathbb{E}\left\{ |X_{\tau}^{x}|^2\right\},\quad 
    \mbox { and thus }
    \quad \mathbb{E}\left\{ | X^{x}_{\tau'}- X^{x}_{\tau}|^2 \right\} = \mathbb{E}\left\{ |X^{x}_{\tau'}|^2\right\} - \mathbb{E}\left\{ |X^{x}_{\tau}|^2\right\}. 
\end{equation*}
By employing the martingale problem for the Markov chain $(X_n^{x})_{n\geq 0}$, we get that for any finite stopping time $T$,
 $ \mathbb{E}\left\{ |X^{x}_T|^2 \right\} 
 =|x|^2
+\mathbb{E}
\left\{ \sum_{i=0}^{T-1}\widetilde{r}^2(X^{x}_i)\right\}, $ hence
$
 \mathbb{E}
 \left\{ | X^{x}_{\tau'}- X^{x}_{\tau}|^2\right\} =  \mathbb{E}\left\{ \sum\limits_{i=\tau}^{ \tau'-1}\widetilde{r}^2(X^{x}_i)\right\}$
 and therefore
 \begin{align*}
\left|\mathbb{E}\left\{ g(X^{x}_{\tau'})\right\} -\mathbb{E}\left\{ g(X^{x}_{\tau})\right\} \right|&
\leq |g|_\alpha
\mathbb{E} \left\{ \sum_{i=\tau}^{\infty}\widetilde{r}^2(X^{x}_i)\right\}  
^{\alpha/2}
\leq  
d^{\alpha/2}|g|_\alpha
\mathbb{E}\left\{ \sum_{i=N^{x}_\varepsilon}^{\infty}\tau^{x}_{i+1}-\tau^{x}_i\right\}  
^{\alpha/2}\\
&\leq d^{\alpha/2}|g|_\alpha\cdot \mathbb{E}\left\{ \tau^x_{\partial D}-\tau^{x}_{N^{x}_\varepsilon}\right\} ^{\alpha/2}=d^{\alpha/2}|g|_\alpha\cdot \mathbb{E}\left\{  v\left(B^x_{\tau^{x}_{N^{x}_\varepsilon}}\right) \right\}  ^{\alpha/2}\\
&\leq d^{\alpha/2}|g|_\alpha\cdot |v|_{\infty}^{\alpha/2}(\varepsilon).
\end{align*}
\end{proof}

\begin{proof}[\textit{\textbf{Proof of \Cref{prop:errorsdirichlet}}}]
Recall that $\left(B^x_{\tau^{x}_n}\right)_n$ and $(X^{x}_n)_n$ are equal in law and hence $B^x_{\tau^{x}_{N^{x}_\varepsilon}}$ and $X^{x}_{N^{x}_\varepsilon}$ are also equal in law.
In particular, using that $g$ is continuous on $\overline{D}$ and that $\tau^{x}_{n}$ converges a.s. to some $\tau_\infty^x$ such that $\tau^x_{D_\varepsilon}\leq \tau_\infty^x\leq \tau^x_{D^c}$ as already mentioned in \Cref{rem:beta0distance},
\begin{equation*}
    \mathbb{E}\{ g(B^x_{\tau^x_\infty})\} =\lim_n \mathbb{E}\{ g(B^x_{\tau^{x}_{n}})\} =\lim_n  \mathbb{E}\{ g(X^{x}_n)\}  \quad \mbox{ for all } x\in D.
\end{equation*}
Also, if $T$ is a finite random time then
$$\lim_n \mathbb{E}\{ g(X^{x}_{n \vee T})\} =\lim_n \{\mathbb{E}\{ g(X^{x}_{n});\;T< n]+\mathbb{E}\{ g(X^{x}_T);\;T\geq n\}\}=\mathbb{E}\{g(B^x_{\tau^x_\infty})\}, \ x\in D.
$$
Now let us fix $\varepsilon\geq \varepsilon_0$ and argue as follows:
\begin{align*}
    &|u(x)-u_M(x)|\\
    &\leq \left|\mathbb{E}\left\{g(B^x_{\tau^x_{D^c}})-g(B^x_{\tau^x_\infty})\right\}\right|
    +\lim\limits_n\left|\mathbb{E}\left\{ g(X^{x}_n)-g(X^{x}_M)\right\} \right|\\
    &\leq |g|_\alpha\mathbb{E}\left\{\|B^x_{\tau^x_{D^c}}-B^x_{\tau^x_\infty}\|^\alpha\right\}
    +\lim\limits_n\left|\mathbb{E}
    \left\{ g(X^{x}_n)-g(X^{x}_M)  ; \; N_\varepsilon^{x}\leq M\right\} \right| 
    +\lim\limits_n\left|\mathbb{E}
    \left\{ g(X^{x}_n)-g(X^{x}_M)  ; \; N_\varepsilon^{x}> M\right\} \right|\\
    &\leq |g|_\alpha\mathbb{E}\left\{\|B^x_{\tau^x_{D^c}}-B^x_{\tau^x_\infty}\|^2\right\}^{\alpha/2} 
    +\lim\limits_n\left|\mathbb{E}
    \left\{ g(X^{x}_{n\vee M\vee N^{x}_\varepsilon})-g\left(X^{x}_{M \vee N_\varepsilon^{x}}\right); \; N_\varepsilon^{x}\leq M\right\} \right| + 2|g|_\infty \mathbb{P}(N_\varepsilon^{x}> M)\\
    &\leq |g|_\alpha d^{\alpha/2}|v|^{\alpha/2}_{\infty}(\varepsilon) 
    +\lim\limits_n\left|\mathbb{E}
    \left\{ g(X^{x}_{n\vee M\vee N^{x}_\varepsilon})-g\left(X^{x}_{M \vee N_\varepsilon^{x}}\right); \; N_\varepsilon^{x}\leq M\right\} \right| + 2|g|_\infty \mathbb{P}(N_\varepsilon^{x}> M),
\end{align*}
where we have used that ${\sf dist}(B^x_{\tau^x_\infty},\partial D)\leq \varepsilon$.
Further, using \Cref{prop6} together with \Cref{prop:nsteps_general}, we get
\begin{align*}
    |u(x)-u_M(x)|
    &\leq 2d^{\alpha/2}|g|_\alpha\cdot v^{\alpha/2}(x,\varepsilon) + 2|g|_\infty \mathbb{P}(N_\varepsilon^{x}> M)
    \leq 2d^{\alpha/2}|g|_\alpha\cdot v^{\alpha/2}(x,\varepsilon) + 4|g|_\infty e^{-\frac{\beta^2\varepsilon^2}{3{\sf diam}(D)^2}M},\\
    &\leq 2d^{\alpha/2}|g|_\alpha\cdot |v|^{\alpha/2}_{\infty}(\varepsilon) + 4|g|_\infty e^{-\frac{\beta^2\varepsilon^2}{3{\sf diam}(D)^2}M}.
\end{align*}
\end{proof}

\subsection{Proofs for \Cref{ss:meantailestimates}}
\label{pf:4}

\begin{proof}[\textit{\textbf{Proof of \Cref{prop:meanestimates}}}]
We have 
\begin{align*}
    \mathbb{E}&\left\{ \left | u (\cdot)-u_M^{N}(\cdot) \right |^2_{L^2(D)} \right\}  \\
    &\leq 2\left | u (\cdot)-u_M(\cdot) \right |^2_{L^2(D)}+ 2\mathbb{E}
    \left\{ \left | u_M (\cdot)-u_M^{N}(\cdot) \right |^2_{L^2(D)} \right\}\\
    & \leq 2\lambda (D) \sup\limits_{x\in D} |u(x)-u_M(x)|+2\int\limits_D\mathbb{E}\left\{\left | u_M (\cdot)-u_M^{N}(\cdot) \right |^2 \right\}\;dx\\
    &=2\lambda (D) \sup\limits_{x\in D} |u(x)-u_M(x)|^2+\frac{2}{N}\int\limits_D\mathbb{E}\left\{\left | u_M (\cdot)-u_M^{1}(\cdot) \right |^2 \right\}\;dx\\
    &\leq 2\lambda (D) \sup\limits_{x\in D} |u(x)-u_M(x)|^2+\frac{4}{N}\left(\int\limits _D\left[|g|^2_\infty+\frac{1}{d^2}M|f|^2_\infty{\sf diam}(D)^2 \mathbb{E}\left\{\sum\limits_{k=1}^M \widetilde{r}^2(X_{k-1}^{x,1})\right\}\;\right] dx\right)\\
    &\leq 2\lambda (D) \sup\limits_{x\in D} |u(x)-u_M(x)|^2+\frac{4}{N}\left(\int\limits _D\left[|g|^2_\infty+\frac{1}{d^2}M|f|^2_\infty{\sf diam}(D)^2 \mathbb{E}\left\{\tau^x_{D^c}\right\}\;\right] dx\right)\\
    \intertext{and by \Cref{coro:v_inf}}
    &\leq 2\lambda (D)\left[ \sup\limits_{x\in D} |u(x)-u_M(x)|^2+\frac{2(|g|^2_\infty+\frac{1}{d^3}M|f|^2_\infty {\sf diam}(D)^4}{N}\right].
\end{align*}
\end{proof}

Before proving the main result,  \Cref{thm:main}, we need several preliminary Lemmas.
\begin{lem} \label{lem:lipM}
Let $\varepsilon>0$ and $\beta\in(0,1]$.
If $f$ and $g$ are $\alpha$-H\"older for some $\alpha\in [0,1]$, and $M,N\geq 1$, then
\begin{equation*}
    |u_M^{N}(x)-u_M^{N}(y)|\leq \left(|g|_\alpha+\frac{{\rm diam}(D)^2|f|_\alpha+2{\rm diam}(D)|f|_\infty}{d}\right)(2+|\widetilde{r}|_1)^M\left(|x-y|^\alpha \vee |x-y|\right),
\end{equation*}
for all $x,y\in D$ almost surely.
In particular,
\begin{equation*}
    |u_M(x)-u_M(y)|\leq \left(|g|_\alpha+\frac{{\rm diam}(D)^2|f|_\alpha+2{\rm diam}(D)|f|_\infty}{d}\right)(2+|\widetilde{r}|_1)^M\left(|x-y|^\alpha \vee |x-y|\right) \quad \mbox{for all } x,y\in D.
\end{equation*}
\end{lem}

\begin{proof}
Clearly, it is sufficient to prove the estimate for $|u_M^{i}(x)-u_M^{i}(y)|\leq 1$, independently of $i$, $1\leq i\leq N$.
To this end, since $\widetilde{r}$ is Lipschitz
\begin{align*}
    |X^{x,i}_{M}-X^{y,i}_{M}|&\leq|X^{x,i}_{M-1}-X^{y,i}_{M-1}|+|\widetilde{r}(X^{x,i}_{M-1})-\widetilde{r}(X^{x,i}_{M-1})|\\
    &\leq (1+|\widetilde{r}|_1)|X^{x,i}_{M-1}-X^{y,i}_{M-1}|\\
    &\leq (1+|\widetilde{r}|_1)^M|x-y|  \quad \mbox{ for all } x,y\in D, M\geq 0.
\end{align*}
Therefore,
\begin{align*}
    |u_M^{i}(x)-u_M^{i}(y)|&\leq |g|_\alpha|X^{x,i}_{M}-X^{y,i}_{M}|^\alpha+\frac{1}{d}\sum\limits_{k=1}^M 2|f|_\infty{\rm diam}(D)|\widetilde{r}(X^{x,i}_{k-1})-\widetilde{r}(X^{y,i}_{k-1})|\\
    &\;\phantom{\leq|g|_\alpha|X^{x,i}_{M}-X^{y,i}_{M}|^\alpha}+\frac{1}{d}\sum\limits_{k=1}^M{\rm diam}(D)^2|f|_\alpha[|X^{x,i}_{k-1}-X^{y,i}_{k-1}|+|\widetilde{r}(X^{x,i}_{k-1})-\widetilde{r}(X^{x,i}_{k-1})|]^\alpha\\
    &\leq |g|_\alpha(1+|\widetilde{r}|_1)^{M\alpha}|x-y|^\alpha+\frac{2{\rm diam}(D)|f|_\infty|x-y|}{d}\sum\limits_{k=1}^M (1+|\widetilde{r}|_1)^{k-1}\\
    &\;\phantom{\leq|g|_\alpha|X^{x,i}_{M}-X^{y,i}_{M}|^\alpha}+\frac{{\rm diam}(D)^2|f|_\alpha|x-y|^\alpha}{d}\sum\limits_{k=1}^M(1+|\widetilde{r}|_1)^{\alpha (k-1)}\\
    &\leq \left(|g|_\alpha+\frac{{\rm diam}(D)^2|f|_\alpha+2{\rm diam}(D)|f|_\infty}{d}\right)(2+|\widetilde{r}|_1)^M  \left(|x-y|^\alpha \vee |x-y|\right).
\end{align*}
\end{proof}

The next lemma is the well-known Hoeffding's inequality:
\begin{lem} \label{lem:Hoef}
Suppose that $(Z_i)_{i\geq 1}$ are iid real random variables such that $a_i\leq Z_i\leq b_i$ for all $i$. Then for all $N\in \mathbb{R}$ and $\gamma \geq 0$
$$
\mathbb{P}\left(\left|\mathbb{E}\{ Z_1\} -\frac{1}{N}\sum_{i=1}^N Z_i\right|\geq \gamma\right)\leq 2 e^{-\frac{2N^2\gamma^2}{\sum_{i=1}^N (b_i-a_i)^2}}.
$$
\end{lem}

Using Hoeffding's inequality we immediately get the following estimate.
\begin{coro}\label{coro5}
For all $N, M\in \mathbb{N}$, and $\gamma \geq 0$ we have
$$
\mathbb{P}\left(\left | u_{M} (x)-u_M^{N}(x) \right | \geq \gamma \right) \leq 2 e^{-\frac{N\gamma^2}{(|g|_\infty+M{\sf diam}(D)^2|f|_\infty/d)^2}}, \quad x\in D.
$$
\end{coro}

\begin{proof}
The result follows directly from \Cref{lem:Hoef} since 
$$
\left|g(X^{x,i}_{M})+\frac{1}{d}\sum\limits_{k= 1}^M \widetilde{r}^2(X^{x,i}_{k-1})f\left(X^{x,i}_{k-1}+\widetilde{r}(X^{x,i}_{k-1})Y^i\right)\right|_\infty\leq |g|_\infty+M{\sf diam}(D)^2|f|_\infty/d.
$$
\end{proof}

Finally, we are in the position to prove the main theorem.
\begin{proof}[\textit{\textbf{Proof of \Cref{thm:main}}}]
First of all, assume without loss of generality that $D\subset [0,{\sf diam}(D)]^d$, and for each $M,K\geq 1$ consider the grid
$$
\Xi=\Xi(K,M,\alpha,|\widetilde{r}|_1):=\left\{\frac{i {\sf diam}(D)}{K(1+|\widetilde{r}|_1)^{\lceil M/\alpha\rceil}} :\; 1\leq i\leq K(2+|\widetilde{r}|_1)^{\lceil M/\alpha\rceil} \right\}^d\cap D.
$$
For $x\in D$ such that $x\in \left[\frac{i_1{\sf diam}(D)}{K(2+|\widetilde{r}|_1)^{\lceil M/\alpha\rceil}} ,\frac{(i_1+1){\sf diam}(D)}{K(2+|\widetilde{r}|_1)^{\lceil M/\alpha\rceil}} \right)\times \cdots \times \left[\frac{i_d{\sf diam}(D)}{K(2+|\widetilde{r}|_1)^{\lceil M/\alpha\rceil}}, \frac{(i_d+1){\sf diam}(D)}{K(2+|\widetilde{r}|_1)^{\lceil M/\alpha\rceil}} \right)$ we set
\begin{equation*}
    x^{\Xi}:=\left(\frac{i_1{\sf diam}(D)}{K(2+|\widetilde{r}|_1)^{\lceil M/\alpha\rceil}} ,\cdots,\frac{i_d{\sf diam}(D)}{K(2+|\widetilde{r}|_1)^{\lceil M/\alpha\rceil}} \right).
\end{equation*}
Note that
\begin{fleqn}
\begin{align*}
    \sup_{x\in D} \left | u(x)-u_M^{N}(x)\right |
    &\leq \sup_{x\in D} \left | u(x)-u_M(x)\right |+\sup_{x\in D} \left | u_M(x)-u_M^{N}(x)\right |\\
    &\leq \sup_{x\in D} \left | u(x)-u_M(x)\right | + \sup_{x\in F} \left | u_M(x)-u_M^{N}(x)\right |\\
    &\;\phantom{ \leq \sup_{x\in D} \left | u(x)-u_M(x)\right | }+ 2\left(|g|_\alpha+\frac{{\rm diam}(D)^2|f|_\alpha+2{\rm diam}(D)|f|_\infty}{d}\right)\left(\frac{{\sf diam}(D)}{K}\right)^\alpha,
\end{align*}
\end{fleqn}
where the last inequality follows from \Cref{lem:lipM} and by the fact that 
$$
|x-x_{\Xi}|\leq \frac{{\sf diam}(D)}{K(2+|\widetilde{r}|_1)^{M/\alpha}} \mbox{ for all } \; x\in D.
$$
Consequently, by setting
\begin{equation*}
    \widetilde{\gamma}:=\gamma-\sup_{x\in D} \left | u(x)-u_M(x)\right | - 2\left(|g|_\alpha+\frac{{\rm diam}(D)^2|f|_\alpha+2{\rm diam}(D)|f|_\infty}{d}\right)\left(\frac{{\sf diam}(D)}{K}\right)^\alpha
\end{equation*}
and using union bound inequality we have
\begin{align*}
   \mathbb{P}&\left( \sup_{x\in D} \left | u(x)-u_M^{N}(x)\right| \geq \gamma \right)\leq\mathbb{P}\left( \sup_{x\in \Xi} \left | u_{M}(x)-u_M^{N}(x)\right | \geq \widetilde{\gamma}\right)\leq \sum\limits_{x\in \Xi} \mathbb{P}\left( \left | u_{M}(x)-u_M^{N}(x)\right | \geq \widetilde{\gamma}\right),
\end{align*}
so, the two desired estimates now follow by
\Cref{thm:stepdet} and \Cref{coro5}. 

The assertions about the particular domains are clear. 

For the statement about the expectation stated in \eqref{eq:main:e}, we only have to involve the following Lemma.  

\begin{lem}
If $X$ is a non-negative random variable with the property that there exist constants and $c_1,c_2,A\ge0$ such that
\begin{equation}
\mathbb{P}(X\ge t)\le  2e^{c_1-c_2((t-A)^+)^2}\text{ for all }t\ge0, 
\end{equation}
then 
\[
\mathbb{E}[X]\le A+\frac{\sqrt{c_1+\log(2)}+1}{\sqrt{c_2}}.  
\]
\end{lem}

\begin{proof}
Start with $\lambda\ge A$ and write 
\[
\mathbb{E}[X]=\int_0^\infty \mathbb{P}(X\ge t)dt=\int_0^\lambda \mathbb{P}(X\ge t)dt+\int_\lambda^\infty \mathbb{P}(X\ge t)dt\le \lambda+\int_\lambda^\infty  2 e^{c_1-c_2(t-A)^2}=\lambda+\int_{\lambda-A}^\infty 2e^{c_1-c_2 t^2}dt.
\]
Optimizing over $\lambda\ge A$, yields the optimum point as  
\[
\lambda^*=A+\sqrt{\frac{c_1+\log(2)}{c_2}}
\]
which in turn yields 
\[
\int_{\lambda^*-A}^\infty 2 e^{c_1-c_2 t^2}dt\le 2e^{c_1}\int_{\lambda^*-A}^\infty \frac{t}{\lambda^*-A}e^{-c_2 t^2}dt=\frac{1}{2\sqrt{(c_1+\log(2))c_2}}\le \frac{1}{\sqrt{c_2}}.
\]
This concludes the estimate.  
\end{proof}

The rest of the statements in the theorem are straightforward now.  

\end{proof}

\begin{proof}[\textit{\textbf{Proof of \Cref{coro:noKepsilon}}}]
Note that since $\widetilde{r}$ is a $(\beta,\varepsilon_0)$-distance, it is also a $(\beta,\varepsilon)$-distance since $\varepsilon_0\leq \varepsilon$.

Now, using \eqref{eq:K} we get that
\begin{equation*}
    2\left(2|g|_\alpha+\frac{{\rm diam}(D)^2|f|_\alpha+2{\rm diam}(D)|f|_\infty}{d}\right)\left(\frac{{\sf diam}(D)}{K}\right)^\alpha\leq \frac{\gamma}{4}.
\end{equation*}
Also, because ${\sf adiam}(D)<\infty$, by \Cref{prop:extime} we have that $|v|_\infty(\varepsilon)\leq \varepsilon {\sf adiam}(D)$. 
Therefore, since $\varepsilon\leq 1$ is given by \eqref{eq:epsilon} { and using that $\varepsilon\leq \varepsilon^{\alpha/2}$}, we get that 
\begin{equation*}
    2d^{\alpha/2}|g|_\alpha\cdot |v|^{\alpha/2}_{\infty}(\varepsilon)+ |f|_\infty |v|_{\infty}(\varepsilon)\leq \frac{\gamma}{4}.
\end{equation*}
Now, if $M$ is as in \eqref{eq:Mii}, then
\begin{equation*}
    (4|g|_\infty+\frac{2}{d}{\sf diam}(D)^2|f|_\infty) e^{-\frac{\beta^2\varepsilon^2}{3{\sf diam}(D)^2}M}\leq \frac{\gamma}{4}.
\end{equation*}
Furthermore, if $D$ is $\delta$-defective convex and $M$ is chosen as in \eqref{eq:Mi}, we have that
\begin{equation*}
    (4|g|_\infty+\frac{2}{d}{\sf diam}(D)^2|f|_\infty)a_d^M\sqrt{\frac{{\sf diam}(D)}{\varepsilon}}\leq \frac{\gamma}{4},
\end{equation*}
All the choices above ensure that $A(D,M,K,d,\varepsilon)$ from \Cref{thm:main} satisfies
\begin{equation*}
A(D,M,K,d,\varepsilon)\le \frac{3\gamma}{4}.
\end{equation*}
Taking into account the above inequality and the estimate \eqref{eq:main1}, it is a dircet check to see that the right hand side of \eqref{eq:main1} is less than $\eta$ if $N$ satisfies \eqref{eq:N}, which concludes the proof.
\end{proof}

\subsection{Proofs for \Cref{ss:extension}}
\label{pf:6}

\begin{proof}[\textit{\textbf{Proof of \Cref{lem:extension}}}]
We have:
\begin{align*}
\frac{|G(x)-G(y)|}{|x-y|^\alpha}&\le \frac{|\psi\left(\frac{1}{\eps_0}r(x))\right)(g(\po(x))-g(\po(y))|}{|x-y|^\alpha}+\frac{|((\psi\left(\frac{1}{\eps_0}r(x))\right)-\psi\left(\frac{1}{\eps_0}r(y))\right))g(\po(y))|}{|x-y|^\alpha} \\
&\le \frac{|\po(x))-\po(y))|^\alpha |g|_\alpha}{|x-y|^\alpha}+\frac{|g|_{\infty}}{\eps_0}\frac{|r(x)-r(y)|}{|x-y|^\alpha}\\
&\le|\nabla\po|_{\infty}^\alpha |g|_\alpha+\frac{|g|_{L^\infty}}{\eps_0}\frac{|x-y|^\alpha |\sf{diam}(D)|^{1-\alpha}}{|x-y|^\alpha}.
\end{align*}
\end{proof}

\begin{proof}[\textit{\textbf{Proof of \Cref{lem:sizeonestep}}}]
First note that if $\phi=W^L\circ\sigma\circ\cdots\sigma W^1$, where $W^i$ is of the form $W^i(z)=A^iz+b^i$, then using the relation $x=\sigma(x)-\sigma(-x)$ we have
\[\phi(x)v=\begin{pmatrix} v& v \end{pmatrix}\sigma\circ\begin{pmatrix} W^L\\ -W^L \end{pmatrix} \circ\sigma\circ\cdots\sigma\circ W^1x, \quad x\in \mathbb{R}^d.
\]
The assertions follow now from \Cref{lem:augment} and \Cref{lem:additions_NNs}.
\end{proof}

\begin{proof}[\textit{\textbf{Proof of \Cref{cor:NNforGext}}}]
We first note that
\begin{equation*}
\left|\frac{1}{\eps_0}r-\frac{1}{\eps_0}\phi_{r}\right|_{\infty}\le\frac{\delta_r}{\eps_0}\quad \mbox{and} \quad
\left|\psi\left(\frac{1}{\eps_0}r\right)-\phi_{\psi}\left(\frac{1}{\eps_0}r\right)\right|_{\infty}\le \delta_\psi.
\end{equation*}
Therefore, by the triangle inequality we get

\begin{align*}
\left|\phi_{\psi}\left(\frac{r}{\eps_0}\right)(x)-\phi_{\psi}\left(\frac{\phi_{r}}{\eps_0}\right)(x)\right|&\le 2\delta_\psi+|\psi'|_{\infty}\left|\frac{1}{\eps_0}r(x)-\frac{1}{\eps_0}\phi_{r}(x)\right|\le 2\delta_\psi+|\psi'|_{\infty}\frac{\delta_r}{\eps_0},\\
\intertext{hence}
\left|\psi\left(\frac{1}{\eps_0}r\right)-\phi_{\psi}\left(\frac{\phi_{r}}{\eps_0}\right)\right|_{\infty}&\le 3\delta_\psi+\frac{\delta_d}{\eps_0}.\\
\intertext{Reasoning analogously we get}
|g(\po)-\phi_{g}\circ\phi_{\pi}|_{\infty}&\le 3\delta_g+|\nabla g|_{\infty}\delta_\pi,
\intertext{so again by the triangle inequality}
\left|\psi\left(\frac{1}{\eps_0}r\right)g(\po)-\phi_{\psi}\left(\frac{\phi_{r}}{\eps_0}\right)\phi_{g}\circ\phi_{\pi}\right|_\infty&\leq \left(3\delta_\psi+\frac{\delta_d}{\eps_0}\right)|g|_\infty+ \left(3\delta_g+|\nabla g|_{\infty}\delta_\pi \right)(\delta_\psi+1). 
\end{align*}
Let now $\Pi$ be the ReLU DNN given in \Cref{coro:product_NNs} with $\epsilon_p:=\left(3\delta_\psi+\frac{\delta_d}{\eps_0}\right)|g|_\infty+ \left(3\delta_g+|\nabla g|_{\infty}\delta_\pi \right)(\delta_\psi+1)=\overline{\delta}/2$, and $c:=|g|_\infty$, so that
\begin{equation*}
    \left|G-\Pi\left(\phi_{\psi}\left(\frac{\phi_{r}}{\eps_0}\right),\phi_{g}\circ\phi_{\pi}\right)\right|_\infty\leq\overline{\delta}.
\end{equation*}
Now, by taking $\phi_G:=\Pi\left(\phi_{\psi}\left(\frac{\phi_{r}}{\eps_0}\right),\phi_{g}\circ\phi_{\pi}\right)$, the statement follows directly from \Cref{coro:product_NNs} and \Cref{lem:composition_NNs}.
\end{proof}

\begin{proof}[\textit{\textbf{Proof of \Cref{thm:mainNN}}}]

i) The fact that $\widetilde{u}_M^N(\cdot)(\omega)$ can be realized, for each $\omega \in \Omega$, as a ReLU DNN, denoted in the statement by $\mathbb{U}_M^N(\omega,\cdot)$, is a direct consequence of \Cref{lem:sizeWoS}, \Cref{coro:product_NNs}, \Cref{lem:composition_NNs}, and \Cref{lem:additions_NNs}.

Next, let $u_M^N$ be given by \eqref{eq:MCestimator} and consider the modification of $u_M^N$ by replacing $g$ and $f$ with $\phi_g$ and $\phi_f$, namely
\begin{equation*}
    \widehat{u}_{M}^{N}(x):=\frac{1}{N}\sum_{i=1}^N\left[ \phi_g(X^{x,i}_{M})+\frac{1}{d}\sum\limits_{k=1}^M \widetilde{r}^2(X^{x,i}_{k-1})\phi_f\left(X^{x,i}_{k-1}+\widetilde{r}(X^{x,i}_{k-1})Y^i\right)\right], x\in D.
\end{equation*}
Then using assumption $a.4)$ we can easily deduce that
$|\widehat{u}_{M}^{N}-\widetilde{u}_{M}^{N}|_\infty\leq \frac{M\epsilon_p}{d}(1+2|f|_\infty)\leq \gamma/6$.
Therefore,
\begin{align*}
    \sup\limits_{x\in D}\left | u(x)-\mathbb{U}_M^{ N}(\cdot,x)\right|&=|u-\widetilde{u}_M^N|_\infty\leq |u-u_M^N|_\infty +|u_M^N-\widehat{u}_M^N|+|\widehat{u}_M^N-\widetilde{u}_M^N|\\
    &\leq |u-u_M^N|_\infty+\epsilon_g+\frac{M{\rm diam}(D)^2\epsilon_f}{d}+\gamma/6\\
    &\leq |u-u_M^N|_\infty+\gamma/2,
\end{align*}
where the last inequality follows from assumptions $a.1)$, $a.2)$.
Consequently,
\begin{align*}
    \mathbb{P}\left( \sup_{x\in D} \left | u(x)-\mathbb{U}_M^{ N}(\cdot,x)\right| \geq \gamma \right)&\leq \mathbb{P}\left(|u-u_M^N|_\infty \geq \gamma/2 \right), 
\end{align*}
so, we can employ \Cref{coro:noKepsilon} to conclude the first assertion.

Let us proceed at proving ii).
First of all, note that
\[
\mathcal{L}(\widetilde{r})=\mathcal{L}(\phi_r)+1,\;\mathcal{W}(\widetilde{r})=\mathcal{W}(\phi_r),\;{\rm size}(\widetilde{r})\leq{\rm size}(\phi_r)+2.
\]
Then, by \Cref{lem:composition_NNs} and \Cref{lem:sizeWoS} we have
\begin{align}
    \mathcal{L}(\phi_g(X^{\cdot,i}_{M}))&= \mathcal{L}(\phi_g)+\mathcal{L}(X^{\cdot,i}_{M})=\mathcal{L}(\phi_g)+M(\mathcal{L}(\widetilde{r})+1))+1=\mathcal{L}(\phi_g)+M(\mathcal{L}(\phi_r)+2))+1,\nonumber\\
    \mathcal{W}(\phi_g(X^{\cdot,i}_{M}))&\leq \max(\mathcal{W}(\phi_g),\mathcal{W}(X^{\cdot,i}_{M}),2d)\leq \max(\mathcal{W}(\phi_g),2d+\max(d,\mathcal{W}(\phi_r))),\nonumber\\
    {\rm size}(\phi_g(X^{\cdot,i}_{M}))&\leq 2{\rm size}(\phi_g)+2{\rm size}(X^{\cdot,i}_{M})\leq 2{\rm size}(\phi_g)+ 4dM[4d+\mathcal{W}(\widetilde{r})+\mathcal{L}(\widetilde{r})+2]+d+2M{\rm size}(\widetilde{r})\nonumber\\
    &\leq 2{\rm size}(\phi_g)+ 4dM[4d+\mathcal{W}(\phi_r)+\mathcal{L}(\phi_r)+3]+2d+4M[{\rm size}(\phi_r)+2]\nonumber\\
    &\in \mathcal{O}({\rm size}(\phi_g)+M{\rm size}(\phi_r)+Md\max(d,\mathcal{W}(\phi_r),\mathcal{L}(\phi_r))).\label{eq:size1}
\end{align}
Further, for each $k\geq 0$, by \Cref{lem:composition_NNs} we get
\begin{align}
    {\rm size}\left(\phi_f\left(X^{\cdot,i}_{k}+\widetilde{r}(X^{\cdot,i}_{k})Y^i\right)\right)&\leq 2 {\rm size}\left(\phi_f\right)+2{\rm size}\left(X^{\cdot,i}_{k}+\widetilde{r}(X^{\cdot,i}_{k})Y^i\right)\nonumber\\
    \intertext{and since $X^{\cdot,i}_{k}+\widetilde{r}(X^{\cdot,i}_{k})Y^i$ has the same size as $X^{\cdot,i}_{k+1}$, we can continue with}
    &=2 {\rm size}\left(\phi_f\right)+2{\rm size}\left(X^{\cdot,i}_{k+1}\right)\nonumber\\
    &\leq 2{\rm size}(\phi_f)+ 4d(k+1)[4d+\mathcal{W}(\phi_r)+\mathcal{L}(\phi_r)+3]+2d+4(k+1)[{\rm size}(\phi_r)+2]. \label{eq:size2}
\end{align}

The next step is to use \Cref{coro:product_NNs} to get
\begin{align}
    {\rm size}\left(\Pi\left(\widetilde{r}(X^{\cdot,i}_{k}),\widetilde{r}(X^{\cdot,i}_{k})\right)\right)&\leq 8 {\rm size}(\widetilde{r}(X^{\cdot,i}_{k}))+\mathcal{O}(\lceil\log (\epsilon_p^{-1})+\log(c)\rceil)\nonumber\\
    &\leq 16{\rm size}(\widetilde{r})+ 16{\rm size}(X^{\cdot,i}_{k})+\mathcal{O}(\lceil\log (\epsilon_p^{-1})+\log(c)\rceil)\nonumber\\
    &\leq 16{\rm size}(\phi_r)+ 32dk[4d+\mathcal{W}(\phi_r)+\mathcal{L}(\phi_r)+3]+16d+32k[{\rm size}(\phi_r)+2]\nonumber\\
    &\phantom{\leq 16{\rm size}(\phi_r)\;}+\mathcal{O}(\lceil\log (\epsilon_p^{-1})+\log(c)\rceil),\label{eq:size3}
\end{align}
so that by \eqref{eq:size2} and \eqref{eq:size3} together with \Cref{coro:product_NNs} we obtain

\begin{align}
{\rm size}&\left(\Pi\left(\Pi\left(\widetilde{r}(X^{\cdot,i}_{k}),\widetilde{r}(X^{\cdot,i}_{k})\right),\phi_f\left(X^{\cdot,i}_{k}+\widetilde{r}(X^{\cdot,i}_{k})Y^i\right)\right)\right)\nonumber\\
&\leq 4{\rm size}\left(\Pi\left(\widetilde{r}(X^{\cdot,i}_{k}),\widetilde{r}(X^{\cdot,i}_{k})\right)\right)+4{\rm size}\left(\phi_f\left(X^{\cdot,i}_{k}+\widetilde{r}(X^{\cdot,i}_{k})Y^i\right)\right)+\mathcal{O}(\lceil\log (\epsilon_p^{-1})+\log(c)\rceil)\nonumber\\
&\leq 24{\rm size}(\phi_r)+ 128dk[4d+\mathcal{W}(\phi_r)+\mathcal{L}(\phi_r)+3]+64d+128k[{\rm size}(\phi_r)+2]\nonumber\\
&\phantom{\leq 24{\rm size}(\phi_r)\;}+8{\rm size}(\phi_f)+ 16d(k+1)[4d+\mathcal{W}(\phi_r)+\mathcal{L}(\phi_r)+3]+8d+16(k+1)[{\rm size}(\phi_r)+2]\\
&\phantom{\leq 24{\rm size}(\phi_r)\;}+\mathcal{O}(\lceil\log (\epsilon_p^{-1})+\log(c)\rceil)\nonumber\\
&\in \mathcal{O}\left(dk\max(d,\mathcal{W}(\phi_r),\mathcal{L}(\phi_r))+k{\rm size}(\phi_r)+{\rm size}(\phi_f)+\lceil\log (\epsilon_p^{-1})+\log(c)\rceil\right). \label{eq:size4}
\end{align}

Finally, corroborating \eqref{eq:size1} with \eqref{eq:size4} and by applying \Cref{lem:additions_NNs}, we obtain for $\omega\in \Omega$ that
\begin{align*}
    {\rm size}&(\mathbb{U}_M^{ N}(\omega,\cdot))\\
    &\in \mathcal{O}\left(MN\left[dM\max(d,\mathcal{W}(\phi_r),\mathcal{L}(\phi_r))+M{\rm size}(\phi_r)+{\rm size}(\phi_f)+{\rm size}(\phi_g)+\lceil\log (\epsilon_p^{-1})+\log(c)\rceil\right]\right),
\end{align*}
and by assumption $a.4)$ we deduce that
\[
{\rm size}(\mathbb{U}_M^{ N}(\omega,\cdot))\in \mathcal{O}\left(MN\left[dM\max(d,\mathcal{W}(\phi_r),\mathcal{L}(\phi_r))+M{\rm size}(\phi_r)+{\rm size}(\phi_f)+\log \left(\frac{1}{\gamma d}\right)\right]\right),
\]
where the tacit constant depends on $\max({\rm diam}(D),|f|_\infty)$.

{
In particular,
\begin{align*}
    M&\in \mathcal{O}\left(d^2\gamma^{-4/\alpha}\log\left(\frac{1}{\gamma}\right)\right), \\
    N&\in \mathcal{O}\left(d^2\gamma^{-8/\alpha-2}\log^2\left(\frac{1}{\gamma}\right)\left[d^3\gamma^{-4/\alpha}\log\left(\frac{1}{\gamma}\right)+\log\left(\frac{1}{\eta}\right)\right]\right), \\
    {\rm size}(\mathbb{U}_M^{ N}(\omega,\cdot))&\in \mathcal{O}\left(d^7\gamma^{-16/\alpha-4}\log^4\left(\frac{1}{\gamma}\right)\left[d^3\gamma^{-4/\alpha}\log\left(\frac{1}{\gamma}\right)+\log\left(\frac{1}{\eta}\right)\right]{\rm S}\right),
\end{align*}
where \[{\rm S}:=\left[\max(d,\mathcal{W}(\phi_r),\mathcal{L}(\phi_r))+{\rm size}(\phi_r)+{\rm size}(\phi_g)+{\rm size}(\phi_f)\right]\]
and the tacit constant depends on $|g|_\alpha,|g|_\infty,|f|_\infty,{\rm diam}(D),{\rm adiam}(D),\delta,\alpha,\log(2+|\phi_r|_1)$.

Now, if $D$ is $\delta$-defective convex, then
\begin{align*}
    M&\in \mathcal{O}\left(d\log\left(\frac{d}{\gamma}\right)\right), \\
    N&\in \mathcal{O}\left(\frac{1}{\gamma^2}\log^2\left(\frac{d}{\gamma}\right)\left[d^2\log\left(\frac{d}{\gamma}\right)\log(2+|\phi_r|_1)+\log\left(\frac{1}{\eta}\right)\right]\right), \\
    {\rm size}(\mathbb{U}_M^{ N}(\omega,\cdot))&\in \mathcal{O}\left(\frac{d^3}{\gamma^2}\log^4\left(\frac{d}{\gamma}\right)\left[d^2\log\left(\frac{d}{\gamma}\right)+\log\left(\frac{1}{\eta}\right)\right]{\rm S}\right),
\end{align*}
where ${\rm S}$ is as above and the tacit constant depends on $|g|_\alpha,|g|_\infty,|f|_\infty,{\rm diam}(D),{\rm adiam}(D),\delta,\alpha,\log(2+|\phi_r|_1)$.
}
\end{proof}

{
\section{Numerical results}
\label{S:numerics}

Throughout this section we numerically test some key theoretical results obtained in this paper. In order to take advantage of the highly parallelizable properties of the proposed Monte Carlo methods, we have implemented the algorithms using GPU parallel computing within the PyTorch framework. 

Concretely, we consider two types of domains, namely:
\begin{equation*}
 D_{\sf c}:=[-1,1]^d, d\geq 1 \quad \mbox{ and } \quad D_{\sf ac}:= [-1,1]^d \setminus \left\{x\in \mathbb{R}^d : |x|_1:=|x_1|+\dots+|x_d|\leq 0.5 \right\},
\end{equation*}
which are illustrated by \Cref{fig:D} for the case $d=2$.

\medskip 
\begin{figure}[H]

\begin{subfigure}{0.5\textwidth}
\centering
\includegraphics[scale=0.5]{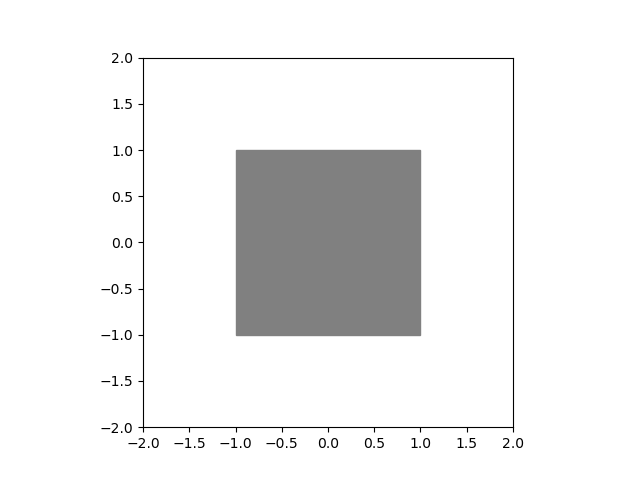}
\caption{$D_{\sf c}=[-1,1]^2$}
\label{fig:Hypercube}
\end{subfigure}
\begin{subfigure}{0.5\textwidth}
\centering
\includegraphics[scale=0.5]{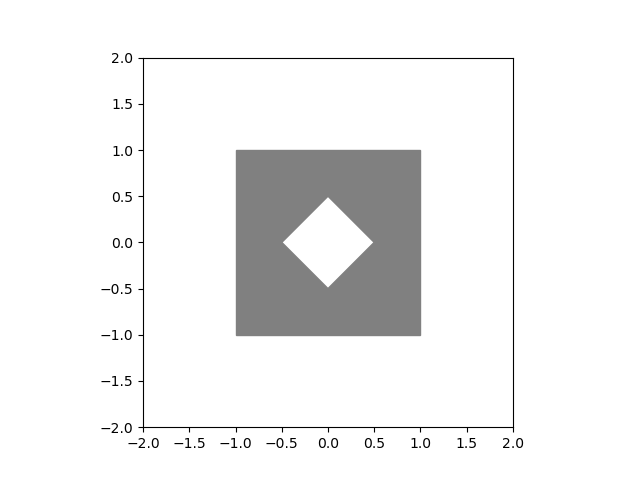}
\caption{$D_{\sf ac}=\left\{x\in [-1,1]^2 : |x_1|+|x_2|\leq 0.5 \right\}$}
\label{fig:AnnularH}
\end{subfigure}
\caption{$D_{\sf c}$ and $D_{\sf ac}$ for $d=2$}
\label{fig:D}
\end{figure}
  
We have performed the following numerical tests by:

\medskip
\noindent{\bf Test 1.} Here we numerically verify the estimate \eqref{eq:exptail}, namely that 
\begin{equation} \label{eq:test 1}
\mathbb{P}(r(X^{x}_M)\geq \varepsilon)\leq \left(1-\frac{\beta^2(1-\delta)}{4d}\right)^M \sqrt{\frac{r(x)}{\varepsilon}}, \quad x\in D, 
\end{equation}
by fixing $\varepsilon$ and varying $d$ and $M$ for both $D_{\sf c}$ and $D_{\sf ac}$. 
Note that $D_{\sf c}$ is $0$-defective convex since it is in fact convex, so \eqref{eq:test 1} is valid with $\delta=0$ by \Cref{prop:logepsilon}. Furthermore, for this test, the WoS is constructed with the exact distance to the boundary, so we shall take $\beta=1$.
On the other hand, $D=D_{\sf ac}$ is not defective convex since if we denote the boundary corner $(0.5,0,0,\dots,0)\in D_{\sf ac}$ by $x_0$, then $\Delta r(x)=(d-1)|x-x_0|^{-1}$ is unbounded on the cone
$\left\{(r,x)\in D_{\sf ac} : r\in (0.5, 0.75), |x|\leq r\right\}$, hence \eqref{eq:dconvex} can not hold.
Nevertheless, for $\varepsilon$ fixed we speculated in \Cref{rem:any D is defective} that a similar asymptotic behavior as in \eqref{eq:test 1} still holds with respect to $M$ and $d$; this is numerically tested in \Cref{fig:N_steps_M_ac} below. 

Let us introduce the notations
\begin{align*}
    &{\sf U}_{\sf bound}(d,M,x,\varepsilon):=\left(1-\frac{1}{4d}\right)^M \sqrt{\frac{r(x)}{\varepsilon}},\\imagini
    &\mathbb{P}_N(d,M,x,\varepsilon):=\frac{1}{N}\sum\limits_{1\leq i \leq N} 1_{\left[r\geq \varepsilon\right]}(X_M^{x,i}), \quad \mbox{ where } X_M^{x,1},\cdots, X_M^{x,N} \sim X_M^{x} \mbox{ are independent}, N\geq 1,
\end{align*}
so that the inequality to be tested becomes
\begin{equation}\label{eq:Test 1}
\mathbb{P}_N(d,M,x,\varepsilon)\approx\mathbb{P}(r(X^{x}_M)\geq \varepsilon) \leq {\sf U}_{\sf bound}(d,M,x,\varepsilon) \quad \mbox{ 
 for } N \mbox{ sufficiently large}.
\end{equation}
We obtained the following results:

\medskip 
\begin{figure}[H]

\begin{subfigure}{0.5\textwidth}
\centering
\includegraphics[scale=0.4]{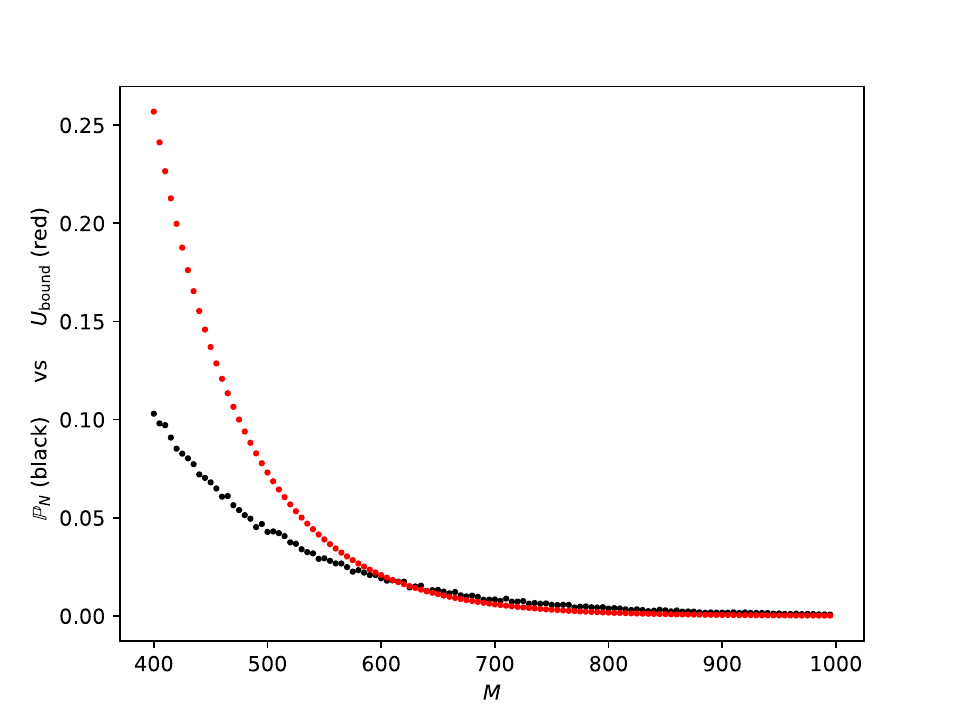}
\caption{$M=5\times i, i\in \{80, \dots, 200\}$}
\label{fig:N_steps}
\end{subfigure}
\begin{subfigure}{0.5\textwidth}
\centering
\includegraphics[scale=0.4]{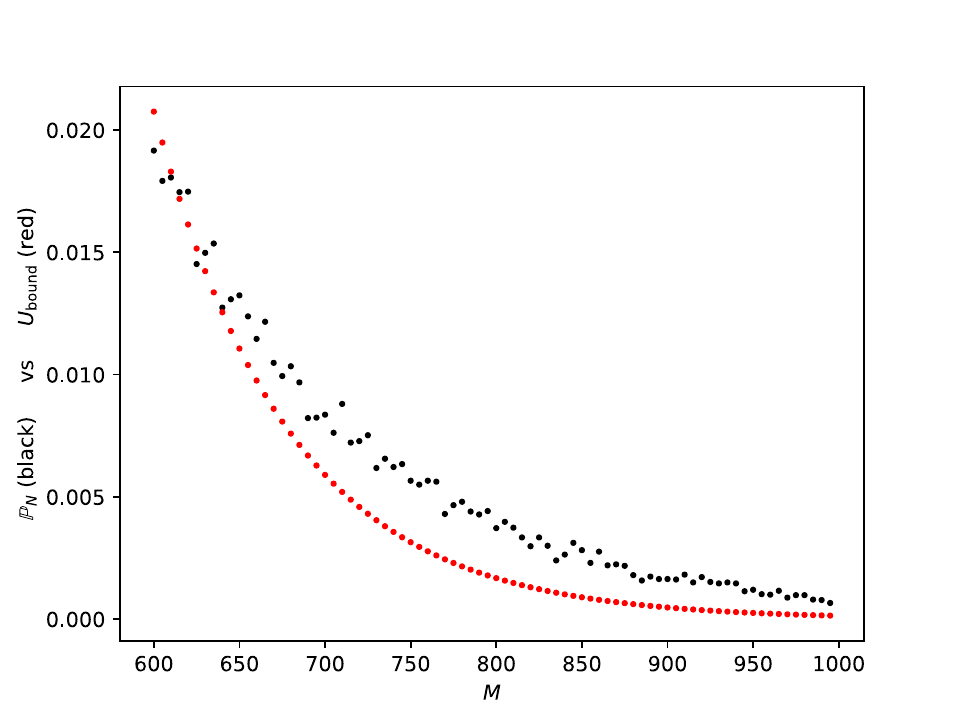}
\caption{$M=5\times i, i\in \{120, \dots, 200\}$}
\label{fig:N_steps_zoom}
\end{subfigure}
\caption{${\sf U}_{\sf bound}(d,\cdot,x_0,\varepsilon)$ (red) vs $\mathbb{P}_N(d,\cdot,x_0,\varepsilon)$ (black) for $D=D_{\sf c}$, $d=20,\varepsilon=10^{-3}, N=5\times 10^4$, whilst $x_0$ is arbitrarily chosen in $D$ such that $|x_0|=0.5$.}
\label{fig:N_steps_M}
\end{figure}

\begin{figure}[H]

\begin{subfigure}{0.5\textwidth}
\centering
\includegraphics[scale=0.4]{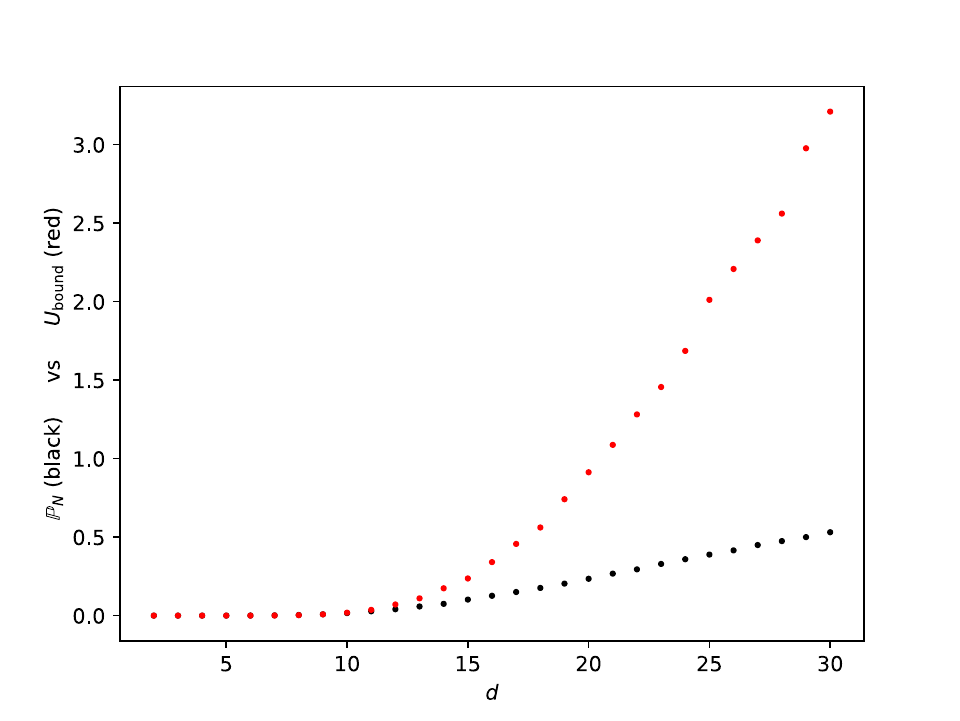}
\caption{$d\in \{2, \dots, 30\}$, $M=300$}
\label{fig:N_steps_D}
\end{subfigure}
\begin{subfigure}{0.5\textwidth}
\centering
\includegraphics[scale=0.4]{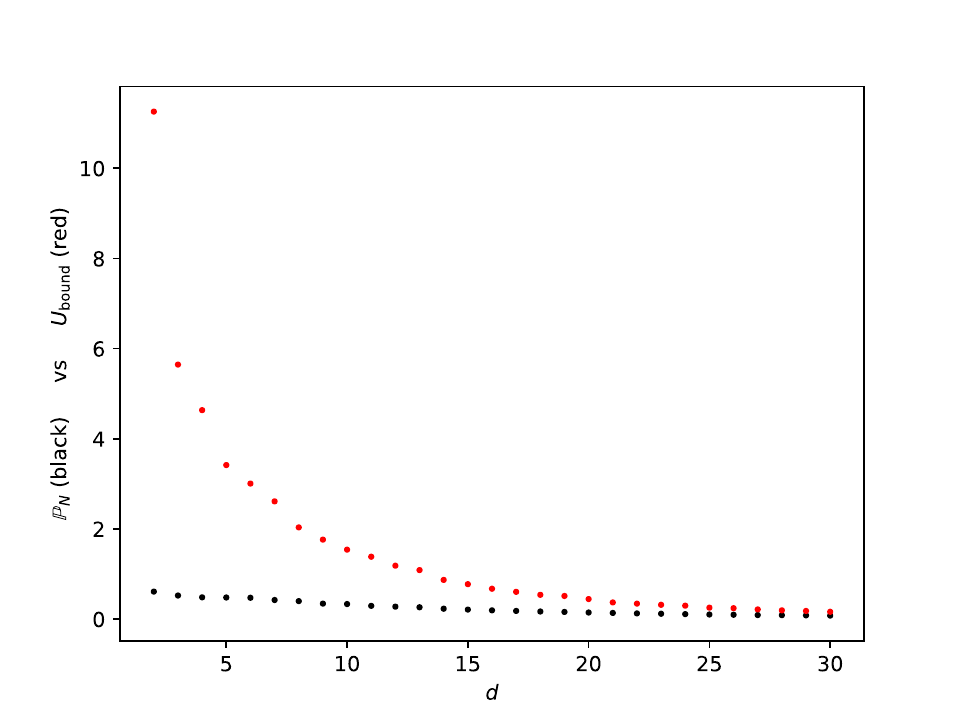}
\caption{$d\in \{2, \dots, 30\}$, $M=4 \times d^{1+0.5}$}
\label{fig:D and N_steps}
\end{subfigure}
\caption{${\sf U}_{\sf bound}(\cdot,\cdot,x_0,\varepsilon)$ (red) vs $\mathbb{P}_N(\cdot,\cdot,x_0,\varepsilon)$ (black) for $D=D_{\sf c}$, $\varepsilon=10^{-3}, N=5\times 10^4$, whilst $x_0$ is arbitrarily chosen in $D$ such that $|x_0|=0.5$, for each dimension $d$.}
\label{fig:N_steps_d}
\end{figure}

\begin{figure}[H]

\begin{subfigure}{0.5\textwidth}
\centering
\includegraphics[scale=0.4]{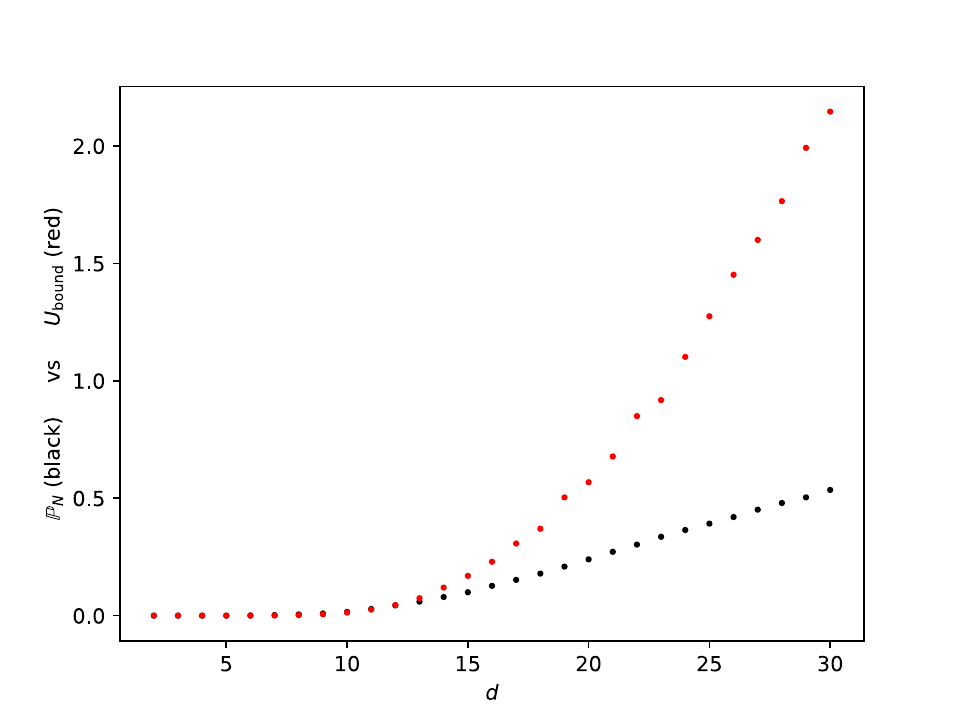}
\caption{$d\in \{2, \dots, 30\}$, $M=300$}
\label{fig:N_steps_ac}
\end{subfigure}
\begin{subfigure}{0.5\textwidth}
\centering
\includegraphics[scale=0.4]{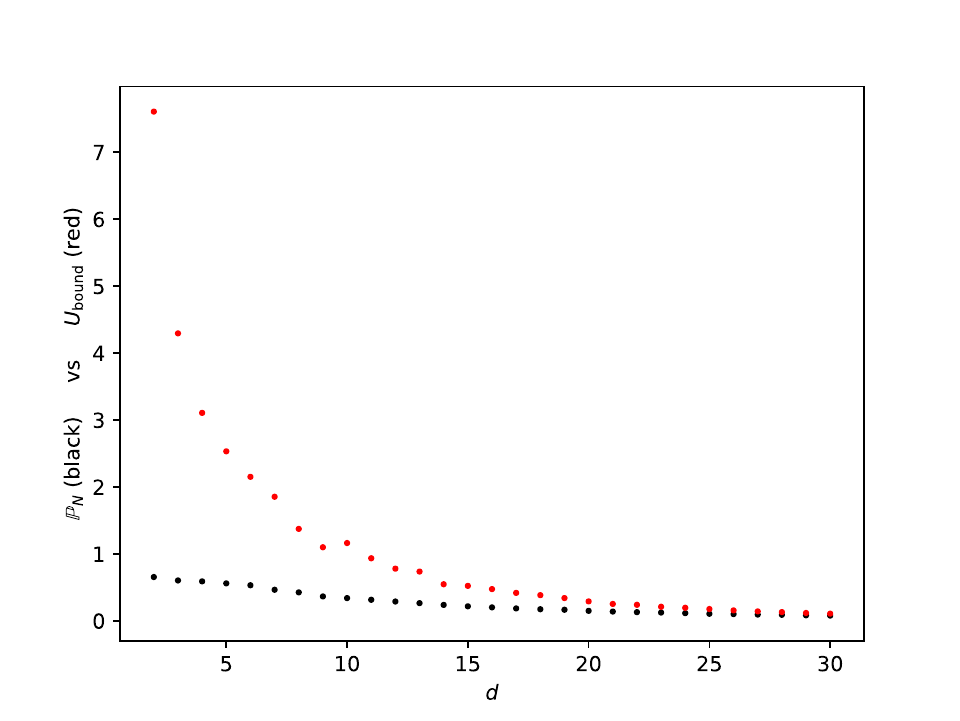}
\caption{$d\in \{2, \dots, 30\}$, $M=4 \times d^{1+0.5}$}
\label{fig:N_steps_M_ac}
\end{subfigure}
\caption{${\sf U}_{\sf bound}(\cdot,\cdot,x_0,\varepsilon)$ (red) vs $\mathbb{P}_N(\cdot,\cdot,x_0,\varepsilon)$ (black) for $D=D_{\sf ac}$, $\varepsilon=10^{-3}$, $N=5\times 10^4$, whilst $x_0$ is arbitrarily chosen in $D$ such that $|x_0|=0.7$, for each dimension $d$.}
\label{fig:N_steps_d_ac}
\end{figure}

{\bf Comments on Test 1.} 
\begin{enumerate}
    \item[$\bullet$] The numerical results depicted in \Cref{fig:N_steps_M} and \Cref{fig:N_steps_d} validate the upper bound \eqref{eq:Test 1}. Note that in \Cref{fig:N_steps_zoom}, namely for large values of $M$, inequality \eqref{eq:Test 1} becomes quite sharp, yet slightly reversed; the apparent reversed inequality is just a consequence of the Monte Carlo error that appears in the approximation $\mathbb{P}_N(d,M,x,\varepsilon)\approx\mathbb{P}(r(X^{x}_M)\geq \varepsilon)$ corresponding to $N=5\times 10^4$ Monte Carlo samples.
    \item[$\bullet$] The numerical evidence illustrated in \Cref{fig:N_steps_d_ac} shows that at least for one relevant example of non $\delta$-defective convex domain, namely for $D_{\sf ac}$, by fixing $\varepsilon$, the estimate \eqref{eq:Test 1} is still in force; in particular, this sustains the idea expressed in \Cref{rem:any D is defective}.
    \item[$\bullet$] Finally, recall that our main estimates obtained in the previous sections in the case of defective convex domains require $M\in \mathcal{O}(d\log(d/\gamma))$. This requirement is suggested also by comparing \Cref{fig:N_steps_D} versus \Cref{fig:D and N_steps} as well as  \Cref{fig:N_steps_ac} versus \Cref{fig:N_steps_M_ac}; in \Cref{fig:D and N_steps} and \Cref{fig:N_steps_M_ac} we have chosen $M\in \mathcal{O}(d^{1+0.5})$ instead of $M\in \mathcal{O}(d\log(d/\gamma))$ only because in the latter case it is harder to nicely visualize the numerical results. 
\end{enumerate}

\medskip
\noindent{\bf Test 2.} 
Here the goal is to test the approximation of the solution $u$ to \eqref{e:0} by simulating its Monte Carlo estimator \eqref{e:5}. 
For simplicity we take the source term $f=1$, hence we deal with
\begin{equation}\label{pde:numerics}
\begin{cases}
\frac{1}{2}\Delta u=-1 \,\mbox{ in } D \\
u|_{\partial D}=g, 
\end{cases}, \mbox{ where } D=D_{\sf ac}.
\end{equation}
In order to validate the numerical results, we consider a particular explicit solution to \eqref{pde:numerics}, namely
\begin{align}\label{eq:exact_u}
    &u(x)= x_1^2+\dots +x_k^2 - x_{k+1}^2-\dots -x_d^2-x_1^2, \quad d=2k, x=(x_1,\cdots,x_d)\in \mathbb{R}^d,\\
    &g=u|_{\partial D}, \quad D=D_{\sf c} \mbox{ or } D_{\sf ac}. 
\end{align}
Further, we need to introduce some notations. 
For $M,N,L,E \geq 1$ and $W_i,1\leq i\leq L$ independent and uniformly distributed on $D$, we introduce the notations
\begin{align*}
    &{\sf Err}_{M,N}(x):=\left | u(x)-u_M^{N}(x)\right|, x\in D, \quad  |{\sf Err}_{M,N}|_{L^1(D/|D|)}:= |D|^{-1}\int_D {\sf Err}_{M,N}(x) \;dx\\
    & |{\sf Err}_{M,N,L}|_{L^1(D/|D|)}:=1/L\sum_{1\leq i\leq L} {\sf Err}_{M,N}^{(i)}(W_i), \quad |{\sf Err}_{M,N,L}|_{L^\infty(D)}:=\max_{1\leq i\leq L} {\sf Err}_{M,N}^{(i)}(W_i)\\ 
    \intertext{where ${\sf Err}_{M,N}^{(i)}(\cdot), 1\leq i\leq L$ are iid copies of ${\sf Err}_{M,N}(\cdot)$, independent of $(W_i)_{1\leq i\leq L}$,}
    & |{\sf Err}_{M,N,L,E}|_{L^\infty(D)}:=1/E\sum_{1\leq j\leq E}|{\sf Err}_{M,N,L}^{(j)}|_{L^\infty(D)},
\end{align*}
where $|{\sf Err}_{M,N,L}^{(j)}|_{L^\infty(D)}, 1\leq j\leq E$ are independent copies of $|{\sf Err}_{M,N,L}|_{L^\infty(D)}$.
In particular, by the law of large numbers we immediately have
\begin{equation}\label{eq:test21}
    \lim_{L\to \infty}|{\sf Err}_{M,N,L}|_{L^1(D/|D|)} =\mathbb{E}\left\{|{\sf Err}_{M,N}|_{L^1(D/|D|)}\right\} \mbox{ almost surely},
\end{equation}
whilst by a similar argument as in \cite[Lemma 4.3]{BeBeGrJe21} one can show that 
\begin{equation}\label{eq:test22}
    \lim_{E\to \infty}\lim_{L\to \infty} |{\sf Err}_{M,N,L,E}|_{L^\infty(D)} =\mathbb{E}\left\{\sup_{x\in D}|{\sf Err}_{M,N}(x)|\right\} \mbox{ almost surely}.
\end{equation}

The aim of this test is to numerically validate the approximation $u\approx u_{M}^N$ by computing the mean errors $\mathbb{E}\left\{|{\sf Err}_{M,N}|_{L^1(D/|D|)}\right\} $ and $\mathbb{E}\left\{\sup_{x\in D}|{\sf Err}_{M,N}(x)|\right\}$. 
To this end, justified by \eqref{eq:test21} and \eqref{eq:test22}, we shall simulate $|{\sf Err}_{M,N,L}|_{L^1(D/|D|)}(\omega)$ and $|{\sf Err}_{M,N,L,E}|_{L^\infty(D)}(\omega)$ for $L$ large, e.g. $L=1000$ or $L=2000$. In order to reduce the computational burden, the value of $E$ is going to be taken relatively small, e.g. $E=5$ or $E=10$. However, we point out that choosing a small value for $E$ shall not alter the reliability of the Monte Carlo estimate, mainly due to the fact that one can show that the distribution of $|{\sf Err}_{M,N,L}|_{L^\infty(D)}$ is concentrated.

\begin{figure}[H]

\begin{subfigure}{0.5\textwidth}
\centering
\includegraphics[scale=0.4]{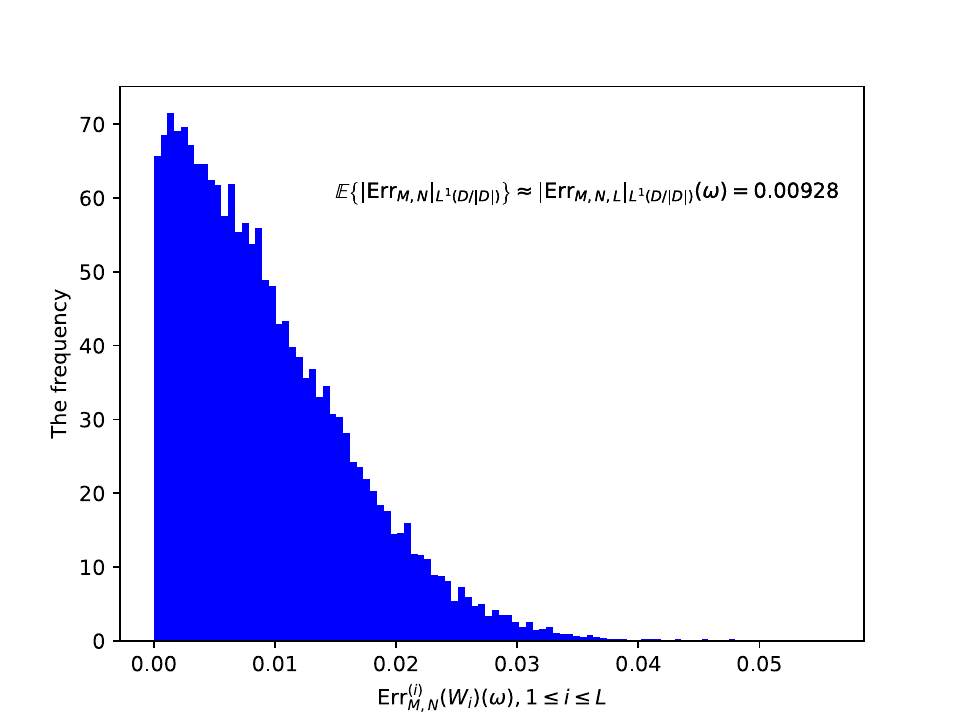}
\caption{The histogram of ${\sf Err}_{M,N}^{(i)}(\omega), 1\leq i\leq L=2\times 10^4$}
\label{fig:Err_MNL_d=10}
\end{subfigure}
\begin{subfigure}{0.5\textwidth}
\centering
\includegraphics[scale=0.4]{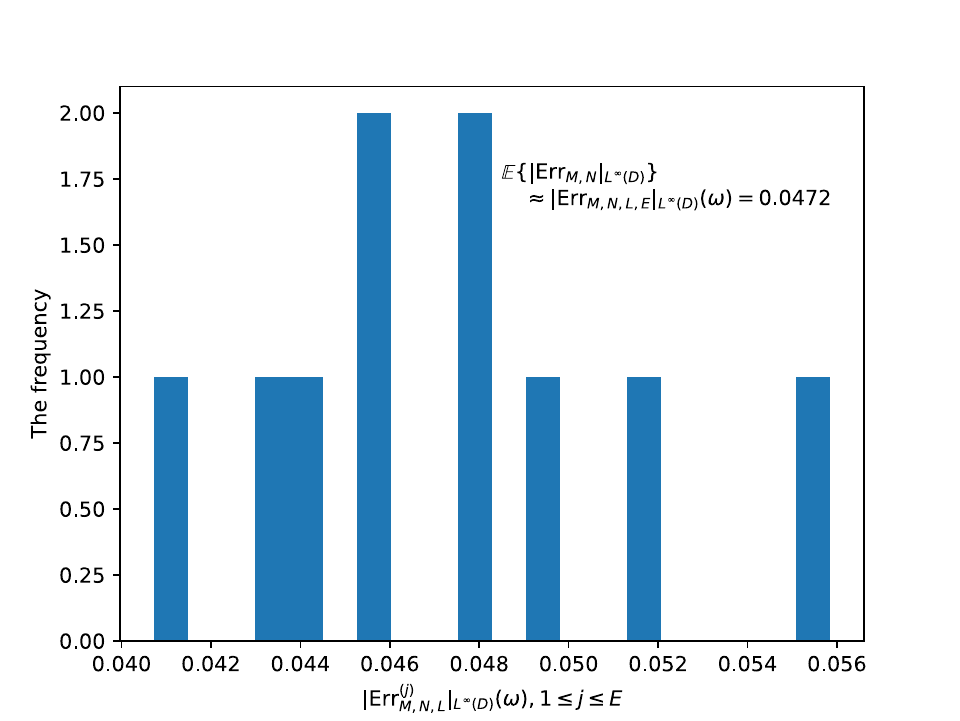}
\caption{The histogram of $|{\sf Err}_{M,N,L}^{(j)}|_{L^\infty(D)}, 1\leq j\leq E=10$, $L=2\times 10^3$}
\label{fig:Err_MNL_inf_d=10}
\end{subfigure}
\caption{ The histograms of ${\sf Err}_{M,N}^{(i)}, 1\leq i\leq L$ and $|{\sf Err}_{M,N,L}^{(j)}|_{L^\infty(D)}, 1\leq j\leq E$ for $D=D_{\sf ac}$, $d=10$, $N=10^5$, $M=500$, whilst $L$ and $E$ are specified for each subfigure, accordingly.}
\label{fig:N_steps_d=10_ac}
\end{figure}

\begin{figure}[H]

\begin{subfigure}{0.5\textwidth}
\centering
\includegraphics[scale=0.4]{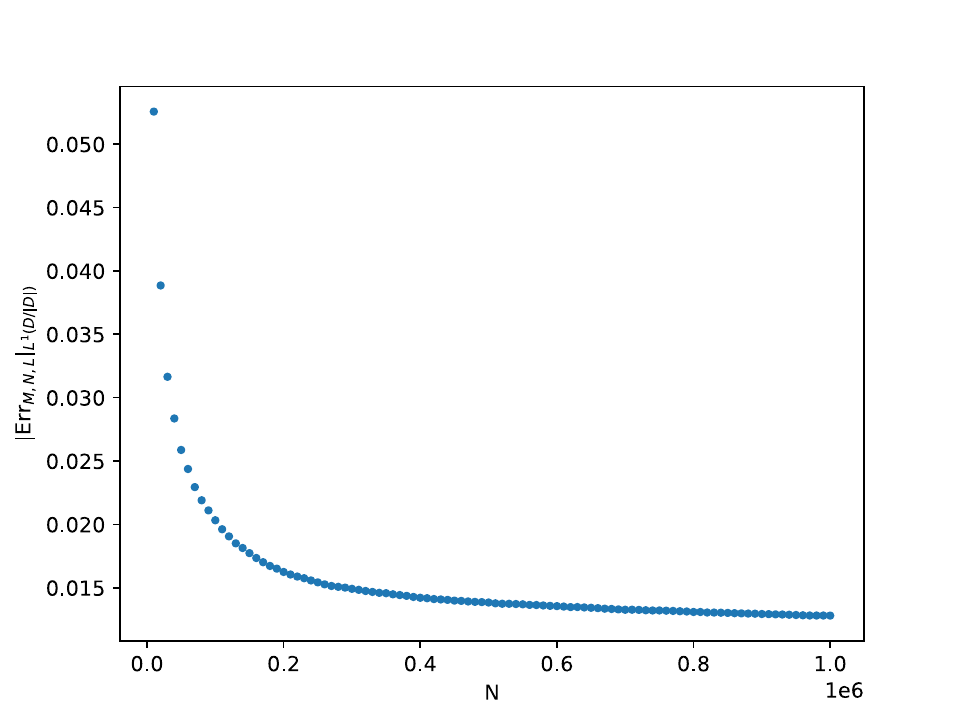}
\caption{The evolution of $|{\sf Err}_{M,N,L}|_{L^1(D/|D|)}(\omega)$ w.r.t. $N$}
\label{fig:Err_MNL_l1_d=100}
\end{subfigure}
\begin{subfigure}{0.5\textwidth}
\centering
\includegraphics[scale=0.4]{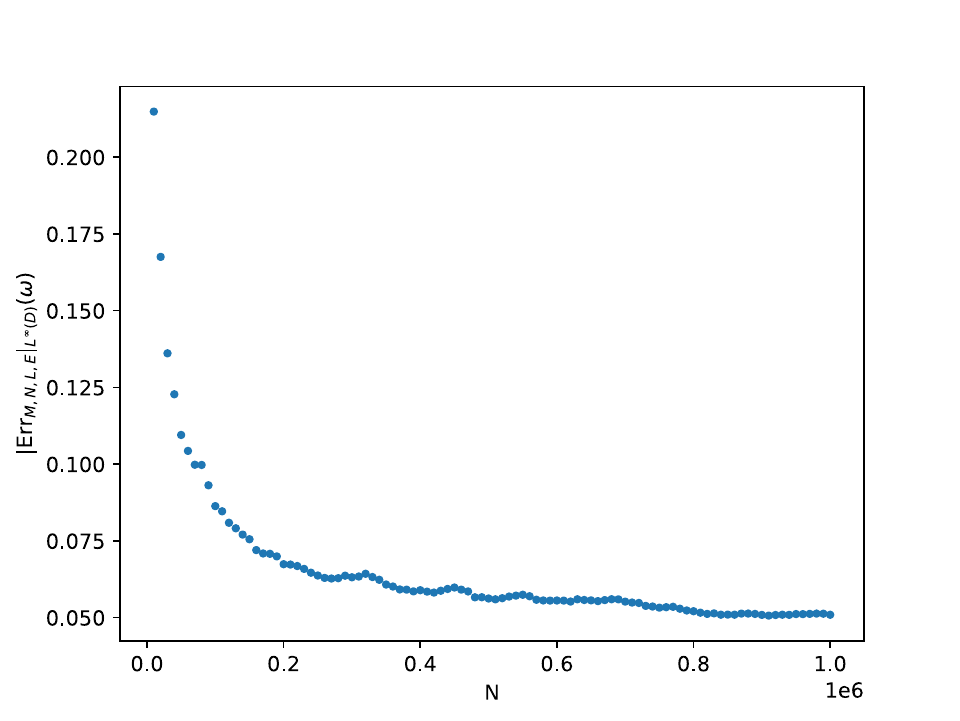}
\caption{The evolution of $|{\sf Err}_{M,N,L,E}|_{L^\infty(D)}(\omega)$ w.r.t. $N$}
\label{fig:Err_MNLE_inf_d=100}
\end{subfigure}
\caption{ The evolution of $|{\sf Err}_{M,N,L}|_{L^1(D/|D|)}(\omega)$ and $|{\sf Err}_{M,N,L,E}|_{L^\infty(D)}(\omega)$ w.r.t. $N$, for $D=D_{\sf ac}$, $d=100$, $M=500$, $L=1000$, and $E=5$.}
\label{fig:N_steps_d=100_ac}
\end{figure}

\begin{figure}[H]

\begin{subfigure}{0.5\textwidth}
\centering
\includegraphics[scale=0.4]{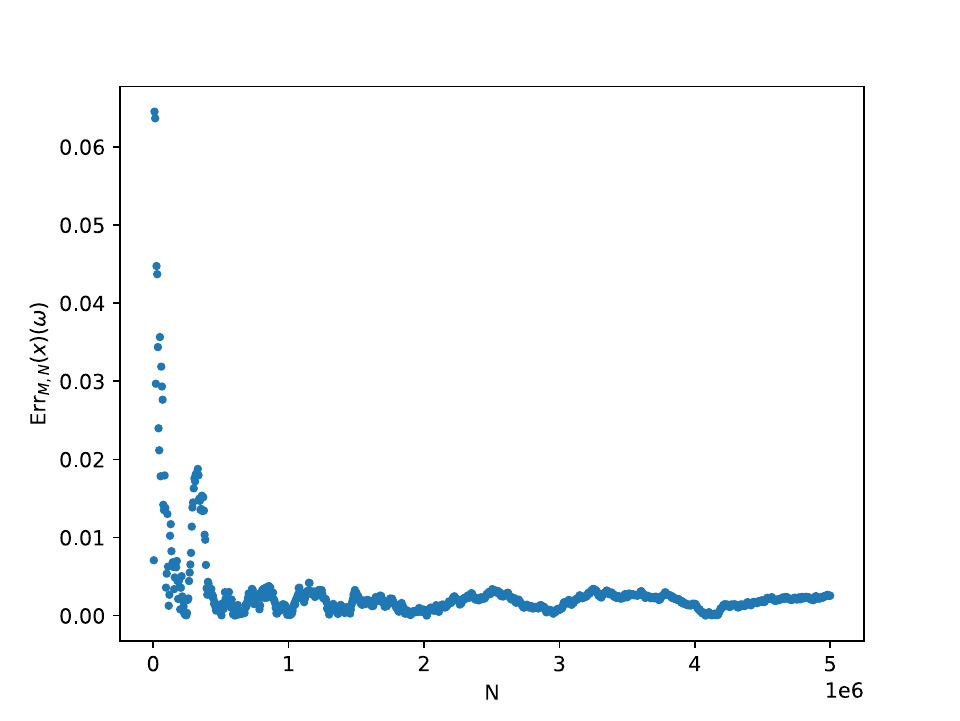}
\caption{The evolution of ${\sf Err}_{M,N}(x)(\omega)$ w.r.t. $N$}
\label{fig:Err_MN_omega_d=100}
\end{subfigure}
\begin{subfigure}{0.5\textwidth}
\centering
\includegraphics[scale=0.4]{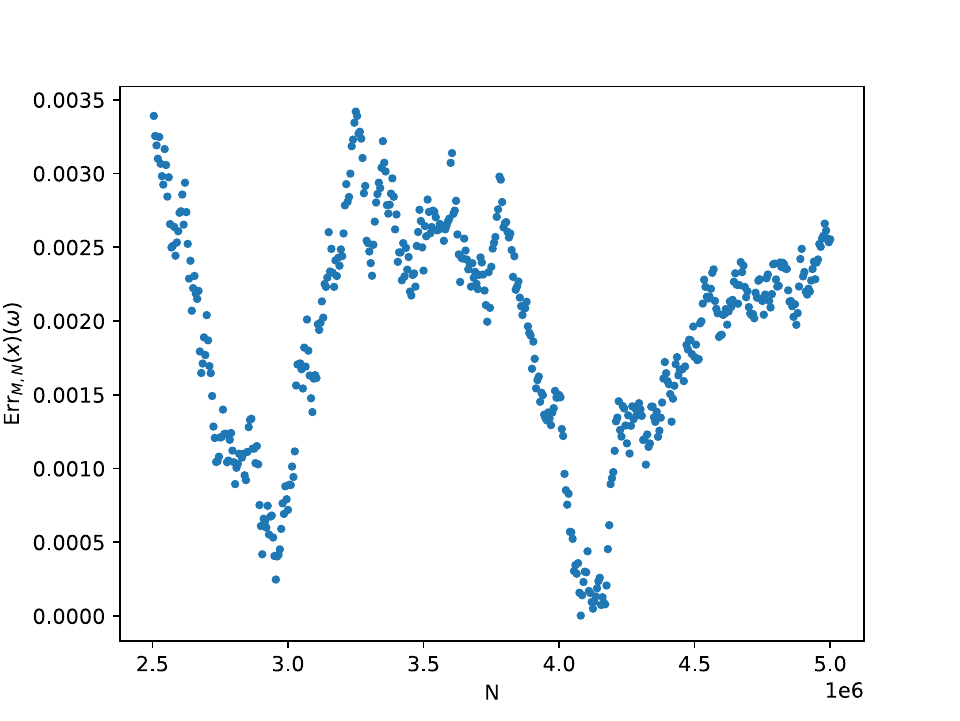}
\caption{Zoom in of \Cref{fig:Err_MN_omega_d=100}}
\label{fig:Err_MN_omega_zoom_d=100}
\end{subfigure}
\caption{ The evolution of ${\sf Err}_{M,N}(x)(\omega)$ w.r.t. $N$, for $D=D_{\sf ac}$, $d=100$, $M=2000$, and an arbitrarily fixed location $x\in D$, at which the exact solution $u(x)$ is $-0.357$.}
\label{fig:Err_MN_x_omega_d=100}
\end{figure}

{\bf Comments on Test 2.} 
\begin{enumerate}
    \item[$\bullet$] The histograms depicted in \Cref{fig:N_steps_d=10_ac} for the Poisson problem \eqref{pde:numerics} for $d=10$ and $D=D_{\sf ac}$ confirm that the random variables given by the normalized $L^1$-error $|{\sf Err}_{M,N}|_{L^1(D/|D|)}$ and the $L^\infty$-error $|{\sf Err}_{M,N,L}|_{L^\infty(D)}$ are small and concentrated if $M$ and $N$ are chosen sufficiently large, as stipulated by the theoretical results \Cref{thm:main} and \Cref{coro:noKepsilon}.
    \item[$\bullet$] The numerical tests depicted by \Cref{fig:N_steps_d=100_ac} for the Poisson problem \eqref{pde:numerics} for $d=100$ and $D=D_{\sf ac}$ confirm that the errors $\mathbb{E}\left\{|{\sf Err}_{M,N}|_{L^1(D/|D|)}\right\}$  and $\mathbb{E}\left\{\sup_{x\in D}|{\sf Err}_{M,N}(x)|\right\}$, approximated by \eqref{eq:test21} and \eqref{eq:test22} respectively, are decreasing to a small value as the number of WoS trajectories $N$ increases. The limit error attained when $N$ goes to infinity is not zero as it depends on $M$, but it decreases to zero as the latter parameter is also increased to infinity. This is discussed in the next comment.
    \item[$\bullet$] Furthermore, the results illustrated by \Cref{fig:Err_MN_x_omega_d=100} for $d=100$ show that the error ${\sf Err}_{M,N}(x)$ for an arbitrary chosen location $x\in D=D_{\sf ac}$ becomes much smaller than the errors obtained in \Cref{fig:N_steps_d=100_ac}, as the number of WoS steps $M$ is increased.
    \item[$\bullet$] Concerning the dependence of $M$ and $N$ with respect to $d$, our numerical tests revealed that on the one hand the choice of $M$ required in  \Cref{thm:main} or \Cref{coro:noKepsilon} is quite optimal, and on the other hand that the value of $N$ can be in fact taken much smaller than the one required in the same theoretical results. As a consequence suggested by this numerical evidence, one could expect that the width of the DNNs provided by \Cref{thm:mainNN} can be significantly reduced.
    \item[$\bullet$] Finally, let us emphasize that the numerical results obtained during {\bf Test 2} are conducted for the domain $D_{\sf ac}$ which is nor defective convex, neither satisfies the uniform exterior ball condition, hence the test turned out to be successful even on a worse domain geometry.
\end{enumerate}

}
\vspace{2mm}

\noindent \textbf{Acknowledgements.} Lucian Beznea and Oana Lupascu-Stamate were supported by a grant of the Ministry of Research, Innovation and Digitization, CNCS - UEFISCDI,
project number PN-III-P4-PCE-2021-0921, within PNCDI III.
Iulian Cimpean acknowledges support from the 
project  PN-III-P1-1.1-PD-2019-0780, within PNCDI~III. The work of Arghir Zarnescu has  been partially supported by the Basque Government through the BERC 2022-2025 program and by the Spanish State Research Agency through BCAM Severo Ochoa excellence accreditation Severo Ochoa CEX2021-00114 and through project PID2020-114189RB-I00 funded by Agencia Estatal de Investigaci\'on (PID2020-114189RB-I00 / AEI / 10.13039/501100011033).


\bibliographystyle{plain}
\bibliography{paper_revised1.bib}

@Inbook{Protter2005,
author="Protter, Philip E.",
title="Stochastic Differential Equations",
bookTitle="Stochastic Integration and Differential Equations",
year="2005",
publisher="Springer Berlin Heidelberg",
address="Berlin, Heidelberg",
pages="249--361",
isbn="978-3-662-10061-5",
doi="10.1007/978-3-662-10061-5_6",
url="https://doi.org/10.1007/978-3-662-10061-5_6"
}

@article{GoGrJeKo22,
author = {Gonon, Lukas and Grohs, Philipp and Jentzen, Arnulf and Kofler, David and Šiška, David},
title = "{Uniform error estimates for artificial neural network approximations for heat equations}",
journal = {IMA Journal of Numerical Analysis},
volume = {42},
number = {3},
pages = {1991-2054},
year = {2021}
}

@article{JuJeKrAnWu20,
  title={Overcoming the curse of dimensionality in the numerical approximation of semilinear parabolic partial differential equations},
  author={Hutzenthaler, Martin and Jentzen, Arnulf and Kruse, Thomas and Anh Nguyen, Tuan and von Wurstemberger, Philippe},
  journal={Proceedings of the Royal Society A},
  volume={476},
  number={2244},
  pages={20190630},
  year={2020},
  publisher={The Royal Society Publishing}
}

@article{GoGrJeKoSi22,
  title={Uniform error estimates for artificial neural network approximations for heat equations},
  author={Gonon, Lukas and Grohs, Philipp and Jentzen, Arnulf and Kofler, David and {\v{S}}i{\v{s}}ka, David},
  journal={IMA Journal of Numerical Analysis},
  volume={42},
  number={3},
  pages={1991--2054},
  year={2022},
  publisher={Oxford University Press}
}

@article{BeBeGrJe21,
author = {Beck, Christian and Becker, Sebastian and Grohs, Philipp and Jaafari, Nor and Jentzen, Arnulf},
title = {Solving the {K}olmogorov {P}{D}{E} by Means of {D}eep Learning},
year = {2021},
issue_date = {Sep 2021},
publisher = {Plenum Press},
address = {USA},
volume = {88},
number = {3},
issn = {0885-7474},
url = {https://doi.org/10.1007/s10915-021-01590-0},
doi = {10.1007/s10915-021-01590-0},
journal = {J. Sci. Comput.},
month = {sep},
numpages = {28}
}

@article{BeHuJeKu23,
title = {An overview on deep learning-based approximation methods for partial differential equations},
journal = {Discrete and Continuous Dynamical Systems - B},
volume = {28},number = {6},pages = {3697-3746},
year = {2023},
issn = {1531-3492},
doi = {10.3934/dcdsb.2022238},
url = {/article/id/639b367cb2114e413cb39c48},
author = {Christian Beck and Martin Hutzenthaler and Arnulf Jentzen and Benno Kuckuck}
}

@article{FaToWa11,
author = {Arash Fahim and Nizar Touzi and Xavier Warin},
title = {{A probabilistic numerical method for fully nonlinear parabolic PDEs}},
volume = {21},
journal = {The Annals of Applied Probability},
number = {4},
publisher = {Institute of Mathematical Statistics},
pages = {1322 -- 1364},
year = {2011},
doi = {10.1214/10-AAP723},
URL = {https://doi.org/10.1214/10-AAP723}
}

@article{SaTa95,
url = {https://doi.org/10.1515/mcma.1995.1.1.1},
title = {Integral Formulation of the Boundary Value Problems and the Method of Random Walk on Spheres},
author = {K.K. SABELFELD and D. TALAY},
pages = {1--34},
volume = {1},
number = {1},
journal = {Monte Carlo Methods and Applications},
doi = {doi:10.1515/mcma.1995.1.1.1},
year = {1995},
lastchecked = {2023-06-15}
}

@article{TaTu90,
author = { Denis   Talay  and  Luciano   Tubaro },
title = {Expansion of the global error for numerical schemes solving stochastic differential equations},
journal = {Stochastic Analysis and Applications},
volume = {8},
number = {4},
pages = {483-509},
year  = {1990},
publisher = {Taylor & Francis},
doi = {10.1080/07362999008809220},
}

@article{HuJeKu19,
author = {Hutzenthaler, Martin and Jentzen, Arnulf and Kruse, Thomas},
journal={ Foundations of Computational Mathematics },
year = {2022},
pages = {905–966},
title = {Overcoming the curse of dimensionality in the numerical approximation of parabolic partial differential equations with gradient-dependent nonlinearities}
}

@article{HuJeWu20,
author = {Martin Hutzenthaler and Arnulf Jentzen and von Wurstemberger Wurstemberger},
title = {{Overcoming the curse of dimensionality in the approximative pricing of financial derivatives with default risks}},
volume = {25},
journal = {Electronic Journal of Probability},
number = {none},
publisher = {Institute of Mathematical Statistics and Bernoulli Society},
pages = {1 -- 73},
year = {2020},
doi = {10.1214/20-EJP423},
URL = {https://doi.org/10.1214/20-EJP423}
}

@article{HuJeKrNg20,
author = {Hutzenthaler, Martin and Jentzen, Arnulf and Kruse, Thomas and Nguyen, Tuan},
year = {2020},
month = {04},
pages = {},
title = {A proof that rectified deep neural networks overcome the curse of dimensionality in the numerical approximation of semilinear heat equations},
volume = {1},
journal = {SN Partial Differential Equations and Applications},
doi = {10.1007/s42985-019-0006-9}
}

@inproceedings{Gall,
  title={Spatial Branching Processes, Random Snakes, and Partial Differential Equations},
  author={J. F. Le Gall},
  year={1999}
}

@article{beznea,
title = {Nonlinear PDEs and measure-valued branching type processes},
journal = {Journal of Mathematical Analysis and Applications},
volume = {384},
number = {1},
pages = {16-32},
year = {2011},
note = {Special Issue on Stochastic PDEs in Fluid Dynamics, Particle Physics and Statistical Mechanics},
issn = {0022-247X},
doi = {https://doi.org/10.1016/j.jmaa.2010.10.034},
url = {https://www.sciencedirect.com/science/article/pii/S0022247X10008711},
author = {Lucian Beznea and Andrei-George Oprina},
}

@article{bsde,
author = {Cheridito, Patrick and Soner, H. and Touzi, Nizar and Victoir, Nicolas},
year = {2005},
month = {10},
pages = {},
title = {Second order backward stochastic differential equations and fully non-linear parabolic PDEs}
}

@article{Baudoin,
author = {Fabrice Baudoin},
title = {{Stochastic analysis on sub-Riemannian manifolds with transverse symmetries}},
volume = {45},
journal = {The Annals of Probability},
number = {1},
publisher = {Institute of Mathematical Statistics},
pages = {56 -- 81},
year = {2017},
doi = {10.1214/14-AOP964},
URL = {https://doi.org/10.1214/14-AOP964}
}

@article{deaconu,
title = {The walk on moving spheres: A new tool for simulating Brownian motion’s exit time from a domain},
journal = {Mathematics and Computers in Simulation},
volume = {135},
pages = {28-38},
year = {2017},
note = {Special Issue: 9th IMACS Seminar on Monte Carlo Methods},
issn = {0378-4754},
doi = {https://doi.org/10.1016/j.matcom.2015.07.004},
url = {https://www.sciencedirect.com/science/article/pii/S0378475415001366},
author = {M. Deaconu and S. Herrmann and S. Maire},
}

@article{BaTa,
author = {Bally, Vlad and Talay, Denis},
year = {1994},
month = {07},
pages = {},
title = {The Law of the Euler Scheme for Stochastic Differential Equations: I. Convergence Rate of the Distribution Function},
volume = {104},
journal = {Probability Theory and Related Fields},
doi = {10.1007/BF01303802}
}

@article{BuChSy,
author = {Krzysztof Burdzy and Zhen-Qing Chen and John Sylvester},
title = {{The heat equation and reflected Brownian motion in time-dependent domains}},
volume = {32},
journal = {The Annals of Probability},
number = {1B},
publisher = {Institute of Mathematical Statistics},
pages = {775 -- 804},
year = {2004},
doi = {10.1214/aop/1079021464},
URL = {https://doi.org/10.1214/aop/1079021464}
}

@article{Hsu,
author = {Zhou, Yijing and Cai, Wei and Hsu, Elton},
year = {2015},
month = {02},
pages = {},
title = {Computation of Local Time of Reflecting Brownian Motion and Probabilistic Representation of the Neumann Problem},
volume = {15},
journal = {Communications in Mathematical Sciences},
doi = {10.4310/CMS.2017.v15.n1.a11}
}

@article{DeLe,
author = {Deaconu, Madalina and Lejay, Antoine},
year = {2006},
month = {03},
pages = {135-151},
title = {A Random Walk on Rectangles Algorithm},
volume = {8},
journal = {Methodology And Computing In Applied Probability},
doi = {10.1007/s11009-006-7292-3}
}

@article{thal,
  title={Exponential integrability and exit times of diffusions on sub-Riemannian and metric measure spaces},
  author={Thalmaier, Anton and Thompson, James},
  journal={Bernoulli},
  volume={26},
  number={3},
  pages={2202--2225},
  year={2020},
  publisher={Bernoulli Society for Mathematical Statistics and Probability}
}

@article{sabtal,
title = {Integral Formulation of the Boundary Value Problems and the Method of Random Walk on Spheres},
author = {K.K. Sabelfeld and D. Talay},
pages = {1--34},
volume = {1},
number = {1},
journal = {Monte Carlo Methods and Applications},
doi = {doi:10.1515/mcma.1995.1.1.1},
year = {1995},
lastchecked = {2022-09-07}
}

@article{arnli,
author = {Marc Arnaudon and Xue-Mei Li},
title = {{Reflected Brownian motion: selection, approximation and linearization}},
volume = {22},
journal = {Electronic Journal of Probability},
number = {none},
publisher = {Institute of Mathematical Statistics and Bernoulli Society},
pages = {1 -- 55},
year = {2017},
doi = {10.1214/17-EJP41},
URL = {https://doi.org/10.1214/17-EJP41}
}

@article{LeStTr,
author = {Haesung Lee, Wilhelm Stannat, Gerald Trutnau},
title = {{Analytic theory of Itô-stochastic differential equations with non-smooth coefficients}},
volume = {},
journal = {arxix},
number = {none},
publisher = {},
pages = {},
year = {2022},
doi = {},
URL = {https://arxiv.org/abs/2012.14410}
}

@ARTICLE{Max,
       author = {{Konarovskyi}, Vitalii and {Marx}, Victor and {von Renesse}, Max},
        title = "{Spectral gap estimates for Brownian motion on domains with sticky-reflecting boundary diffusion}",
      journal = {arXiv e-prints},
         year = 2021,
        month = may,
          eid = {arXiv:2106.00080},
        pages = {arXiv:2106.00080},
archivePrefix = {arXiv},
       eprint = {2106.00080},
 primaryClass = {math.PR},
       adsurl = {https://ui.adsabs.harvard.edu/abs/2021arXiv210600080K}
}

@article{Kuran,
  title={The convexity of a domain and the superharmonicity of the signed distance function},
  author={D.H. Armitage and {\"U}. Kuran},
  journal={Proceedings of the Amer. Math. Soc.},
  volume={93},
  number={4},
  pages={598--600},
  year={1985}
}

@article{Ge92,
  title={The probabilistic solution of the {D}irichlet problem for $1/2{\Delta} + \langle a, \nabla\rangle + b$ with singular coefficients},
  author={W.D. Gerhard},
  journal={Journal of Theoretical Probability},
  volume={5},
  number={3},
  pages={503--520},
  year={1992},
  publisher={Springer}
}

@book{GTbook,
  title={Elliptic Partial Differential Equations of Second Order},
  author={D. Gilbarg and N.S. Trudinger},
  year={2001},
  publisher={Springer}
}

@article{Mu56,
  title={Some continuous Monte Carlo methods for the Dirichlet problem},
  author={M.E. Muller},
  journal={The Annals of Mathematical Statistics},
  pages={569--589},
  year={1956},
  publisher={JSTOR}
}

@article{Ky17,
author = {A.E. Kyprianou and A.Osojnik and T. Shardlow},
year = {2017},
month = {09},
pages = {},
title = {Unbiased `walk-on-spheres' Monte Carlo methods for the fractional Laplacian},
volume = {38},
journal = {IMA Journal of Numerical Analysis},
doi = {10.1093/imanum/drx042}
}

@article{binder2012rate,
  title={The rate of convergence of the walk on spheres algorithm},
  author={I. Binder and M. Braverman},
  journal={Geometric and Functional Analysis},
  volume={22},
  number={3},
  pages={558--587},
  year={2012},
  publisher={Springer}
}

@article{yarotsky2017error,
  title={Error bounds for approximations with deep ReLU networks},
  author={D. Yarotsky},
  journal={Neural Networks},
  volume={94},
  pages={103--114},
  year={2017},
  publisher={Elsevier}
}

@article{cybenko1989approximation,
  title={Approximation by superpositions of a sigmoidal function},
  author={G. Cybenko},
  journal={Mathematics of Control, Signals and Systems},
  volume={2},
  number={4},
  pages={303--314},
  year={1989},
  publisher={Springer}
}

@article{EDGB_2019,
  title={Deep neural network approximation theory},
  author={D. Elbr{\"a}chter and D. Perekrestenko and P.  Grohs and H. B{\"o}lcskei},
  journal={IEEE Transactions on Information Theory},
  volume={67},
  number={5},
  pages={2581--2623},
  year={2021},
  publisher={IEEE}
}

@article{pinkus1999approximation,
  title={Approximation theory of the MLP model in neural networks},
  author={A. Pinkus},
  journal={Acta Numerica},
  volume={8},
  pages={143--195},
  year={1999},
  publisher={Cambridge Univ. Press}
}

@article{GrHe21,
    author = {P. Grohs and L. Herrmann},
    title = "{Deep neural network approximation for high-dimensional elliptic PDEs with boundary conditions}",
    journal = {IMA Journal of Numerical Analysis},
    volume = {42},
    number = {3},
    pages = {2055-2082},
    year = {2021},
    month = {05},
    issn = {0272-4979},
    doi = {10.1093/imanum/drab031},
    url = {https://doi.org/10.1093/imanum/drab031},
    eprint = {https://academic.oup.com/imajna/article-pdf/42/3/2055/44955748/drab031.pdf},
}

@inproceedings{NEURIPS2021_7edccc66,
 author = {T. Marwah and Z.C. Lipton  and A. Risteski},
 booktitle = {Advances in Neural Information Processing Systems},
 editor = {M. Ranzato and A. Beygelzimer and Y. Dauphin and P.S. Liang and J. Wortman Vaughan},
 pages = {15044--15055},
 publisher = {Curran Associates, Inc.},
 title = {Parametric Complexity Bounds for Approximating PDEs with Neural Networks},
 url = {https://proceedings.neurips.cc/paper/2021/file/7edccc661418aeb5761dbcdc06ad490c-Paper.pdf},
 volume = {34},
 year = {2021}
}

@article{valenzuela2022new,
  title={A new approach for the fractional Laplacian via deep neural networks},
  author={N. Valenzuela},
  journal={arXiv preprint arXiv:2205.05229},
  year={2022}
}

@article{lagaris1998artificial,
  title={Artificial neural networks for solving ordinary and partial differential equations},
  author={I.E. Lagaris and A.  Likas and D.I. Fotiadis},
  journal={IEEE transactions on neural networks},
  volume={9},
  number={5},
  pages={987--1000},
  year={1998},
  publisher={IEEE}
}

@article{lagaris2000neural,
  title={Neural-network methods for boundary value problems with irregular boundaries},
  author={I.E Lagaris and  A.C.  Likas and D.G. Papageorgiou},
  journal={IEEE Transactions on Neural Networks},
  volume={11},
  number={5},
  pages={1041--1049},
  year={2000},
  publisher={IEEE}
}

@article{malek2006numerical,
  title={Numerical solution for high order differential equations using a hybrid neural network—optimization method},
  author={A. Malek and R. Shekari Beidokhti},
  journal={Applied Mathematics and Computation},
  volume={183},
  number={1},
  pages={260--271},
  year={2006},
  publisher={Elsevier}
}

@article{han2018solving,
  title={Solving high-dimensional partial differential equations using deep learning},
  author={J. Han and  A. Jentzen and E. Weinan},
  journal={Proceedings of the National Academy of Sciences},
  volume={115},
  number={34},
  pages={8505--8510},
  year={2018},
  publisher={National Acad. Sciences}
}

@article{raissi2019physics,
  title={Physics-informed neural networks: A deep learning framework for solving forward and inverse problems involving nonlinear partial differential equations},
  author={M. Raissi and  P. Perdikaris and G. Karniadakis},
  journal={Journal of Computational Physics},
  volume={378},
  pages={686--707},
  year={2019},
  publisher={Elsevier}
}

@article{SIRIGNANO20181339,
title = {{D}{G}{M}: A deep learning algorithm for solving partial differential equations},
journal = {Journal of Computational Physics},
volume = {375},
pages = {1339-1364},
year = {2018},
issn = {0021-9991},
doi = {https://doi.org/10.1016/j.jcp.2018.08.029},
url = {https://www.sciencedirect.com/science/article/pii/S0021999118305527},
author = {J. Sirignano and K. Spiliopoulos}
}

@article{jentzen2018proof,
  title={A proof that deep artificial neural networks overcome the curse of dimensionality in the numerical approximation of Kolmogorov partial differential equations with constant diffusion and nonlinear drift coefficients},
  author={A. Jentzen and D. Salimova  and T. Welti},
  journal={Communications in Mathematical Sciences},
  volume={19},
  number={5},
  pages={1167–1205},
  year={2021}
}

@article{pardoux1999forward,
  title={Forward-backward stochastic differential equations and quasilinear parabolic PDEs},
  author={Pardoux, Etienne and Tang, Shanjian},
  journal={Probability Theory and Related Fields},
  volume={114},
  number={2},
  pages={123--150},
  year={1999},
  publisher={Springer}
}

@book{gobet2016monte,
  title={Monte-Carlo methods and stochastic processes: from linear to non-linear},
  author={Gobet, Emmanuel},
  year={2016},
  publisher={Chapman and Hall/CRC}
}

@book{Hsu2,
  title={Stochastic analysis on manifolds},
  author={Hsu, Elton P},
  number={38},
  year={2002},
  publisher={American Mathematical Soc.}
}

@book{Stroock,
  title={An introduction to the analysis of paths on a Riemannian manifold},
  author={Stroock, Daniel W},
  number={74},
  year={2000},
  publisher={American Mathematical Soc.}
}

@article{von2005transport,
  title={Transport inequalities, gradient estimates, entropy and Ricci curvature},
  author={von Renesse, Max-K and Sturm, Karl-Theodor},
  journal={Communications on pure and applied mathematics},
  volume={58},
  number={7},
  pages={923--940},
  year={2005},
  publisher={Wiley Online Library}
}

@article{talay2,
  title={Probabilistic interpretation and random walk on spheres algorithms for the Poisson-Boltzmann equation in molecular dynamics},
  author={Bossy, Mireille and Champagnat, Nicolas and Maire, Sylvain and Talay, Denis},
  journal={ESAIM: Mathematical Modelling and Numerical Analysis},
  volume={44},
  number={5},
  pages={997--1048},
  year={2010},
  publisher={EDP Sciences}
}

@article{talay3,
  title={One-dimensional parabolic diffraction equations: pointwise estimates and discretization of related stochastic differential equations with weighted local times},
  author={Martinez, Miguel and Talay, Denis},
  journal={Electronic Journal of Probability},
  volume={17},
  pages={1--30},
  year={2012},
  publisher={Institute of Mathematical Statistics and Bernoulli Society}
}

@article{lejay2013new,
  title={New Monte Carlo schemes for simulating diffusions in discontinuous media},
  author={Lejay, Antoine and Maire, Sylvain},
  journal={Journal of computational and applied mathematics},
  volume={245},
  pages={97--116},
  year={2013},
  publisher={Elsevier}
}
\end{document}